\documentclass[10pt]{article}
\usepackage[english,activeacute]{babel}
\usepackage{epsfig}

\usepackage[latin1]{inputenc}
\usepackage[T1]{fontenc}
\usepackage{stmaryrd}
\usepackage{amsmath,amsfonts,amssymb,mathrsfs,amsthm}
\usepackage{xcolor}
\usepackage{dsfont}
\usepackage{graphicx}
\usepackage[mathscr]{eucal}
\usepackage{makeidx}
\usepackage{verbatim}
\usepackage{graphics,graphicx}
\usepackage{textcomp}
\usepackage{float}
\usepackage[colorlinks=true]{hyperref}
\usepackage{txfonts} 
\usepackage{mathrsfs} 

\newtheorem{lemma}{Lemma}[section]
\newtheorem{theorem}[lemma]{Theorem}
\newtheorem{proposition}[lemma]{Proposition}
\newtheorem{corollary}[lemma]{Corollary}

\theoremstyle{definition}
\newtheorem{example}[lemma]{Example}
\newtheorem{definition}[lemma]{Definition}
\newtheorem{remark}[lemma]{Remark}

\makeatletter
\def\keywords{
    \vspace{1ex}
    \noindent
    \if@twocolumn
      \small{\bf  Keywords}\/---$\!$    \else
      \begin{center}\small\ {\bf Keywords}\end{center}\quotation\small
    \fi}
\def\endkeywords{\vspace{0.6em}\par\if@twocolumn\else\endquotation\fi
    \normalsize\rm}
\makeatother

\renewcommand{\O}{\ensuremath{\mathcal O}}

\renewcommand{\L}{\ensuremath{\mathcal L}}
\newcommand{\calN}{\ensuremath{\mathcal N}}

\DeclareMathOperator{\Lie}{Lie}

\DeclareMathOperator{\Int}{Int}

\newcommand{\mb}[1]{\ensuremath{\mathbb{#1}}}
\newcommand{\N}{{\mb{N}}}

\newcommand{\R}{{\mb{R}}}

\newcommand{\w}{\ensuremath{w}}
\newcommand{\W}{\ensuremath{W}}
\newcommand{\y}{\ensuremath{y}}
\newcommand{\z}{\ensuremath{z}}

\newcommand{\eps}{\varepsilon}

\newcommand{\M}{\ensuremath{\mathcal M}}
\newcommand{\T}{\ensuremath{\mathbb T}}
\newcommand{\D}{\ensuremath{\mathscr D}}

\let \Re \relax
\DeclareMathOperator{\Re}{Re}

\newcommand{\ovl}[1]{\overline{#1}}

\DeclareMathOperator{\supp}{supp}

\DeclareMathOperator{\diag}{diag}
\DeclareMathOperator{\dist}{dist}
\DeclareMathOperator{\vois}{Vois}

\DeclareMathOperator{\vect}{span}

\renewcommand{\d}{\ensuremath{\partial}}

\newcommand{\dsp}{\displaystyle}

\let \div \relax
\DeclareMathOperator{\div}{div}

\DeclareMathOperator{\length}{length}

\newcommand{\E}{\mathscr E}
\newcommand{\Z}{\mathbb Z}

\renewcommand{\H}{\ensuremath{\mathcal H}}

\def\e{{\varepsilon}}

\newcommand\bna{\begin{eqnarray*}} 
\newcommand\ena{\end{eqnarray*}}

\newcommand\bnan{\begin{eqnarray}} 
\newcommand\enan{\end{eqnarray}}

\newcommand\bnp{\begin{proof}} 
\newcommand\enp{\end{proof}}

\newcommand\bneq{\begin{eqnarray*}\left\lbrace \begin{array}{rcl}}
\newcommand\eneq{\end{array} \right.\end{eqnarray*}}
\newcommand\bneqn{\begin{eqnarray}\left\lbrace \begin{array}{rcl}}
\newcommand\eneqn{\end{array} \right.\end{eqnarray}}

\newcommand\nor[2]{\left\|#1\right\|_{#2}}
\newcommand\sgn{\textnormal{sgn}}

\numberwithin{equation}{section}

\newtheorem{hypo}{Assumption}[section]

\author{Camille Laurent\footnote{CNRS UMR 7598 and UPMC Univ Paris 06, Laboratoire Jacques-Louis Lions, F-75005, Paris, France, email: laurent@ann.jussieu.fr} and Matthieu L\'eautaud\footnote{\'Ecole Polytechnique, Centre de Math\'ematiques Laurent Schwartz UMR7640,  91128 Palaiseau cedex France. 
Most of this research was done when the second author was in CRM, CNRS UMI 3457, Université de Montréal, Case Postale 6128, Succursale Centre-Ville, Montréal (QC) Canada H3C 3J7 and Universit\'e Paris Diderot, IMJ-PRG, UMR 7586, B\^atiment Sophie Germain, 75205 Paris Cedex 13 France, email: matthieu.leautaud@polytechnique.edu.} 
}

\begin{document}
\title{Tunneling estimates and approximate controllability for hypoelliptic equations}
\maketitle
\textbf{Abstract.} 
This article is concerned with quantitative unique continuation estimates for equations involving a ``sum of squares'' operator $\L$ on a compact manifold $\M$ assuming: $(i)$ the Chow-Rashevski-H\"ormander condition ensuring the hypoellipticity of $\L$, and $(ii)$ the analyticity of $\M$ and the coefficients of $\L$.

The first result is the tunneling estimate $\|\varphi\|_{L^2(\omega)} \geq Ce^{- \lambda^{\frac{k}{2}}}$ for normalized eigenfunctions $\varphi$ of $\L$ from a nonempty open set $\omega\subset \M$, where $k$ is the hypoellipticity index of $\L$ and $\lambda$ the eigenvalue.

The main result is a stability estimate for solutions to the hypoelliptic wave equation $(\d_t^2+\L)u=0$: for $T>2 \sup_{x \in \M}(\dist(x,\omega))$ (here, $\dist$ is the sub-Riemannian distance), the observation of the solution on $(0,T)\times \omega$ determines the data. The constant involved in the estimate is $Ce^{c\Lambda^k}$ where $\Lambda$ is the typical frequency of the data.

We then prove the approximate controllability of the hypoelliptic heat equation $(\d_t+\L)v=\mathds{1}_\omega f$ in any time, with appropriate (exponential) cost, depending on $k$. In case $k=2$ (Grushin, Heisenberg...), we further show approximate controllability to trajectories with polynomial cost in large time.

We also explain how the analyticity assumption can be relaxed, and a boundary $\d \M$ can be added in some situations.

Most results turn out to be optimal on a family of Grushin-type operators.

The main proof relies on the general strategy developed by the authors in~\cite{LL:15}. 

\tableofcontents
\setcounter{tocdepth}{1}

\section{Introduction and main results}

\subsection{Introduction}

Let $\M$ be a smooth compact connected $d$-dimensional manifold without boundary. We denote $\mathcal{X}^{\infty}$ the space of smooth vector fields on $\M$ (with real coefficients), which we identify to derivations on $\M$. We assume $\M$ is endowed with a smooth positive density measure $ds $, so that we may integrate functions on $\M$\footnote{See e.g.~\cite[Chapter~16 p427]{Lee:book}: given a local chart $(U_\kappa, \kappa)$ of $\M$, we have $\int_{U_\kappa} u \ ds=  \int_{\kappa(U_\kappa)} u\circ \kappa^{-1}(y) \varphi^\kappa(y) dy$ for an appropriate smooth positive function $\varphi^\kappa$, and for any $u \in C^0_c(U_\kappa)$.}.
We may then define the space $L^2(\M) = L^2(\M, ds)$ of square integrable functions with respect to this measure.
For $X \in \mathcal{X}^{\infty}$, we define by $X^*$ its formal dual operator for the duality of $L^2(\M)$, that is\footnote{Note that in the local chart $(U_\kappa, \kappa)$ we have $X^\kappa = \sum_j a_j^\kappa (x) \d_j$, and thus $(X^\kappa)^* = \sum_j -a_j^\kappa (x) \d_j - \d_j a_j^\kappa - \frac{\d_j \varphi^\kappa}{\varphi^\kappa}a_j^\kappa$, which is a vector field (namely $-X$) plus a multiplication operator (namely $-\div_{ds}(X)$, see Remark~\ref{r:sub-riem-lap} below).}, 
 $$
 \int_\M X^*(u)(x) v(x)ds(x) =   \int_\M u(x) X(v)(x)ds(x), \quad \text{ for any } u, v \in C^\infty(\M) .
 $$

\medskip
Given $m\in \N$ and $m$ vector fields\footnote{The assumption $1\leq m\leq d$ is sometimes made in the references we use, but can always be removed.} $X_1,\cdots,X_m \in \mathcal{X}^{\infty}$, we are interested in properties of the following (non-positive) second order operator, associated to the $X_i$'s (namely the so-called type~I H\"ormander operator)\footnote{See Remark~\ref{r:sub-riem-lap} below for a discussion on general sub-Riemannian Laplacians.}
\bnan
\label{def:L}
\L =   \sum_{i=1}^m X_i^*X_i .
\enan
Note that this operator is formally symmetric nonnegative, when defined on functions in $C^{\infty}(\M)$, since we have
\bnan
\label{e:LXi}
( \L u,u)_{L^2(\M)}=  \sum_{i=1}^m\nor{X_i u}{L^2(\M)}^2.
\enan

Both from the geometric control and the operator theoretic points of view, it is in this context natural to consider iterated Lie brackets of the vector fields $X_i$. 
We refer for instance to the following classical article~\cite{Bellaiche} and textbooks~\cite{Montgomery:book,Jean:book,Rifford:book,BarilariAgrachevBoscainBook}.
\begin{definition}
For any family $\mathcal{F}$ of smooth vector fields on $\M$ and $\ell\in \N^*$, we define the subspaces $\Lie^\ell(\mathcal{F})$ of $\mathcal{X}^{\infty} $ by iteration as follows:
\begin{itemize}
\item $\Lie^1(\mathcal{F})$ is the space spanned by $\mathcal{F}$ in $\mathcal{X}^{\infty}$,
\item $\Lie^{\ell+1}(\mathcal{F})=\vect \left(\Lie^\ell(\mathcal{F})\cup \left\{[X,Y];X\in \mathcal{F},Y\in \Lie^\ell(\mathcal{F})\right\}\right)$.
\end{itemize}
For any point $x \in \M$, $\ell\in \N^*$, we denote $\Lie^\ell(\mathcal{F})(x)$ the set of all tangent vectors $X(x)$ with $X\in \Lie^\ell(\mathcal{F})$.
\end{definition}

We shall always assume that the family $(X_i)$ satisfies the Chow-Rashevski-H\"ormander condition (or is ``bracket generating'').
\begin{hypo}
\label{assumLiek}
There exists $\ell \geq 1$ so that for any $x\in \M$, $\Lie^\ell(X_1,\cdots,X_m) (x)=T_x\M$\footnote{Note that it is sufficient to assume that for all $x\in \M$, there is $\ell = \ell(x)\in \N$ such that this holds. The upper semi-continuity of $x \mapsto \ell(x)$ and the compactness of $\M$ then imply the stronger form of Assumption~\ref{assumLiek} as stated.}. Denote then by $k \in \N^*$ the minimal $\ell$ for which this holds.
\end{hypo}
The integer $k$ will sometimes be refered to as the \textit{hypoellipticity index} of $\L$.
Assumption~\ref{assumLiek} is central in control theory and operator theory, for it characterizes both the controllability of the controlled ODE driven by the vector fields $(X_i)$ and the Hypoellipticity of the operator $\L$. Let us now recall these two seminal results, namely the Chow-Rashevski theorem and the H\"ormander theorem, which we both use in the sequel. 
\begin{theorem}[Chow~\cite{Chow:39}, Rashevski~\cite{Rash:38}]
\label{t:chow}
Under Assumption \ref{assumLiek}, the following statement holds: for any $x_0, x_1 \in \M$, any $T>0$, there exist $u_i \in L^1(0,T)$ for $i \in \{1, \cdots, m\}$ such that the unique solution of 
\bnan
\label{e:chow-rach}
\dot{\gamma}(t) = \sum_{i=1}^m u_i(t) X_i(\gamma(t)) , \quad \gamma(0) = x_0
\enan
satisfies $\gamma(T)=x_1$.
\end{theorem}
We refer e.g. to~\cite[Chapter~1.4]{Jean:book}, \cite[Chapter~1.4]{Rifford:book} or~\cite[Chapter~2]{Montgomery:book} for statements and proofs of the Chow-Rashevski theorem, and in particular for the definition of the solution of~\eqref{e:chow-rach}. See also~\cite[Chapter~3]{Cor:book} for examples and applications in control theory.
This theorem motivates the following definition.
\begin{definition}[Horizontal path]
\label{d:horizontal-path}
We say that an absolutely continuous function $\gamma: [0,T] \to \M$ is a horizontal path if there exist $u_i \in L^1(0,T; \R)$ for $i = 1, \cdots m$ such that for almost every $t \in [0,T]$, we have $\dot{\gamma}(t) = \sum_{i=1}^m u_i(t) X_i(\gamma(t))$.
\end{definition}
Such a trajectory is in particular absolutely continuous and almost everywhere tangent to the so-called horizontal distribution $\vect(X_1, \cdots , X_m )$. 
The second key role played by Assumption~\ref{assumLiek} in analysis is summarized in the following result.
\begin{theorem}[H\"ormander \cite{Hor67}, Rothschild-Stein \cite{RS:76}]
\label{thmhypoestim}
Under Assumption \ref{assumLiek}, the operator $\L$ in~\eqref{def:L} is hypoelliptic, that is, for all $u \in \D'(\M)$ and $x_0\in \M$, if $\L u \in C^\infty$ near $x_0$ then $u \in C^\infty$ near $x_0$.

Moreover, it is subelliptic of order $\frac{1}{k}$, that is, the following estimates hold: there is $C>0$ such that for any $u\in C^{\infty}(\M)$, we have
\bnan
\nor{u}{H^{\frac{1}{k}}(\M)}^2\leq C\sum_{i=1}^m\nor{X_i u}{L^2(\M)}^2+C\nor{u}{L^2(\M)}^2 , \label{estimhypo} \\
\nor{u}{H^{\frac{1}{k}}(\M)}^2\leq C(\L u,u)_{L^2(\M)}+C\nor{u}{L^2(\M)}^2 , \label{estimhypo2} \\
\nor{u}{H^{\frac{2}{k}}(\M)}^2\leq C\nor{\L u}{L^2(\M)}^2+C\nor{u}{L^2(\M)}^2 \label{estimhypo3}.
\enan
\end{theorem}
The hypoellipticity was shown by H\"ormander \cite{Hor67}, who also provides with a subelliptic estimate with loss (see also~\cite{Kohn:78,Kohn:05} or~\cite[Chapter~2]{HN:book} for a simpler proof). 
The optimal subelliptic estimate~\eqref{estimhypo} with gain of $1/k$ derivatives is proved by~\cite{RS:76} (see also~\cite[p288]{FP:81} for a different proof). More precisely, (even slightly hidden) it is written in \cite{RS:76} Theorem~17 and estimate (17.20) p311 in a local form. It is then easy to globalise on the compact manifold $\M$ to obtain~\eqref{estimhypo} (since commutators of $X_i$ with a smooth cutoff function is a multiplication operator).  

Both estimates~\eqref{estimhypo2} and~\eqref{estimhypo3} may then be deduced from~\eqref{estimhypo}. This is clear for~\eqref{estimhypo2} when recalling~\eqref{e:LXi}.
The proof of~\eqref{estimhypo3} requires a commutator argument (detailed e.g. in~\cite{FP:83}) and is proved in Appendix~\ref{app:Hs-comm}, as well as $H^s$ variants of~\eqref{estimhypo2} and~\eqref{estimhypo3}.
Note that these subelliptic estimates are also obtained in Fefferman-Phong \cite{FP:83} for some wider class of symetric operators, not neccessarily sums of squares, and with a shorter proof.

\bigskip
Since the operator $\L$ is symmetric non-negative, the hypoellipticity of $\L+1$ and the compactness of $\M$ directly imply that $\L$ is essentially selfadjoint (see e.g. Reed-Simon \cite[Theorem X.26]{ReedSimon}). 
Hence, it extends uniquely as a selfadjoint operator (its Friedrich extension)
$$
\L : D(\L) \subset L^2(\M) \to L^2(\M), 
$$
with, according to~\eqref{estimhypo3}, $H^2(\M) \subset D(\L) \subset H^{\frac{2}{k}}(\M)$ (still under Assumption \ref{assumLiek}). 
The operator $\L$ is hence selfadjoint on $L^2(\M)$, with compact resolvent: it admits a Hilbert basis of eigenfunctions $(\varphi_j)_{j \in \N}$, associated with the real eigenvalues $(\lambda_j)_{j \in \N}$, sorted increasingly, that is
\bnan
\label{e:spectral-elts}
\L\varphi_i=\lambda_i\varphi_i, \qquad (\varphi_i , \varphi_j)_{L^2(\M)} = \delta_{ij} , \qquad 0= \lambda_0 < \lambda_1 \leq \lambda_2 \leq \cdots \leq \lambda_j \to + \infty .
\enan
Note that a bootstrap argument in~\eqref{estimhypo3} shows that $\varphi_j \in C^\infty(\M)$.
In particular, the spectral decomposition allows to define solutions of the hypoelliptic wave and heat equations (respectively $(\d_t^2 +\L)v=0$ and $(\d_t +\L)u=0$), which we shall consider in this paper.  

\bigskip

In addition to Assumption~\ref{assumLiek}, we will also assume in the main part of the article that everything is real-analytic. This assumption in not made in Section~\ref{s:non-analytic} though, where we give some results in the non-analytic context.
\begin{hypo}
\label{hypoanal}
The manifold $\M$, the density $ds$, and the vector fields $X_i$ are \textbf{real-analytic}. 
\end{hypo}
In particular, it implies that the operator $\L$ has analytic coefficients in any analytic coordinate set compatible with the manifold $\M$.
Note that under this assumption, the converse of Theorems~\ref{t:chow} and~\ref{thmhypoestim} also hold, namely:
\begin{itemize}
\item Attainability, in the sense of Theorems~\ref{t:chow}, implies Assumption \ref{assumLiek}, see~\cite{Her:63,Nag:66} (see also~\cite{Lobry:70,Sus:73} for generalizations);
\item The hypoellipticity of $\L$, in the sense of Theorem~\ref{thmhypoestim}, implies Assumption \ref{assumLiek}, see~\cite[Theorem~2.2]{Derridj:71} if there is no point $x_0$ where all $X_i$ cancel.
\end{itemize}
The analyticity assumption is further discussed in Sections~\ref{s:without-analyticity} and~\ref{s:comp-to-other} below.

\bigskip
Before stating our main results, let us provide with classical examples of operators that are considered in the present paper.
\begin{example}[Elliptic operators, $k=1$]
\label{ex:elliptic}
In the case $k=1$, then, $\vect(X_1,\cdots,X_m) (x)= T_x \M$ for all $x \in \M$, and the operator $\L$ is elliptic. Most of the results stated in this paper (or stronger versions of them) are already known in this situation (and in greater generality), see~\cite{LR:95,Leb:Analytic,LL:15}. That all Laplace-Beltrami operators can be written under the form~\eqref{def:L} is a consequence of Remark~\ref{assumLiek} below.
\end{example}
\begin{example}[The Grushin operator, $k=2$]
\label{ex:Grushin}
 Consider the torus $\M = (\R/2 \Z) \times  (\R/\Z) $ (which we identify with $[-1, 1[\times [0,1[$ with periodicity conditions), endowed with the Lebesgue measure  $ds=dx_1 dx_2$ and 
$$
\L= - \big( \partial_{x_1}^2 + x_1^{2}\partial_{x_2}^2 \big) = X_1^*X_1+X_2^* X_2, \quad \text{ with } \quad  X_1=\frac{\partial}{\partial x_1} , \quad X_2=x_1 \frac{\partial}{\partial x_2}.
$$
 We have $\vect(X_1,X_2)=\R^2$ if $x_1 \neq 0$, but  on the singular set $x_1 = 0$, we have $\vect(X_1,X_2)= \R X_1$. However, we have $[X_1 ,X_2] = \frac{\partial}{\partial x_2}$, so that $\vect(X_1 ,X_2 , [X_1 ,X_2])= \R^2$ on the whole $\M$, and Assumption~\ref{assumLiek} is satisfied for $k=2$.
Remark that $x_1^{2}$ is not analytic (not even $C^1$) on the torus $ \M$; here it can be replaced e.g. by $\sin(\pi x_1/2)^2$, being analytic and satisfying the same Hypoelliptic property. The original Grushin example will also be discussed later with Dirichlet boundary conditions, in which case it is smooth on $[-1, 1]\times (\R/\Z)$ or $[-1, 1]\times [0,1]$.
\end{example}
\begin{example}[Higher order Grushin operators, $k\in \N$]
\label{ex:Grushin+}
Consider again $\M = (\R/2 \Z) \times  (\R/\Z)$,  $ds=dx_1 dx_2$ and, for $\gamma \in \N$, set
\bnan
\label{e:Grushin+}
\L_\gamma= - \big( \partial_{x_1}^2 + x_1^{2\gamma}\partial_{x_2}^2 \big) = X_1^*X_1+X_2^* X_2, \quad \text{ with } \quad  X_1=\frac{\partial}{\partial x_1} , \quad X_2=x_1^\gamma \frac{\partial}{\partial x_2}.
\enan
Again, $x_1^{2\gamma}$ may be replaced by $\sin(\pi x_1/2)^{2\gamma}$ so that $\L_\gamma$ has analytic coefficients.
 We have $\vect(X_1,X_2)=\R^2$ if $x_1 \neq 0$, but  on the singular set $x_1 = 0$, we have to use iterated Lie brackets: Since $\left[\frac{\partial}{\partial x_1},x_1^{\beta}\frac{\partial}{\partial x_2}\right]= \beta x_1^{\beta-1}\frac{\partial}{\partial x_2}$ for all $\beta\geq 1$, we have, with $\mathcal{F}=\left\{X_1,X_2\right\}$, that 
\begin{itemize}
\item $\Lie^1(\mathcal{F})$ is the space spanned by $\mathcal{F}$ in $\mathcal{X}^{\infty}$;
\item $\Lie^{\ell}(\mathcal{F})=\left\{f=a\frac{\partial}{\partial x_1}+\sum_{i=\gamma-\ell-1}^{\gamma}b_ix_1^{i}\frac{\partial}{\partial x_2}\left|a,b_i\in \R\right.\right\}$ for $1\leq \ell \leq \gamma+1$;
\item $\Lie^{\ell}(\mathcal{F})=\Lie^{\gamma+1}(\mathcal{F}) = \left\{f=a\frac{\partial}{\partial x_1}+B(x_1)\frac{\partial}{\partial x_2}\left|a\in \R,B\in \R^{\gamma}[X]\right.\right\}$ if $\ell\geq \gamma+1$.
\end{itemize}
Hence, for $x=(0, x_2)$, we have $\Lie^{\ell}(\mathcal{F})(x)=\R \frac{\partial}{\partial x_1}$ if $\ell<\gamma+1$ and $\Lie^{\gamma+1}(\mathcal{F})(x)=\R^2$.
In particular, Assumption \ref{assumLiek} is fulfilled with $k=\gamma+1$. Note that we recover Example~\ref{ex:elliptic} in case $\gamma=0$ and Example~\ref{ex:Grushin} in case $\gamma=1$. 
\end{example}

\begin{example}[The Heisenberg operator on the Heisenberg group]
\label{exHeisen}
On $\R^3$ with current point $w=(x,y,s)$, the following two vector fields $X_1 = \d_{x}+2y\d_{s}$ and $X_2 = \d_{y} - 2x \d_{s}$ constitute the model case for contact geometry. Indeed, we have, with $\mathcal{F}=\left\{X_1,X_2\right\}$, that 
\begin{itemize}
\item $\Lie^1(\mathcal{F}) = \vect(\mathcal{F})$ is of dimension $2$ at any point in $\R^3$;
\item $\Lie^{2}(\mathcal{F})=\R^3$ at any point in $\R^3$, since $[X_1,X_2] = -4\d_{s}$.
\end{itemize}
Let us now define a compact context in which these are two analytic vector fields. 

First equip $\R^3$ with the (non-commutative) group law
$$
w\bullet w'=(x,y,s)\bullet (x',y',s') = (x+x',y+y',s+s'-2xy'+2yx') .
$$
With this law, $(\R^3, \bullet)$ (with $\R^3$ endowed with its canonical differential structure) is a Lie group which we denote by $G$. Given $L>0$, the set $\Gamma = L\Z \times L\Z \times L^2\Z$ is a subgroup of $G$, and both vector fields $X_1$ and $X_2$ are left invariant vector fields on $G$, i.e. setting $m_g :G\to G, w \mapsto g\bullet w$, we have $dm_g(X_j(w)) = X_j(m_g(w)) = X_j(g\bullet w)$ for $j=1,2$. The subgroup $\Gamma$ being co-compact, the left quotient $\M :=\Gamma \setminus G$ is a compact three dimensional analytic manifold. Moreover, the vector fields $X_1, X_2$ go to the quotient as analytic vector fields on $\M$. From the computation on $\R^3$, we obtain $\dim \Lie^{1}(\mathcal{F})(w)=2$ and $\Lie^{2}(\mathcal{F})(w)=T_w\M$ for some/any point $w\in \M$. The Haar measure turns out to be the Euclidian measure in the coordinates $(x,y,s)$. We consider the operator $\L = X_1^*X_1+X_2^* X_2=-X_1^2-X_2^2=-\Delta_{\mathbb{H}}$, where $\Delta_{\mathbb{H}}$ is the Kohn Laplacian, for which $k=2$. We refer for instance to~\cite[Section~1.2]{BFKG:12} for more on this example.
\end{example}

This last example belongs to the following general class of constant rank sub-Riemannian structures.
\begin{example}[Lie Groups]
\label{ex:lie-groups}
Assume that $(\M,\bullet)$ is a {\em compact} $d$-dimensional Lie group. Let $\mathds{1}$ be the identity of $(\M,\bullet)$, and write $L:= T_{\mathds{1}}\M$ its Lie algebra. Recall (see e.g.~\cite[Tome~II, p627]{Gode:82}) that there is a unique real-analytic differentiable structure on $\M$ compatible with the action of $\bullet$, with which we endow $\M$.
We write as in the above example $m_g : \M \to \M, x \mapsto g\bullet x$ for the left multiplication. 
Given $m<d$ and $m$ vectors $(e_1, \cdots e_m) \in L^m$, we denote by $(X_1, \cdots, X_m)$ the associated $m$ left-invariant vector fields defined, for $x \in \M$, by $X_j(x) := dm_x(X_j(\mathds{1})) = dm_x(e_j)$.

Now, we assume that the vectors $(e_1, \cdots e_m)$ generate the whole Lie algebra, namely $\Lie(\M)=L$, which implies that the vector fields $(X_1, \cdots, X_m)$ satisfy Assumption~\ref{assumLiek}, for some $k$.

Finally, we remark that, by construction, both the vector fields $X_j$ and the Haar measure $ds$ of $\M$ are real-analytic and left invariant. All our results shall hence apply to the associated operator $\L$.
\end{example}
Finally, let us mention that hypoelliptic operators appear naturally in several physical and mathematical contexts such as stochastic processes and the theory of functions of several complex variables. We refer to \cite[Chapter 2]{Bramanti:14} for a presentation of some of these applications. 

\subsection{Main results}
\label{s:results} 
Our main results under Assumptions~\ref{assumLiek} and~\ref{hypoanal} are of three different types:
\begin{enumerate}
\item Tunneling estimates for eigenfunctions $\varphi_j$ of $\L$ (Section~\ref{s:tunneling}); 
\item Quantitative approximate observability (and associated controllability) of the hypoelliptic wave equation $(\d_t^2 +\L)v=0$ from a subset $\omega \subset \M$ (Section~\ref{s:approx-control-wave});
\item Quantitative approximate observability (and associated controllability) of the hypoelliptic heat equation $(\d_t +\L)u=0$ from $\omega$ (Section~\ref{s:approx-control-heat});.
\end{enumerate}
Also, we provide with a class of examples (which are generalizations of those considered in Example~\eqref{ex:Grushin+}) where all these results hold as well without the analyticity Assumption~\ref{hypoanal} (Section~\ref{s:without-analyticity}). 

All results obtained depend explicitely on the hypoellipticity index $k$ of the operator considered, i.e. the minimal number of iterated brackets necessary to span the whole tangent space, given by Assumption~\ref{assumLiek}.
We finally prove with an example that the results are optimal in general.

\subsubsection{Eigenfunction tunneling}
\label{s:tunneling}
Our first result is the following.
\begin{theorem}
\label{t:spec-ineq}
Let $\omega$ be a nonempty open subset of $\M$. Then, there is $C,c>0$ such that every eigenfunction $\varphi_i$ of $\L$ associated to the eigenvalue $\lambda_i$ satisfies
\bnan
\label{e:eigenfunction-tunneling}
\| \varphi_j\|_{L^2(\M)}\leq Ce^{c\lambda_j^{k/2}}\|\varphi_j\|_{L^2(\omega)} .
\enan 
\end{theorem}
This estimate may be read as $\|\varphi_j\|_{L^2(\omega)} \geq \frac{1}{C} e^{-c\lambda_j^{k/2}}$ for all normalized eigenfunctions, and hence quantizes the possible vanishing rate of eigenfunctions on any subdomain $\omega$. 

In the case $k=1$, i.e. when $\L$ is an elliptic operator, the analyticity assumption~\ref{hypoanal} is not needed and the result follows from the Donnelly-Fefferman paper~\cite{DF:88}. In this situation, it also holds on a manifold with boundary for Dirichlet eigenfunctions~\cite{LR:95} (see also~\cite{LR:97} for other boundary conditions). 

\bigskip
We shall also deduce from estimates of \cite[Section~2.3]{BeauchardCanGugl} that the tunneling estimate~\eqref{e:eigenfunction-tunneling} is optimal in the following particular setting (close to Example~\eqref{ex:Grushin+}).
\begin{example}[Higher order Grushin operators on the rectangle]
\label{ex:Grushin++}
Consider the manifold {\em with boundary} $\M=[-1,1] \times [0,1]$ or $\M=[-1,1]\times(\R/\Z)$, endowed with the Lebesgue measure $dx$, and for $\gamma>0$, define the operator $\L_\gamma =  - \big( \partial_{x_1}^2 + x_1^{2\gamma}\partial_{x_2}^2 \big)$ as in~\eqref{e:Grushin+} with Dirichlet conditions on $\d \M$. If  $\gamma \in \N$, then the operator $\L_\gamma$ is hypoelliptic of order $k=\gamma+1$ (i.e. Assumption~\ref{assumLiek} is fulfilled with $k=\gamma+1$).
\end{example}

\begin{proposition}
\label{Prop:BCG}
Consider, for $\gamma>0$ the situation of Example~\ref{ex:Grushin++}.
Assume that $\overline{\omega} \cap \left\{x_1=0\right\} = \emptyset$. Then there exists $C,c_0>0$ and a sequence  $(\lambda_j , \varphi_j)$ of eigenvalues and associated eigenfunctions of $\L_{\gamma}$ with $\lambda_j \to +\infty$ such that 
\bna
\| \varphi_j\|_{L^2(\omega)} \leq  Ce^{ - c_0 \lambda_j^{\frac{\gamma+1}{2}}}\|\varphi_j\|_{L^2(\M)} .
\ena 
\end{proposition}
We recall that if $\gamma \in \N^*$, then $\L_\gamma$ is hypoelliptic of order $k=\gamma+1$, so that Proposition~\ref{Prop:BCG} shows that, in general, one cannot expect a better estimate than that of Theorem~\ref{t:spec-ineq}. We shall also prove (see Section~\ref{s:without-analyticity}) that Estimate~\eqref{e:eigenfunction-tunneling} holds as well in a setting containing those of Example~\ref{ex:Grushin+} and Example~\ref{ex:Grushin++}, thus providing a genuine converse of Proposition~\ref{Prop:BCG} (for $\gamma \in \N^*$).

Note that in the analytic context, the {\em qualitative} uniqueness: 
\bna
\Big( \L \varphi  = \lambda \varphi \text{ on }\M , \quad \varphi = 0 \text{ on } \omega \Big) \implies \varphi \equiv 0\text{ on }\M ,
\ena
was proved by Bony~\cite{Bo:69}, as a consequence of the Holmgren-John theorem. Removing the analyticity assumption, even for such a qualitative unique continuation property, remains a very subtle issue, as discussed in Section~\ref{s:comp-UC} below.

\subsubsection{Quantitative approximate observability of the hypoelliptic wave equation}
\label{s:approx-control-wave}
To state our main result here, we need to introduce the appropriate notions of Sobolev spaces and sub-Riemannian distance, which are adapted to the analysis of the operator $\L$.

All along the paper, we shall use the functional calculus given, for appropriate functions $f$ and $u$, by
\bnan
\label{e:fct-calcul}
f(\L) u = \sum_{j \in \N} f(\lambda_j) (u , \varphi_j)_{L^2(\M)} \varphi_j . 
\enan
This allows for instance to define the operators $(1+\L)^\frac{s}{2} : C^\infty(\M) \to C^\infty(\M)$, which, by duality, may be extended as operators $(1+\L)^\frac{s}{2}   : \D'(\M) \to \D'(\M)$.
We next define the Sobolev spaces 
$$
 \H^s_\L =  \{ u \in \D'(\M) ,  \ (1+\L)^\frac{s}{2}   u \in L^2(\M) \} , \quad  s\in \R ,
$$
and associated norms
\bna
\nor{u}{\H^s_\L}=\nor{(1+\L)^\frac{s}{2} u}{L^2(\M)} , \quad  s\in \R .
\ena
 
Let us now also introduce basic notions of sub-Riemannian geometry needed to formulate our main result. We refer to~\cite{Bellaiche,Montgomery:book,Jean:book,Rifford:book,ABB:EMS16,BarilariAgrachevBoscainBook} for an introductions to sub-Riemannian geometry, as well as for further developments.
The so-called sub-Riemannian metric associated to $(X_1, \cdots ,X_m)$ is defined, for $x\in \M$ and $v\in T_x\M$, by 
\bnan
\label{def-g}
g(x,v):=
\left\{
\begin{array}{ll} \dsp
\inf\left\{\sum_{i=1}^m u_i^2\left|(u_1,\cdots,u_m)\in \R^m,\sum_{i=1}^mu_i X_i(x)=v\right.\right\}& \text{if } v \in \vect(X_i(x), i \in \{1, \cdots, m\} ), \\
+ \infty & \text{if not} .
\end{array}
\right.
\enan
This defines for any $x \in \M$ a positive definite quadratic form $g(x, \cdot)$ on the the horizontal space $\vect(X_1(x), \cdots, X_m(x))$. Remark that, if finite, the infimum is in fact a minimum, and is realized by a unique vector $(u_1,\cdots ,u_m)\in \R^m$.
Given $\gamma : [0,1] \to \M$ an absolutely continuous path, we define its length accordingly by
\bna
\length(\gamma):=\int_0^1 \sqrt{g(\gamma(t),\dot{\gamma}(t))}dt.
\ena
The fact that this quantity is finite implies that $\dot{\gamma}(t) \in \vect \big(X_1(\gamma(t)), \cdots, X_m(\gamma(t)) \big)$ for almost all $t \in [0,1]$. Also, it is always finite if $\gamma$ is a horizontal path (in the sense of Definition~\ref{d:horizontal-path}).

Then, this allows to define a sub-Riemannian (also called Carnot-Carath\'eodory) distance on $\M$ by
\bna
d_{\L}(x_0,x_1)=\inf\left\{\length(\gamma)\left|\gamma \textnormal{ horizontal path, }\gamma(0)=x_0,\gamma(1)=x_1 \right.\right\}.
\ena
The Chow-Rashevski Theorem~\ref{t:chow} implies that, under Assumption \ref{assumLiek}, the distance defined by $d_\L$ is always finite.

With these definitions in hand, we may now state our main result, which concerns the quantitative unique continuation (or quantitative approximate observability) for the Hypoelliptic wave equation
\bneqn
\label{hypoelliptic-wave}
\partial_t^2 u+ \L u&=&0\\
(u,\partial_t u)|_{t=0}&=&(u_0,u_1) .
\eneqn
 
\begin{theorem}
\label{thmwavehypo}
Let $\L$ as above satisfying Assumptions \ref{assumLiek} and \ref{hypoanal}.
Assume that $\omega$ is a non empty open set of $\M$ and let $T>\sup_{x\in \M} d_\L(x,\omega)$.
Then, there exist $\kappa , C ,\mu_0 >0$ such that we have 
\bnan
\label{th-estimate-k}
\nor{(u_0,u_1)}{L^2\times \H^{-1}_\L} \leq C e^{\kappa \mu^k}\nor{u}{L^2(]-T,T[\times \omega)} +\frac{1}{\mu}\nor{(u_0,u_1)}{\H^{1}_\L\times L^2}
\enan
for all $\mu\geq \mu_0$ and for any $(u_0,u_1)\in \H^{1}_\L\times L^2$, and associated solution $u$ solution of~\eqref{hypoelliptic-wave} on $]-T,T[$.
\end{theorem}
Note first that this estimate could be stated equivalently for all $\mu >0$ (see e.g.~\cite[Lemma~A.3]{LL:15}). We chose to keep the above formulation to underline the interesting case (being $\mu$ large).
Note also that this theorem can be equivalently rewritten under one of the following two formulations (see e.g.~\cite[Lemma~A.3]{LL:15})
\bnan
\label{th-estimate-Lambda}
\nor{(u_0,u_1)}{\H^{1}_\L\times L^2}\leq Ce^{c\Lambda^k} \nor{u}{L^2(]-T,T[\times \omega)},
\quad \text{ with } \Lambda=\frac{\nor{(u_0,u_1)}{\H^{1}_\L\times L^2}}{\nor{(u_0,u_1)}{L^2\times \H^{-1}_\L}} ,
\enan
or 
\bnan
\label{th-estimate-log}
\nor{(u_0,u_1)}{L^2\times \H^{-1}_\L}\leq C \frac{\nor{(u_0,u_1)}{\H^1_\L \times L^2}}{\log\left(\frac{\nor{(u_0,u_1)}{\H^1_\L \times L^2}}{\nor{u}{L^2(]-T,T[\times \omega)}}+1\right)^{\frac{1}{k}}} ,
\enan
where, in the last expression, the function $x \mapsto \left(\log(1+ \frac{1}{x})\right)^{-\frac{1}{k}}$ has to be extended by zero at $x= 0^+$.

Again, in the particular situation of Example~\ref{ex:Grushin++}, i.e. for the operators~\eqref{e:Grushin+}, the sequence of eigenfunctions of Proposition~\ref{Prop:BCG} shows that the exponent $e^{\kappa \mu^k}$ in~\eqref{th-estimate-k} (resp. $e^{c\Lambda^k}$ in \eqref{th-estimate-Lambda} and $\log^{-\frac{1}{k}}$ in~\eqref{th-estimate-log}) cannot be improved in general.

\begin{remark}
That $d_\L$ is the relevant distance function in view of the study of the Hypoelliptic partial differential operator $\L$ comes from the fact that the sub-Riemannian metric $g(x,v)$ and the principal symbol $\ell (x, \xi)$ of the operator $\L$ are linked through the Legendre transform $\frac12 g(x,v) = \max_{\xi \in T^*_x \M}(\left< \xi , v\right> - \frac12 \ell(x, \xi ))$ (see~\cite[Section~1.2]{Bellaiche} or Appendix~\ref{app:davide}): $d_\L$ is thus the appropriate distance when analyzing properties of $\L$.

Moreover, the assumption on the time $T>\sup_{x\in \M} d_\L(x,\omega)$ is optimal because of the finite speed of propagation satisfied by equation \eqref{hypoelliptic-wave}. Indeed, Hypoelliptic wave equations also satisfy a finite speed of propagation similar to the classical wave equation but with the Riemannian distance replaced by the sub-Riemannian distance $d_\L$. This (not obvious) fact  was proved by Melrose \cite{Melrose:86} (see also the remarks in \cite[Section 4]{Jerison:87} for the link between the distance defined in \cite{Melrose:86} and $d_\L$).
\end{remark}

\medskip
As a corollary of this result (see~\cite{Robbiano:95} or \cite[Appendix]{LL:17}), we obtain the approximate controllability of the Hypoelliptic wave equation, as well as an estimate of the cost of approximate controls. Here, we only state approximate controllability to zero, which is equivalent to approximate controllability to the whole state space $\H^{1}_\L \times L^2$ on account to the reversibility of the equation.

\begin{corollary}[Cost of approximate control]
\label{corwavehypo}
For any $T> 2 \sup_{x\in \M} d_\L(x,\omega)$, there exist $C,c>0$ such that for any $\eps >0$ and any $(u_0 ,u_1) \in \H^{1}_\L \times L^2$, there exists $g \in L^2((0,T) \times \omega)$ with 
$$
\|g\|_{L^2((0,T) \times \omega)} \leq C e^{\frac{c}{\eps^k}} \nor{(u_0 ,u_1)}{\H^{1}_\L \times L^2} ,
$$
such that the solution of 
\bneq
(\d_t^2 + \L) u  = \mathds{1}_\omega g   && \text{ in } (0,T) \times \M , \\
(u ,\d_t u)|_{t=0} = (u_0 , u_1) ,&& \text{ in }\M ,
\eneq
satisfies $\nor{(u ,\d_t u)|_{t=T}}{L^2 \times \H^{-1}_\L} \leq \eps\nor{(u_0 ,u_1)}{\H^{1}_\L \times L^2}$. 
\end{corollary}

To the authors' knowledge, these results are the first ones concerning the approximate observability/controllability of hypoelliptic waves. They furnish not only the approximate observability/controllability but also an (optimal in general) estimate of the cost.

In the elliptic case $k=1$, these can be obtained by the theory developed by Lebeau in~\cite{Leb:Analytic} (even on a manifold with boundary). However, in this (elliptic) case, the analyticity assumption can be removed, as proved by the authors in~\cite{LL:15}. This followed a long series of papers concerning the {\em qualitative} unique continuation~\cite{RT:72,Lerner:88, Robbiano:91, Hormander:92, Tataru:95} (see also~\cite{RZ:98,Hor:97,Tataru:99} for more general operators), i.e. the property:
\bna
\big( (\d_t^2+ \L) u= 0 \text{ on }(-T,T) \times \M, \quad u = 0 \text{ on } (-T,T) \times \omega \Big) \implies u \equiv 0 ,
\ena
and another series of papers~\cite{Robbiano:95,Phung:10,Tatarunotes} concerning variants of Estimate~\eqref{th-estimate-k} (still in the elliptic case $k=1$) which are not optimal with respect to the minimal time and the exponent of $\mu$.
We refer to the introductions of~\cite{LL:15,LL:16} for a more detailed discussion on this issue.  
Here, in the analytic context, we directly prove the quantitative result but, to our knowledge, even the qualitative result was not known.

We shall see that we prove actually a more general statement in which the term $\nor{(u_0 ,u_1)}{\H^{1}_\L \times L^2}$ in the right-handside of Estimate~\eqref{th-estimate-k} can be changed into $\nor{(u_0 ,u_1)}{\H^{s}_\L \times \H^{s-1}_\L}$ for any $s>0$, if changing the power of $\mu$ accordingly, see Theorem~\ref{thmwavehypo-s} below.

\subsubsection{Quantitative approximate observability of the hypoelliptic heat equation}
\label{s:approx-control-heat}
We now turn to the study of observability properties for solutions of the hypoelliptic heat equation
\bneqn
\label{abstractheat}
\partial_t \y +\L \y=0 , && \text{ in } (0,T) \times \M , \\
\y(0)=\y_0  && \text{ in }\M ,
\eneqn
from a subdomain $\omega \subset \M$. By duality, we are equivalently concerned here with different controllability properties of the following system
\bneqn
\label{e:control-heat}
(\d_t + \L) u  = \mathds{1}_\omega g  ,& & \text{ in } (0,T) \times \M , \\
u(0) = u_0 ,&& \text{ in }\M .
\eneqn
We provide with three main results, still under Assumptions \ref{assumLiek} and \ref{hypoanal}:
\begin{enumerate}
\item For any $k\in \N^*$, we prove an approximate observability result in any time $T>0$ with a frequency-depending constant of order $Ce^{c\Lambda^k}$, where $\Lambda = \frac{\nor{\y_0}{\H^1_\L}}{\nor{\y_0}{L^2}}$, or, equivalently, approximate controllability with cost $e^{\frac{c}{\eps^k}}$. These are the analogues of Theorem~\ref{thmwavehypo} and Corollary~\ref{corwavehypo} for parabolic equations.
\item If we moreover assume the data to be sufficiently smooth (in some Gevrey-type norm with respect to the spectral decomposition of $\L$), then the cost of approximate observability can be improved to a polynomial one, i.e. of the form $\frac{C}{\eps^\beta}$ for some $\beta>0$; this yields approximate controllability in a much weaker topology, but with a much lower cost.
\item Finally, in the very particular case $k=2$ (including Grushin and Heisenberg operators), we prove an approximate observability/controllability property {\em to trajectories} in large time with a polynomial cost. This may be interpreted as a counterpart of the exact controllability to trajectories for the heat equation~\cite{LR:95,FI:96} (case $k=1$). There is no similar result if $k > 2$.
\end{enumerate}

The first result we obtain provides the cost of approximate observability of the whole state space $L^2(\M)$. There is no restriction for the hypoellipticity index $k$, but the (exponential) cost depends on this parameter.
\begin{theorem}
\label{t:approx-control-heat}
For all $T > 0$, there exist $C,c>0$ such that for any $\y_0 \in \H^1_\L$ and associated solution $\y$ of~\eqref{abstractheat}, we have
\bnan
\label{e:approx-control-heat1}
\nor{\y_0}{L^2}^2\leq Ce^{c\Lambda^k} \int_{0}^T \int_{\omega}\left|\y(t,x)\right|^2 dx~dt,\qquad \Lambda = \frac{\nor{\y_0}{\H^1_\L}}{\nor{\y_0}{L^2}} ,
\enan
and, for any $\mu>0$, 
\bnan
\label{e:approx-control-heat2}
\nor{\y_0}{L^2}^2\leq Ce^{c\mu^k} \int_{0}^T \int_{\omega}\left|\y(t,x)\right|^2 dx~dt + \frac{1}{\mu^2}\nor{\y_0}{\H^1_\L}^2 .
\enan
\end{theorem}
That~\eqref{e:approx-control-heat1} and~\eqref{e:approx-control-heat2} are equivalent comes for instance from~\cite[Lemma~A.3]{LL:15}.
Again, in the particular situation of Example~\ref{ex:Grushin++}, i.e. for the operators~\eqref{e:Grushin+}, the sequence of eigenfunctions of Proposition~\ref{Prop:BCG} shows that the exponent $e^{\kappa \mu^k}$ in~\eqref{e:approx-control-heat2} (resp. $e^{c\Lambda^k}$ in \eqref{e:approx-control-heat1}) cannot be improved in general.

This theorem generalizes the results of Fernandez-Cara-Zuazua and Phung \cite{FCZ:00,Phung:04} in the elliptic case $k=1$. Yet, in this framework, the analyticity was not necessary (as in all above stated results in the case $k=1$) and the setting can be relaxed (uniform dependence of the constants with respect to lower order terms and to the time $T$, boundary value problems...). 

As a corollary (see~\cite[Appendix]{LL:17}), we obtain, given an initial state and a target state both belonging to the space $L^2(\M)$, and given a precision $\eps$, the existence of a control function bringing the initial state in an $\eps$-neighborhood of the target (in appropriate topology). We obtain as well an estimate of the cost of the control.
\begin{corollary}[Cost of approximate control to the state space]
For any $T>0$, there exist $C,c>0$ such that for any $\eps >0$ and any $u_0 \in L^2(\M), u_1 \in L^2(\M)$, there exists $g \in L^2((0,T) \times \omega)$ with 
$$
\|g\|_{L^2((0,T) \times \omega)} \leq C e^{\frac{c}{\eps^k}} \nor{e^{-T\L} u_0-u_1}{L^2(\M)} ,
$$
such that the solution of~\eqref{e:control-heat} issued from $u_0$ satisfies 
$$
\nor{u(T) - u_1}{\H^{-1}_\L} \leq \eps \nor{ e^{-T\L} u_0-u_1}{L^2(\M)}.
$$ 
\end{corollary}

\bigskip
To state our second result concerning the hypoelliptic heat equation, we need to introduce the following spectral Gevrey-type norms for functions defined on $\M$: For $\alpha >0$, $\theta \in \R$, we set
\bnan
\label{e:gevrey-norms}
\nor{u}{\alpha,\theta}^2=\sum_{j \in \N} e^{2\theta \lambda_j^{\alpha}} |u_j|^2 , \quad \text{with}\quad u = \sum_{j \in \N} u_j \varphi_j .
\enan
For $\theta \geq 0$, we define $H^{\alpha,\theta}$ to be the subspace of $L^2(\M)$ consisting in functions $u$ such that $\nor{u}{\alpha,\theta}<\infty$. For $\theta <0$, we let $H^{\alpha,\theta}$ be the completed of linear combinations of eigenfunctions for this norm. Remark that taking as usual $L^2(\M)$ as a pivot space, the space $H^{\alpha,-\theta}$ is identified to $\big( H^{\alpha,\theta}\big)'$ for $\theta \geq 0$. Also, according to the hypoellipticity Assumption~\ref{assumLiek} (see Corollary~\ref{cor:HsHsL}) we have $H^{\alpha,\theta} \subset C^\infty(\M)$ for $\theta>0$, so that $H^{\alpha,-\theta}$ is larger than spaces of distributions on $\M$ (and its topology weaker).

We obtain the following result, which assumes the data to be extremely regular and then yields approximate observability with a polynomial cost only. The regularity needed is linked to the hypoellipticity index $k$. 
\begin{theorem}
\label{thm:para-gevrey}
Fix any $k \in \N^*$. There exists $\theta_{0}>0$ such that for any $T>0$ and any $\theta> \theta_{0}$, there exist $C>0$ so that for $\eps >0$, we have for any solution $\y$ to~\eqref{abstractheat},
 \bnan
\label{estimobserheathypohigh}
\nor{\y_0}{L^2}^2
\leq \frac{C}{\e^{\frac{\theta_{0}}{\theta-\theta_{0}}}}\int_{T/2}^T \int_{\omega}\left|\y(t,x)\right|^2~dt~dx+\e\nor{\y(0)}{k/2,\theta}^2.
\enan
\end{theorem}
 Again Proposition~\ref{Prop:BCG} shows that a polynomial cost is optimal for data in $H^{\theta, \frac{k}{2}}$ in the situation of Example~\ref{ex:Grushin++}. Note that we provide with an explicit dependence of the polynomial cost (i.e. the power $\frac{\theta_{0}}{\theta-\theta_{0}}$) with respect to the regularity of the data; in particular, we see how it improves when $\theta$ becomes larger.
As a Corollary of Theorem~\ref{thm:para-gevrey}, we obtain an approximate controllability result with polynomial cost, but the target is well approximated in a very weak topology.
\begin{corollary}[Cost of approximate control to the state space in very weak topology]
With $\theta_0>0$ given as in Theorem~\ref{thm:para-gevrey} (depending only on $\M, \omega, \L$), for any $T >0$ and $\theta >\theta_0$, there exist $C>0$ such that for any $\eps >0$ and any $u_0 \in L^2(\M), u_1 \in L^2(\M)$, there exists $g \in L^2((0,T) \times \omega)$ with 
$$
\|g\|_{L^2((0,T) \times \omega)} \leq \frac{C}{\e^{\frac{\theta_{0}}{\theta-\theta_{0}}}} \nor{e^{-T\L} u_0 - u_1}{L^2(\M)} ,
$$
such that the solution of~\eqref{e:control-heat} issued from $u_0$ satisfies 
$$
\nor{u(T) - u_1}{k/2, -\theta} \leq \eps \nor{e^{-T\L} u_0- u_1}{L^2(\M)}.
$$ 
\end{corollary}
We are not aware of any such results, even for the usual heat equation (i.e. with $k=1$). In this case, our proof also works in the $C^\infty$ context, and in the presence of boundaries, starting from the estimates obtained in \cite{LL:15} or the spectral estimates of \cite{LR:95}.

\bigskip
Our last main result concerning the hypoelliptic heat equation is, as opposed to the first two ones, concerned with final state approximate observability (or equivalently an approximate controllability to trajectories) with a {\em polynomial} cost, and is restricted to the case $k=2$.
\begin{theorem}
\label{thm:parabolic}
Assume that $k=2$. There exist $T_0 ,C>0$ such that for all $\eta >0$, all $T > T_0 + \eta$ and all $\eps >0$, we have for any solution $\y$ to~\eqref{abstractheat},
\bnan
\label{e:heat-approx-eps}
D \nor{\y(T)}{L^2}^2
\leq \frac{1}{\eps^\beta}
\int_{T-\eta}^T \int_{\omega}\left|\y(t,x)\right|^2~dt~dx+\eps \nor{\y(0)}{L^2}^2 , 
\enan
with $D=\min\{ \frac{e^{-C/\eta}}{C}, 1\}$ and $\beta = \frac{T_0}{T-(T_0 +\eta)}$.
\end{theorem}

In particular, we obtain an explicit estimate on how the cost improves as $T$ increases.
This result gives directly the following corollary concerning approximate controllability to trajectories (or, equivalently, to zero) at a polynomial cost (see again~\cite[Appendix]{LL:17}).
\begin{corollary}[Cost of approximate control to trajectories if $k=2$]
\label{cor:parabolic}
Assume that $k=2$, and let $T_0>0$ as in Theorem~\ref{thm:parabolic}. For all $\eta >0$, all $T > T_0 + \eta$ and all $\eps >0$, we have the following statement: for any $u_0, \tilde{u}_0\in L^2$, there exists $g \in L^2((0,T) \times \omega)$ with 
$$
\|g\|_{L^2((0,T) \times \omega)} \leq  \frac{\tilde{C}}{\eps^{\beta}} \nor{u_0- \tilde{u}_0}{ L^2} ,
$$
such that the associated solution $u$ of~\eqref{e:control-heat} satisfies 
$$
\nor{u(T) - e^{-T\L}\tilde{u}_0}{L^2(\M) } \leq \eps\nor{u_0 -\tilde{u}_0}{L^2}, 
$$ 
where $\beta = \frac{T_0}{T-(T_0 +\eta)}$ and $\tilde{C} = \tilde{C}(\eta, T_0,T)$.
\end{corollary}
 Remark that these two results only hold in the case $k=2$. Once again, in the particular situation of Example~\ref{ex:Grushin++}, i.e. for the operators~\eqref{e:Grushin+}, the sequence of eigenfunctions of Proposition~\ref{Prop:BCG} shows that this result cannot hold if $k\geq 3$ (or even, if $\gamma>1$). Let us here be more precise: assume in the context of Example~\ref{ex:Grushin++} that an estimate of the form~\eqref{e:heat-approx-eps} is satisfied for a cost function $\Phi(\eps)$, that is 
\bna
\nor{\y(T)}{L^2}^2
\leq \Phi(\eps) \int_{0}^T \int_{\omega}\left|\y(t,x)\right|^2~dt~dx+\eps \nor{\y(0)}{L^2}^2, \quad \text{ for all }\eps>0 ,
\ena
and test it with $\y(t)=e^{-\lambda_jt} \varphi_j$ (solution of \eqref{abstractheat}), where $\varphi_j$ is given by Proposition~\ref{Prop:BCG}. Then we have for all $\eps>0$,
$$
e^{-2\lambda_j T} \leq \Phi(\eps) \frac{e^{-2c_0 \lambda_j^{k/2}}}{\lambda_j}  + \eps , \quad \text{ for all }\eps>0 .
$$
Fixing then $\eps = \eps_j : = \frac{e^{-2\lambda_j T}}{2} \to 0^+$ and taking logarithm yields $2(c_0 \lambda_j^{k/2}-\lambda_jT)  \leq \log\Phi(\eps_j)$. That $\Phi(\eps) \leq C/\eps^\beta$ (i.e. having a polynomial cost) implies $k \leq 2$.  In the case $k=2$, this implies $\Phi(\eps_j)\geq \frac{1}{(2\eps_j)^{\frac{c_0}{T}-1}}$, so that for $T< c_0$, the polynomial cost cannot be improved. Unfortunately, the constant $T_0$ in the above result is much larger than $c_0$, so that this discussion does not imply neither that the polynomial cost is the optimal one, nor that a minimal time is necessary. 
However, as we shall see in Section~\ref{s:discusison-control} below, it may happen that exact controllabilty holds for no time $T>0$, which may indicate that a polynomial cost is not far from being sharp.

 \subsubsection{Relaxing the analyticity assumption}
\label{s:without-analyticity}

In this section, we provide with a simple family of examples for which the analyticity Assumption~\ref{hypoanal} can be partially removed. Still this family contains those of Examples~\ref{ex:Grushin+} and~\ref{ex:Grushin++}. In this context, most above results hold as well. The motivation is both to relax the analyticity assumption (and replace it with analyticity with respect to one set of variables), and to include the setting of the article~\cite{BeauchardCanGugl}.

\begin{example}[Partially analytic Grushin-type operators]
\label{ex:Grushin+++}
Consider the manifold with boundary $\M=[-1,1]\times(\R/\Z)$, endowed with the Lebesgue measure $dx$, and define, for $f\in C^{\infty}([-1,1] \times (\R/\Z))$, the Grushin type operator 
\bnan
\label{e:def-op-L-partially}
 \L = X_1^*X_1 + X_2^* X_2, \qquad X_1 = \d_{x_1}, \quad X_2 =  f(x_1,x_2) \d_{x_2} , 
\enan
that is 
$$
\L = -\d_{x_1}^2 - f^2 \d_{x_2}^2 - (2 f \d_{x_2} f )\d_{x_1} ,
$$
with Dirichlet conditions on $\d \M$.
We further assume that
\begin{itemize}
\item $f(x_1,x_2)$ is {\em analytic} in the variable  $x_2$ (that is, for any point $x= (x_1 , x_2) \in ]-1,1[\times(\R/\Z)$, $f$ is equal to its partial Taylor expansion at $x_2$ with respect to the variable $x_2$ uniformly in a neighborhood of $x$ in $]-1,1[\times(\R/\Z)$);
\item there exists $\e>0$ such that $f(x_1,x_2)=f(x_1)$ does not depend on $x_2$, and $f(x_1)>0$ for all $x_1 \in [-1,-1+\e]\cup [1-\e,1]$;
\item $X_1$ and $X_2$ satisfy the Chow-Rashevski-H\" ormander Assumption~\ref{assumLiek};
\item $f(x_1,x_2)$ does not depend on $x_2$ in a neighborhood of $\omega$.
\end{itemize}
\end{example}

Note that under these assumptions, the operator $\L$ is elliptic near the boundary $\d \M$.

For instance if $f(x_1,x_2) = f(x_1) \in C^\infty([-1,1])$ does not depend on the variable $x_2$, and we have :
$$
f(x_1) \neq 0 , \text{ for } x_1 \neq 0 \qquad f^{(\alpha)}(0) = 0 ,\text{ for all } \alpha \leq k -2, \qquad \text{ and }  f^{(k-1)}(0) \neq 0 ,
$$
then, the operator $\L$ defined by~\eqref{e:def-op-L-partially}, namely $\L =  -\d_{x_1}^2 - f(x_1)^2 \d_{x_2}^2$, satisfies all assumptions of Example~\ref{ex:Grushin+++} (in particular, it is hypoelliptic of order $k$). This contains the situation of Examples~\ref{ex:Grushin+} and~\ref{ex:Grushin++} (for $\gamma \in \N$ of course).

We prove the following result.
\begin{theorem}
\label{t:partially-anal}
In the context of Example~\ref{ex:Grushin+++}, all results of Theorems~\ref{t:spec-ineq},~\ref{thm:para-gevrey} and~\ref{thm:parabolic} still hold, as well as their corollaries.

Theorem \ref{thmwavehypo} is true with the following estimate instead
\bnan
\label{th-estimate-partial}
\nor{(u_0,u_1)}{L^2\times \H^{-1}_\L} \leq C e^{\kappa \mu}\nor{u}{L^2(]-T,T[\times \omega)} +\frac{1}{\mu}\nor{(u_0,u_1)}{\H^{k}_\L\times \H^{k-1}_\L}.
\enan
Theorem \ref{t:approx-control-heat} is still true but with the estimates
\bnan
\label{e:approx-control-heat1partial}
&&\nor{\y_0}{L^2}^2\leq Ce^{c\Lambda_k} \int_{0}^T \int_{\omega}\left|\y(t,x)\right|^2 dx~dt,\qquad \Lambda_k = \frac{\nor{\y_0}{\H^k_\L}}{\nor{\y_0}{L^2}} , \\
\label{e:approx-control-heat2partial}
&&\nor{\y_0}{L^2}^2\leq Ce^{c\mu} \int_{0}^T \int_{\omega}\left|\y(t,x)\right|^2 dx~dt + \frac{1}{\mu^2}\nor{\y_0}{\H^k_\L}^2 .
\enan
\end{theorem}
Note that since $f$ does not vanish near $\d \M$, the metric $g$ defined as in \eqref{def-g} is Riemannian near $\d \M$ and the notions of length and distance defined above can be extended up to the boundary.

We explain in Section~\ref{s:non-analytic} how the proofs in the completely analytic case need to be modified. The version that are true in the partially analytic are actually some particular cases of general estimates with all Sobolev scales for measuring the typical frequency. For instance, \eqref{th-estimate-partial} is a particular case of Theorem \ref{thmwavehypo-s} with $s=k$. We refer to the discussion in Section~\ref{subsectrkintro} below.

\begin{remark}
Under appropriate assumption on $f$, it is classical to extend Theorem~\ref{t:partially-anal} to the same situation as in Example~\ref{ex:Grushin+++}, but on the domain $\M= [-1,1]_{x_1}\times [-1,1]_{x_2}$ with Dirichlet boundary conditions by using symmetry arguments. Also, the case of the domain $\M= (\R/\Z)^2$ is simpler.
\end{remark}
\begin{remark}
All observability results of Theorem~\ref{t:partially-anal} also hold if the internal observation term $\nor{u}{L^2(]-T,T[\times \omega)}$ is replaced by a boundary observation $\nor{\d_n u}{L^2(]-T,T[\times \Gamma)}$, where $\Gamma$ is a nonempty open subset of $\d \M$ and $\d_n$ denotes the normal derivative to $\d\M$. See~\cite[Section~5]{LL:15}. In turn, they imply their boundary controllability counterparts. 
We did not state these results for the sake of brevity.
\end{remark}
\subsection{Comparison to other works}
\label{s:comp-to-other}
\subsubsection{Previous works on unique continuation for sum-of-squares operators}
\label{s:comp-UC}
Observability inequalities (as those provided by Theorems~\ref{t:spec-ineq},~\ref{thmwavehypo},~\ref{t:approx-control-heat},~\ref{thm:para-gevrey} and~\ref{thm:parabolic} above, or as~\eqref{obserGrushinIntro} below) are quantitative estimates of the unique continuation property for the operator involved (namely $\L-\lambda$, $\d_t^2+\L$, and $\d_t+\L$ respectively). Hence, when studying such inequalities, it is natural to compare our results with the known unique continuation properties for such operators. When the ellipticity condition is dropped, i.e., when $k>1$, this property seems to be a very intricate problem, even under the simplest form 
\bna
\big(  \L u  = 0   \text{ on }\M , \quad u = 0 \text{ on } \omega \big) \implies u =  0\text{ in a neighborhood of }\omega \text{ in } \M .
\ena
To our knowledge, the most general such result was proved by Bony \cite{Bo:69}, and holds under both the Chow-Rashevski-H\"ormander condition and the assumption that the coefficients of the operator are {\em analytic}. Therefore, our assumptions~\ref{assumLiek} and~\ref{hypoanal} (except in the partially analytic case of Theorem~\ref{t:partially-anal}) are essentially the same as in this paper. In particular, Theorem~\ref{t:spec-ineq} could be read as quantification of Bony's result. Also, the proof of the result of Bony mainly relies on the Holmgren-John theorem and is quite indirect. Here, we need to make a new proof of his result, that we are also able to quantify. 

Some attempts have been done to relax this analyticity assumption. Watanabe \cite{W:82} proved the unique continuation property for $C^{\infty}$ coefficients in dimension $d=2$. 
Yet, later on, Bahouri~\cite{Ba:86} proved a surprising general non-uniqueness result: for a large class of  sum-of-squares operators $\L$ with  $C^{\infty}$ coefficients, and satisfying Assumption~\ref{assumLiek}, there is $C^\infty$ potentials $V$ such that $\L+V$ does not satisfy the local unique continuation property.
These counterexamples to unique continuation contain for instance in dimension $d=3$ and $d=4$ the case where the dimension spanned by the vector fields is of dimension $d-1$ (Heisenberg-like situations). Moreover, this result suggests that a classical Carleman estimate approach cannot work for all hypoelliptic operators. Also, it strongly suggest that the  (complete or partial) analyticity assumptions that we make are not completely artificial.

This analyticity assumptions might be completely removed in some specific situations where the operator is elliptic outside of a submanifold, see the comments of Bahouri \cite[p140]{Ba:86}. This was proved in the paper~\cite{Ga:93} by Garofalo for specific examples. Colombini-Del Santo-Zuily \cite{CDZ:93} also treated some related classes of degenerate elliptic operators having a specific form with respect to a hypersurface. Nevertheless, even in these situations, the quantitative estimates that we obtain are optimal as stated in Proposition \ref{Prop:BCG}.

All these results are concerned with the unique continuation property for operators like $\L$ (``degenerate elliptic operators''). We are not aware of works studying the unique continuation property for operators like $\d_t^2 + \L$ or $\d_t + \L$ (``hyperbolic, resp. parabolic operators with a degenerate elliptic part''), except in the context of control theory, that we review in the next section.

\subsubsection{Previous works on the controllability of the hypoelliptic heat equation}
\label{s:discusison-control}
The investigation of the controllability of hypoelliptic operators is for the moment quite at an early stage and has been mainly restricted to some specific operators or classes of operators. A striking result concerning the parabolic observation problem~\eqref{abstractheat}, where $\L = \L_\gamma$ is given by Example~\ref{ex:Grushin++} (i.e. higher order Grushin operators on the rectangle, with Dirichlet boundary conditions), was proved by Beauchard, Cannarsa and Guglielmi \cite{BeauchardCanGugl}.
The authors are interested in the following observablity inequality, equivalent to the controllability to zero (and hence to trajectories)
\bnan
\label{obserGrushinIntro}
\nor{\y (T)}{L^2(\M)}^2\leq C \int_0^T\int_{\omega}|\y (t,x)|^2 ~dt dx , \quad \text{for all } \y_0\in L^2(\M) \text{ and $\y$ solution of~\eqref{abstractheat}.}
\enan
\begin{theorem}[Beauchard, Cannarsa and Guglielmi \cite{BeauchardCanGugl}]
Assume $\L = \L_\gamma$ is given by Example~\ref{ex:Grushin++}.
\begin{enumerate}
\item If $\gamma\in [0,1[$, then the observability inequality \eqref{obserGrushinIntro} holds true for any nonempty open set $\omega \subset \M$ in any time $T > 0$.
\item If $\gamma=1$ and if $\omega=]a, b[ \times ]0,1[$ where $0 < a < b < 1$, then there exists $T^*\geq a^2/2$ such that
\begin{itemize}
\item for every $T > T ^*$ the observability inequality \eqref{obserGrushinIntro} holds true,
\item for every $T < T^*$ the observability inequality \eqref{obserGrushinIntro} is false.
\end{itemize}
\item If $\gamma>1$ and $\omega \subset (0,1)\times (0,1)$, then the observability inequality \eqref{obserGrushinIntro} never holds true, in any time $T > 0$.
\end{enumerate}
\end{theorem}
In the case $\gamma=1$ (i.e. $k=2$) and with a symmetric observation region $\omega=(]-b,  -a[ \cup ]a, b[ )\times ]0,1[$ with $0<a<b\leq 1$, it has been recently proved by Beauchard, Miller and Morancey \cite{BeauMillMor15} that $T^*=a^2/2$ is actually the critical time.
This result is quite surprising since parabolic type equations often display an infinite speed of propagation. The controllability of parabolic evolutions thus usually holds in an arbitrary small time independent on the geometry; appearance of a minimum controllability time is hence unusual. Yet, the proof uses a lot the specific geometry of $\omega$ as a vertical strip. Indeed, another very striking result was recently proved by Koenig in the case $\gamma = 1$ (i.e. for the Grushin operator): if $\omega$ is disjoint from an horizontal strip, null-controllability never holds (in any time). 
\begin{theorem}[Koenig~\cite{Koe:17}]
\label{thmKoenig}
Let $\L = \L_\gamma$ be given by Example~\ref{ex:Grushin++} with $\gamma=1$. Assume that there is $0 < c < d < 1$ such that $\omega\cap \big(  ]-1, 1[ \times ]c,d[ \big) =\emptyset$. Then, for any $T >0 $ the observability inequality \eqref{obserGrushinIntro} is false.
\end{theorem}
A remarkable consequence of this result, when compared with the two abovementioned ones is that a {\em geometric condition} on the set $\omega$ is needed for the observability estimate~\ref{obserGrushinIntro} to hold. 

Hence, the best result one can then expect in a general situation is a final state approximate observability result with a cost function $\Phi(\eps) \to_{\eps\to 0^+} + \infty$ (or equivalently an approximate controllability to trajectories with cost $\Phi(\eps)$), which is precisely our Theorem~\ref{thm:parabolic} and Corollary~\ref{cor:parabolic} with $\Phi(\eps) = \eps^{-\beta}$. 

Finally, let us also underline that all these result hold in the context of Example~\ref{ex:Grushin++}, that is for the operator $-(\d_{x_1}^2+ x_1^{2\gamma}\d_{x_2}^2)$, which coefficients are analytic (indeed constant) with respect to the variable $x_2$. As such, they fit into the framework of Theorem~\ref{t:partially-anal} as long as $\gamma \in \N$.

More recently, there has been some study of the control of the Heat equation on the Heisenberg group by Beauchard and Cannarsa \cite{BeauCan:17}. It still corresponds to the case $k=2$ as described in Example \ref{exHeisen}. Some phenomenon similar to the Grushin case seem to occur with the existence of a minimal time if the observability is made on a cylinder. This strenghtens the fact suggested by our result that the important parameter is the hypoelliptic index $k$. 

Yet, in both cases of Grushin and Heisenberg (or more generally when $k=2$), it remains to understand what is the geometric property on the observation set $\omega$ that makes the difference between the polynomial cost that is provided by our result without any assumption on $\omega$ (which is likely to be optimal in general as suggested by Theorem \ref{thmKoenig}) and the exact controlability that requires some geometrical assumptions on $\omega$.
\subsubsection{Controllability of other equations driven by degenerate elliptic operators}
To conclude this section, let us mention different works related to the controllability of parabolic equations driven by hypoelliptic or degenerate elliptic operators, that do not fit in the framework of the present article.

\medskip
First, the paper~\cite{Mor:15} by Morancey treats the approximate controllability (or the unique continuation property) for the heat equation associated with the Laplace Beltrami operator of the Grushin sub-Riemannian metric defined in~\cite{BoscLaurent:13}. This operator is equal to the Grushin operator discussed in Example~\ref{ex:Grushin} plus a singular potential on the singular set $x_1=0$. Hence, the analysis of the cost associated to approximate controls is much beyond the scope of the present paper.

\medskip
Second, we only considered here type~I (selfadjoint) H\"ormander operators, that is $\L=\sum_{i=1}^m X_i^*X_i$. Another classical class of hypoelliptic operators consists in type~II H\"ormander operators, namely $\L=\sum_{i=1}^m X_i^* X_i+X_0$, where the vector field (the drift) $X_0$ is necessary to span the full tangent space with iterated Lie brackets. These are no longer selfadjoint operators. The simplest example is the so-called Kolmogorov (or Fokker-Planck) operator $\L=-\partial_v^2 + v\partial_x$. Our results do not apply in this setting, especially because our main theorems only see the principal symbol of the operator. Yet,  recent progress has been made to analyse the observability/controllability of parabolic equations driven by such operators (mainly for some variants of the Komogorov operator, though). We quote for instance the papers \cite{BeauchZuaz, BeauchardKolm:14, BeauchardHellRob}. The observability/controllability problem has also been considered on the whole space $\R^d$. This led to other geometrical problems about how the domain is "spread out" at infinity, see for instance Le Rousseau-Moyano \cite{LM:16} for the Kolmogorov equation and Beauchard-Pravda-Starov~\cite{BeauPravd:16,BeauPravd:16b} for some more general class of quadratic operators. In these last papers, the idea of the proof is to combine some observability of low frequency (where the frequency are the usual Fourier ones or related to the harmonic oscillator) with some decay and regularizing properties of the semigroup. 

It would be very interesting to understand the common features and differences of our results and methods with these ones. At first sight, it seems that in both cases, one important idea (which goes back to Lebeau-Robbiano \cite{LR:95}) is to compare the decay rate of the heat equation with the cost of observability (or control) of low frequency. But, in our paper, we define "frequency" with respect to the hypoelliptic operator $\L$ and then, the decay is natural for high frequency. The hard part is then to understand the cost of observability of low frequencies. In these papers, it seems that they have chosen to define the frequency as the usual Euclidian one (or with respect to a fixed well known operator as the harmonic oscillator in \cite{BeauPravd:16b}). In this situation, the observation at low (Euclidian) frequency does not really see the hypoelliptic operator and follows from more usual Carleman estimates. Yet, the decay rate of high frequency (see for instance Proposition 2.2 of \cite{BeauPravd:16}) and the understanding of the commutation with Fourier cutoff turns out to be much more complicated and to reflect deeply the hypoelliptic properties of the operator. 

Note also that there has been several studies of the unique continuation property for type~II H\"ormander operators (see e.g. \cite{LZ:82} for related operators).

\medskip
Finally, other types of degeneracies have also been studied, as for instance elliptic operators with coefficients vanishing near the boundary of a domain $\M \subset \R^d$. In this case, adaptations of the usual Carleman estimates (combined with appropriate Hardy inequalities) are sometimes tractable. The literature is vast, and we simply mention the recent memoir~\cite{CannaMartiVan:16} and refer the reader to the references therein.

\subsection{Sketch of the proofs and plan of the paper}
\label{s:sketch-plan}

Even though this is not explicit in the discussion above, the cornerstone result of this paper is Theorem~\ref{thmwavehypo}, concerning the hypoelliptic wave equation. All results concerning eigenfunctions (Theorem~\ref{t:spec-ineq}) or the hypoelliptic heat equation (Theorems~\ref{t:approx-control-heat},~\ref{thm:para-gevrey} and~\ref{thm:parabolic}) are then deduced from Theorem~\ref{thmwavehypo}. The proof in the partially analytic case (Theorem~\ref{t:partially-anal}) shall be discussed afterwards. Let us hence first comment the proof of Theorem~\ref{thmwavehypo}.

\medskip
The proof of Theorem~\ref{thmwavehypo} is based on the general strategy developed by the authors in~\cite{LL:15} for quantifying and propagating unique continuation properties. From~\cite{LL:15}, we only use here (except for the partially analytic situation of Theorem~\ref{t:partially-anal}) the ``Holmgren-John'' case, i.e. when the operator has analytic coefficients. It states basically 
\begin{itemize}
\item that an appropriate quantitative (low frequency) estimate holds across any non-characteristic hypersurface;
\item that such local estimates can be propagated, leading towards global ones.
\end{itemize}
In Section~\ref{s:Holmgren-john-LL15}, we review results and tools developed in~\cite{LL:15}; for sake of readability, we specify the latter to the very particular case of {\em second order} operators that are elliptic when restricted to $\zeta_a=0$ (the cotangent variable to the analytic variable, called $\xi_a$ in~\cite{LL:15}), which includes all operators studied in the present article.

Here, when compared to the case of the classical wave equation, two more difficulties arise: one of geometric nature, and one related to the compatibility of energy space associated to $\L$ and those dealt with in~\cite{LL:15}.

Let us first describe the geometric difficulty. The proof is inspired by the case of the classical wave equation given in~\cite[Section~6.1]{LL:15}: the idea is, given a point $x_0 \in \M$, to take any path $\gamma: [0,1]\to \M$ with $\gamma(0)=x_0$ and $\gamma(1)\in \omega$ (observation set), of length sufficiently small, and then to construct a family of appropriate noncharacteristic hypersurfaces in these coordinates near $[-T,T] \times \gamma$. There, we apply the general theorem of~\cite{LL:15}, which allows to bound the solution $u$ to $(\d_t^2-\Delta)u=0$ in a neighborhood of $(t,x) = (0,x_0)$ by $u$ in $[-T,T]\times \omega$.

Here, due to the non definiteness/ellipticity of the operator $\L$, we are not able to construct {\em global} coordinates near {\em any} path $\gamma$ together with appropriate noncharacteristic hypersurfaces, in which to apply the results of~\cite{LL:15}. To overcome this difficulty, we do not consider {\em any} path between $x_0$ and $\omega$, but rather only so called {\em normal geodesics}, that is, projections on $\M$ of hamiltonian curves of the principal symbol of the operator $\L$. The existence of such paths $\gamma$ (minimizing the sub-Riemannian distance) from any point $x_0$ to $\omega$ is a well-known result in sub-Riemannian geometry, proved by Rifford and Tr\'elat~\cite{RiffordTrelatMorse}. Then, locally near a point of $\gamma$, the introduction of normal geodesic coordinates allows us to define local coordinates in which to apply a {\em local} version of our results in~\cite{LL:15}.

A new difficulty, linked to the methods used in~\cite{Tataru:95, Tataru:99, RZ:98,Hor:97, LL:15}, then arises: the whole setting of these papers relies on a splitting of space into analytic and non-analytic coordinates. Hence, most ``patchable estimates'' (linked to a relation $\lhd$, see Section~\ref{s:def-mult}) produced in~\cite{LL:15} require the analytic variable to be global and straight, which is obviously not the case here. To solve this problem we do not rely on the main (neither global, nor local) result of~\cite{LL:15}, but rather on the specific result of~\cite[Theorem~4.11]{LL:15}, which takes into account the possible changes of variables. Having these results in hand allows to prove an estimate of the form
\bnan
\label{e:intro-estim-partial-wave}
 \nor{ u}{L^2(]- \e, \e[\times \M)}\leq C e^{\kappa \mu}\nor{u}{L^2(]-T,T[\times \omega)}+\frac{C}{\mu}\nor{u}{H^{1}(]-T,T[\times \M)} ,
\enan
for $\mu$ large and $u$ solution to $(\d_t^2+\L)u=0$. This estimate is the same as that obtained in~\cite{LL:15} for the wave equation.

\medskip
This leads us to the second main difficulty we have to face in the proof of Theorem~\ref{thmwavehypo}. 
Whereas the left hand-side of~\eqref{e:intro-estim-partial-wave} is bounded from below by the natural $L^2\times \H_\L^{-1}$ norm of the data, the right hand-side is not directly linked to their $\H_\L^1 \times L^2$ norm.
More precisely, the hypoelliptic estimates Rothschild and Stein~\cite{RS:76} (see Theorem~\ref{thmhypoestim} above and Appendix~\ref{app:Hs-comm}) imply that $\nor{u}{H^{1}(]-T,T[\times \M)} \leq C \nor{(u_0,u_1)}{\H_\L^k \times \H_\L^{k-1}}$. This provides a weaker version of Theorem~\ref{thmwavehypo} which has exactly the same form as in the case of the wave equation (cost $e^{\kappa \mu}$), but with the norm $\nor{(u_0,u_1)}{\H_\L^k \times \H_\L^{k-1}}$ in the right hand-side. This weaker version is however interesting for itself since the proof is much less involved, and we prove it in Section~\ref{s:simple-case}. 

To obtain the estimate of Theorem~\ref{thmwavehypo} (and in fact, a family of such estimates with any $\H_\L^s \times \H_\L^{s-1}$, $s>0$, in the right hand-side, see Theorem~\ref{thmwavehypo-s} below), we thus need to work with a version of \eqref{e:intro-estim-partial-wave} still containing frequency cutoff localization and an $e^{-c\mu}$ small remainder (instead of the $1/\mu$ one). These low-frequency-with-exponentially-small-remainder estimates are then combined with the spectral representation of solutions to $(\d_t^2+\L)u=0$ in order to gain back derivatives in the remainder term. Such estimates are close to those we prove in~\cite{LL:17} for the classical wave equation.
These final energy estimates are performed in Section~\ref{s:energy-general-case}, and conclude the proof of Theorem~\ref{thmwavehypo}.

\bigskip
Starting from Theorem~\ref{thmwavehypo}, let us now explain how to deduce the other results of the paper, namely Theorems~\ref{t:spec-ineq},~\ref{t:approx-control-heat},~\ref{thm:para-gevrey} and~\ref{thm:parabolic}. First of all, Theorems~\ref{t:spec-ineq} is simply deduced from Theorem~\ref{thmwavehypo} by using a particular solution to the wave equation~\eqref{hypoelliptic-wave}, namely $u(t,x)= \cos(\sqrt{\lambda_j }t)\varphi_j(x)$. See Section~\ref{s:interlude}.

\medskip
Section~\ref{s:hypo-heat} is devoted to the proofs of Theorems~\ref{t:approx-control-heat},~\ref{thm:para-gevrey} and~\ref{thm:parabolic}, which follow the general idea that the controllability/observability properties for hyperbolic equations implies controllability/observability properties for their parabolic counterpart, see~\cite{Russell:73,Miller:06b,EZ:11,EZ:11s} (see also~\cite{LR:95}). This has been named as ``transmutation methods'' by Luc Miller~\cite{Miller:06b}. Here, we use the method developed in~\cite{EZ:11}. In that paper, Ervedoza and Zuazua deduced the (exact final time) observability of the heat equation (known from~\cite{LR:95,FI:96}) from the approximate observability estimate for waves (namely the analogue of Theorem~\ref{thmwavehypo}) as proved  in~\cite{Phung:10} (with loss) or~\cite{LL:15} (without loss). Their proof consists in constructing an appropriate kernel $k_T(t,)$ such that if $\y(t)$ is a solution to the usual heat equation, $u(t)=\int_0^T k_T(t,s)\y(s)ds$ is a solution to the usual wave equation, to which we can apply the analogue of Theorem~\ref{thmwavehypo}. Because of the exponential cost in term of the frequency ($e^{c\Lambda}$), the resulting estimates are only useful at low frequency: for data having (spectral) frequencies $\sqrt{\lambda_j}\leq \sqrt{\lambda}$, one then obtain observability (or controllability if we think about the dual problem) at cost $e^{c\sqrt{\lambda}}$ as in~\cite{LR:95}. The proof of final state observability then follows from comparing this cost with the heat dissipation for frequencies $\sqrt{\lambda_j}\geq \sqrt{\lambda}$, namely $e^{-t \lambda}$ as in the original proof~\cite{LR:95} (see also~\cite{LeLe:09} or the simplified argument of~\cite{Miller:10}). 

Here, we follow the approach of~\cite{EZ:11} (in particular, we use the same kernel $k$ and its properties) in the proofs of Theorems~\ref{t:approx-control-heat},~\ref{thm:para-gevrey} and~\ref{thm:parabolic}, with the following modifications.

The proof of Theorems~\ref{thm:para-gevrey} and~\ref{thm:parabolic} are vey close to that of~\cite{EZ:11}. However, application of the method of~\cite{EZ:11} yields that the observability of low frequencies $\sqrt{\lambda_j}\leq \sqrt{\lambda}$ costs $e^{c\lambda^{k/2}}$, see Lemma~\ref{l:estim-heat-BF} (remark that Proposition~\ref{Prop:BCG} implies that this is optimal in general). This cost has to be compared to the dissipation for high frequencies $\sqrt{\lambda_j}\geq \sqrt{\lambda}$, namely $e^{-t \lambda}$. Hence, we see that the case $k=1$ (classical heat equation, already discussed), $k = 2$, and $k>2$ display very different features:
\begin{enumerate}
\item In case $k=2$, the cost of observation of low frequencies $e^{c\lambda}$ and the parabolic dissipation for high frequencies $e^{-t\lambda}$ have the same strength: in this case, we need to wait a time long enough so that the dissipation ``beats'' the cost of the observability (essentially $t>c$). Moreover, the iterative procedure devised in~\cite{LR:95} in order to control/observe {\em all} frequencies in finite time cannot converge here: each step would need a time $t>c$. Therefore, we only obtain the approximate controllability result of Theorem~\ref{thm:parabolic}, with a cost improving as time increases. See Section~\ref{s:proof-thm:parabolic} for the proof of Theorem~\ref{thm:parabolic}.
\item In case $k>2$, the dissipation for high frequencies $e^{-t\lambda}$  has no chance to compete with the cost of observation of low frequencies $e^{c\lambda^{k/2}}$.
Assuming that the initial data are in the Gevrey-type space $H^{\theta, \alpha}$ with $\alpha= k/2$ allows to compensate for the cost of low frequencies $e^{c\lambda^{k/2}}$ (the $\theta$ having to be compared to $c$), leading to Theorem~\ref{thm:para-gevrey}. Note that parabolic dissipation at high frequencies does not play any role here: low frequencies are observed thanks to transmutation and high-frequencies are absorbed by the Gevrey norm. The cases $k=1,2$ (in which the Gevrey norm is relatively weaker) in Theorem~\ref{thm:para-gevrey} are a little different and require elements similar to those used in the proof of Theorem~\ref{t:approx-control-heat}.
Note that this type of result seems to be new for the classical heat equation as well (in which case our proof also holds in a much more general setting).
Proof of Theorem~\ref{thm:para-gevrey} is performed in Section~\ref{s:proof-thm:para-gevrey}.
\end{enumerate}

Finally, the proof of Theorem~\ref{t:approx-control-heat} in Section~\ref{s:proof-t:approx-control-heat} relies on the same transmutation technique. However, we do not split the solution into low and high-frequencies, but rather apply the transmutation kernel $k_T(t,s)$ to the full solution $\y$ to the heat equation: $u(t)=\int_0^T k_T(t,s)\y(s)ds$ is a solution to the wave equation. We then prove a fine asymptotics analysis of $\int_0^T k_T(0,s)e^{-\lambda s}ds$ for high frequencies together with convexity estimates to bound the frequency function of $u(0)$ by the frequency function of $\y(0)$, namely $\Lambda = \frac{\nor{\y(0)}{\H^1_\L}}{\nor{\y(0)}{L^2}}$. 
The proof of this result via a direct transmutation method seems to be new, even for the classical heat equation. The usual proofs~\cite{FCZ:00,Phung:04} rather rely on the exact final time observability estimate, which does not hold here in general. However, as opposed to~\cite{FCZ:00,Phung:04}, we do not recover uniform estimates in terms of the control time $T$ as $T\to 0^+$.

\medskip

Finally, in Section~\ref{s:non-analytic}, we prove the partially analytic result of Theorem~\ref{t:partially-anal}. Only the analogue of Theorem~\ref{thmwavehypo} at regularity $\H^{k}_{\L}$ needs to be proved (namely estimate \eqref{th-estimate-partial}), since, as discussed above, all results of Theorems~\ref{t:spec-ineq},~\ref{t:approx-control-heat},~\ref{thm:para-gevrey} and~\ref{thm:parabolic} (under the appropriate form) are corollaries of that of Theorem~\ref{thmwavehypo}.
The situations is almost the same as that of Theorem~\ref{thmwavehypo} except for four main differences.
First, the presence of the boundary makes it complicated to apply globally the geometric result of Rifford and Tr\'elat~\cite{RiffordTrelatMorse}, and we only rely on a local version of it.
Second, the partial analyticity assumption does not allow to make changes of variables. This difficulty is overcome by the very simple geometry of $[-1,1]_{x_1}\times \T_{x_2}$, in which we barely do not perform any change of variable.
Third, the application of the results in~\cite{LL:15} yields an observation term in a mixed $L^2-H^1$ norm; we have to refine this estimate to recover the $L^2$ observation term.
Finally, the available hypoelliptic estimates do not apply directly in the presence of boundary and we have to patch hypoelliptic estimates in the interior with elliptic estimates at the boundary.

\medskip
The paper ends with three appendices, the first of which, Appendix~\ref{s:proof-Prop:BCG} is devoted to the proof of the optimality result of Proposition~\ref{Prop:BCG} using some estimates of~\cite[Section~2.3]{BeauchardCanGugl}. The second part, Appendix~\ref{app:sectsub} provides the proof of several subelliptic estimates that are used throughout the paper. They are consequences of Theorem \ref{thmhypoestim}.
Finally, Appendix~\ref{app:davide} contains an technical computation.

\subsection{Some remarks and further comments}
\label{subsectrkintro}

This section contains several remarks concerning the setting of the present paper and the results we obtain.
\begin{remark}[Sub-Riemannian Laplacians]
\label{r:sub-riem-lap}
Here, we explain why the assumption that $\L$ writes as a sum of $X_j^*X_j$, although seemingly restrictive, contains in fact a general family of intrinsically defined sub-Riemannian Laplacians.

We first define here the sub-Riemannian Laplacian $\Delta_{(U,f), ds}$ associated to a sub-Riemannian structure $(U,f)$ on $\M$ and a smooth density $ds$. We then explain why it can be rewritten under the form \eqref{def:L} for some (sufficiently many) vector fields $X_1, \cdots , X_m$. 

First, we assume that $\M$ is equipped with a general sub-Riemannian structure $(U,f)$, see~\cite[Definition~1.3]{Bellaiche} or~\cite[Definition~3.2]{ABB:EMS16} with $U$ a Euclidean bundle with base $\M$ and $f:U\to T\M$ a smooth map being linear on fibers. This allows to define first a sub-Riemannian metric, that is, a metric on the horizontal distribution $\mathcal{D}$ with $\mathcal{D}_x = f(U_x) \subset T_x \M$ by $g(x,v) = \inf\{|u| , u \in U_x , v = f(x, u)\}$ (where $|\cdot |$ denotes the Euclidean norm in $U$). Second, this provides a sub-Riemannian gradient $\nabla_{(U,f)}$ on $\M$: namely, for $u\in C^\infty(\M)$, $\nabla_{(U,f)}u(x)$ is the unique vector in $\mathcal{D}_x$ such that for all $v \in\mathcal{D}_x$, we have $d_x u (v) = \tilde{g}(x, \nabla_{(U,f)}u(x), v)$ (where $\tilde{g}$ is the bilinear form associated to $g$).

Next, the smooth density $ds$ allows to define the divergence $\div_{ds}$ of a vector field $X \in \mathcal{X}^{\infty}$ by
$$
\left.\frac{d}{dt}(e^{tX})^*(ds)\right|_{t=0} = \div_{ds} (X) ds ,
$$
where $e^{tX}$ denotes the flow of $X$ (or, equivalently, by the formula $X^* = -X -\div_{ds}(X)$).
Hence, a natural definition of the sub-Riemannian Laplacian $\Delta_{(U,f), ds}$ is 
$$
\Delta_{(U,f), ds} u = \div_{ds} \left( \nabla_{(U,f)} u\right) , \quad u \in C^\infty(\M).
$$

Now, according to~\cite[Corollary~3.26]{ABB:EMS16}, the sub-Riemannian structure $(U,f)$ is equivalent to a free one, that is, there exist $m \in \N$ and $m$ vector fields $X_1, \cdots , X_m$  on $\M$ such that the horizontal distribution at $x \in \M$ is given by $\mathcal{D}_x = \vect(X_1(x), \cdots , X_m(x))$, and the metric on this distribution is defined by \eqref{def-g}.
A computation similar to that in Appendix~\ref{app:davide} shows that the sub-Riemannian gradient $\nabla_{(U,f)}$ of a function $u$ is then given by 
$$
\nabla_{(U,f)} u = \sum_{i=1}^m (X_i u) X_i .
$$
Hence, the formula $\div_{ds}(uX) = u \div_{ds}(X)+ Xu$ for $u \in C^\infty(\M)$ and $X \in \mathcal{X}^{\infty}$ yields
$$
\Delta_{(U,f), ds} u = \div_{ds} \left(\sum_{i=1}^m (X_i u) X_i \right) = \sum_{i=1}^m \div_{ds}( X_i) X_i u   +  X_i^2 u
 = - \sum_{i=1}^m X_i^* X_i u .
$$
As a consequence, all results presented in this article remain valid for general, intrinsically defined sub-Riemannian Laplacians $\Delta_{(U,f), ds}$.

In the above discussion, we assume the density $ds$ to be given: the sub-Riemannian Laplacian $\Delta_{(U,f), ds}$ then depends both on the sub-Riemannian structure $(U,f)$ and the density.
One may also wonder whether, given the sub-Riemannian structure $(U,f)$ only, there is an associated intrinsic choice of density $ds$, as in the Riemannian case. This question is an object of current research. While it seems that in the equiregular case, i.e., when the growth vector does not depends on the point, there is a natural intrisic measure (namely the Popp measure, see e.g.~\cite{Montgomery:book}), there is no consensus for what should be the natural one in the general case. For instance, in the Grushin case of Example \ref{ex:Grushin} the metric $g$ defined in~\eqref{def-g} is Riemannian outside of $\{x_1=0\}$. Hence, 
a natural choice would be to take the Riemannian density outside of $\{x_1=0\}$. The associated Laplacian is equal to that of Example \ref{ex:Grushin} plus a potential which is singular on $\{x_1=0\}$. The analysis in ~\cite{BoscLaurent:13} shows that the zone $\{x_1=0\}$ creates a barrier which the information cannot cross. This is in strong contrast with results of the present paper.

To conclude this remark, let us also notice that the whole class of operators studied by Fefferman and Phong in~\cite{FP:83} is not contained in the class of sub-Riemannian Laplacians $\Delta_{(U,f), ds}$ defined above. It would be interesting to investigate the questions of the present paper for such hypoelliptic operators. 
\end{remark}

\begin{remark}[Norm of the observation term]
Note that in the right hand-side of~\eqref{th-estimate-k}, the observation term only comes with a $L^2$ norm (which is not the case in most results in~\cite{LL:15}). This is due to the fact that we are in the context of operators with analytic coefficients with respect to all variables. In the case of partially analytic operators as described in Subsection \ref{s:without-analyticity}, we are able to get observability in $L^2$ using a  refined argument (see Section \ref{subsectL2partial}). Similar arguments are also applied to the classical wave equation in a forthcoming companion paper \cite{LL:17}.
\end{remark}

\begin{remark}[Other levels of $\H^{s}_{\L}$ regularity]
\label{rem:Hsspaces}
As already mentionned, Theorem~\ref{thmwavehypo} is a particular case of general estimates where all Sobolev scales are possible for measuring the typical frequency of the initial datum. Indeed, we prove the following more general result.
\begin{theorem}
\label{thmwavehypo-s}
Let $\L$ as above satisfying Assumptions \ref{assumLiek} and \ref{hypoanal}.
Assume that $\omega$ is a non empty open set of $\M$ and $T>\sup_{x\in \M} d_\L(x,\omega)$.
Then, for any $s>0$, there exist $\kappa , C ,\mu_0$ so that for all $\mu\geq \mu_0$, and all $u$ solution of \eqref{hypoelliptic-wave}, we have
\bnan
\label{e:wavehypo-s}
\nor{(u_0,u_1)}{L^2\times \H^{-1}_{\L}} \leq C e^{\kappa \mu^{k}}\nor{u}{L^2(]-T,T[\times \omega)} +\frac{1}{\mu^{s}}\nor{(u_0,u_1)}{\H^{s}_{\L}\times \H^{s-1}_{\L}}
\enan
\end{theorem}

Theorem~\ref{thmwavehypo} is the case $s=1$ of Theorem~\ref{thmwavehypo-s} and we believed that the frequency functions $\Lambda$ used in Theorem \ref{thmwavehypo} is the more natural presentation. Yet, it turns out that Theorem~\ref{thmwavehypo-s} is actually easier to prove for $s=k$. Below, we first prove this simpler case $s=k$ (Section \ref{s:simple-case}); then we need to prove refined estimates for the general case (Section \ref{s:energy-general-case}). The second part of the proof seems to require additional arguments in the partially analytic case described in Section~\ref{s:without-analyticity}. That is the reason why the estimate \eqref{th-estimate-partial} of Theorem \ref{t:partially-anal} is restricted to the case $s=k$. Nevertheless, as already explained, most of the results about eigenfunction and the hypoelliptic heat equation only use the easier case $s=k$.

Remark also that, in inequalities such as~\eqref{e:wavehypo-s}, deducing the $\H^{s'}_{\L}$ case from the $\H^{s}_{\L}$ case follows from an interpolation argument if $s'>s$. Indeed, for $0\leq s < s'$, denoting $\widetilde{\H}^s = \H^{s}_{\L}\times \H^{s-1}_{\L}$ we have
$$
\nor{U}{\widetilde{\H}^s} \leq \nor{U}{\widetilde{\H}^0}^{1-\frac{s}{s'}} \nor{U}{\widetilde{\H}^{s'}}^{\frac{s}{s'}}
\leq (1-\frac{s}{s'}) A \nor{U}{\widetilde{\H}^0} + \frac{s}{s'} A^{1-\frac{s'}{s}} \nor{U}{\widetilde{\H}^{s'}} , \quad \text{ for all }A >0 .
$$
Taking then $A = \frac12 \mu^s$ and putting this into~\eqref{e:wavehypo-s} yields the same estimate with $s$ replaced by $s'$. Hence, the difficulty in Theorem~\ref{thmwavehypo-s} when compared to the case $s=k$ (which proof is simpler) is only for small $s>0$.

Finally, let us mention that the result of Theorem~\ref{t:approx-control-heat} remains valid as well with the $\H^1_\L$-norm replaced by any $\H^s_\L$-norm, $s>0$, when modifying the powers accordingly.
\end{remark}

\begin{remark}[Constants]
Note that there are mostly two relevant constants in Theorem~\ref{thmwavehypo}, namely the minimal time $2\sup_{x\in \M} d_\L(x,\omega)$, and the constant $\kappa$ in the exponent of~\ref{th-estimate-k}. All constants appearing in the results of Section~\ref{s:approx-control-heat} can be explicitely formulated in terms of these two. For instance, in Theorem~\ref{thm:parabolic}, the constant $T_0$ can be taken as $T_0 >\kappa$, where $\kappa$ is the exponent in~\eqref{th-estimate-k} for some $S> \sup_{x\in \M} d_\L(x,\omega)$.
We refer to Remark~\ref{rem-ctes} below for more on this subject.
\end{remark}

\begin{remark}[Interpolation spaces, see~\cite{Leb:Analytic,LL:17}]
Notice that Estimate~\eqref{e:heat-approx-eps} may be reformulated (after an optimization in $\eps$) as the following interpolation inequality, for $\y$ solution to~\eqref{abstractheat} 
\bna
\nor{\y(T)}{L^2}
\leq C \nor{\y}{L^2((0,T)\times \omega)}^{\frac{T-(T_0+\eta)}{T-\eta}}\nor{\y(0)}{L^2}^{\frac{T_0}{T-\eta}},
\ena
while \eqref{estimobserheathypohigh} can be written
\bna
\nor{\y(0)}{L^2}^2\leq C \nor{\y}{L^2((0,T)\times \omega)}^{\frac{\theta-\theta_{0}}{\theta}}\nor{\y(0)}{k/2,\theta}^{\frac{\theta_0}{\theta}}.
\ena
More generally, in both cases, there exists $\alpha\in ]0,1[$ so that, we have the estimates
\bnan
\label{e:interp-ineq-a-la-lebeau}
\nor{\y_0}{F_1}
\leq C  \nor{\y_0}{F_{obs}}^{\alpha}\nor{\y_0}{F_0}^{1-\alpha} ,
\enan
where $F_0$, $F_T$ and $F_{obs}$ which are defined as the spaces of data obtained as the completion of linear combinations of eigenfunctions of $\L$ for the norms
\begin{equation}
\begin{array}{rcll}
\nor{\y_0}{F_1}&=&\nor{e^{-T\L}\y_0}{L^2}, &\quad \textnormal{ resp. }\nor{\y_0}{F_1} =\nor{\y_0}{L^2} , \\
\nor{\y_0}{F_0}&=&\nor{\y_0}{L^2} ,& \quad \textnormal{ resp. }\nor{\y_0}{F_0} = \nor{\y_0}{k/2,\theta} , \\
\nor{\y_0}{F_{obs}}&=&\nor{\y}{L^2((0,T)\times \omega)} = \nor{e^{-t\L}\y_0}{L^2((0,T)\times \omega)} &\quad  \text{ in both cases}.
\end{array}
\end{equation}
The latter are proper norms as a consequence of uniqueness, backward uniqueness (consequence e.g. of Lemma~\ref{l:jensen-norm} below) and unique continuation property for the hypoelliptic heat equation~\eqref{abstractheat}. Note that we have $F_0 \subset F_{obs}, F_1$  in the first case, and $F_0\subset F_1 \subset F_{obs}$ in the second.

As a consequence of~\eqref{e:interp-ineq-a-la-lebeau} (see for instance \cite[Appendix, Lemma 1]{Leb:Analytic}), there exists $\delta>0$ such that 
\bna
[F_0,F_{obs}]_{\delta}\subset F_1 , 
\ena where $[F_0,F_{obs}]_{\delta}$ is the space of interpolation between $F_0$ and $F_{obs}$. As in Lebeau \cite[Section 3]{Leb:Analytic}, this yields
\bnan
\label{e:lebeau-dual}
 F_1' \subset[F_0',F_{obs}']_{1-\delta}. 
\enan
Now, the duality between~\eqref{abstractheat} and~\eqref{e:control-heat} will allow to identify the spaces $F_0', F_1' , F_{obs}'$, to deduce properties of the controllable and the attainable sets.

\medskip
First, the duality between~\eqref{abstractheat} and~\eqref{e:control-heat} writes 
\bnan
\label{e:duality}
\int_0^T(\mathds{1}_{\omega}\y(T-t) , g )_{L^2(\M)} dt = (\y_0 , u(T))_{L^2(\M)} - (\y(T), u_0)_{L^2(\M)}.
\enan
We define $E_{att} = \{u_1 \in L^2(\M), \text{there exists } g \in L^2((0,T)\times \omega) , \text{ s.t. the associated solution $u$ to \eqref{e:control-heat} with } u(0)=0 \text{ satisfies }u(T)=u_1 \}$ the space of attainable data from zero with $L^2$ control, endowed with the norm
\bna
\nor{u_1}{E_{att}}=\inf \left\{\nor{g}{L^2((0,T)\times \omega)} , g \in L^2((0,T)\times \omega) \text{ s.t. the associated solution $u$ to \eqref{e:control-heat} with } u(0)=0 \text{ satisfies }u(T)=u_1\right\} .
\ena
From~\eqref{e:duality}, we obtain for all $u_1 \in E_{att} \subset L^2$ and all $\y_0 \in L^2$
$$ 
\left| (\y_0 , u_1)_{L^2(\M)} \right| = \left| \int_0^T(\mathds{1}_{\omega}\y(T-t) , g )_{L^2(\M)} dt  \right| \leq \|\y_0\|_{F_{obs}} \|u_1\|_{E_{att}}. 
$$
Hence, the $L^2(\M)$ scalar product extends uniquely as a duality product $\langle \y_0 , u_1 \rangle_{F_{obs} , E_{att}}$, allowing to identify $F_{obs}'$ with $E_{att}$. With the identification, we have as well $F_1= e^{T\L}L^2(\M)$ and $F_0 = L^2(\M)$ (resp. $F_1 = L^2(\M)$ and $F_0 =H^{k/2, \theta}$) so that $F_1' \approx e^{-T\L}L^2(\M)$ and $F_0' \approx L^2(\M)$ (resp. $F_1' \approx L^2(\M)$ and $F_0' \approx H^{k/2, -\theta}$). With~\eqref{e:lebeau-dual}, this yields
\bna
e^{-T \L}L^2\subset [E_{att},L^2]_{1-\delta} , \\
\text{resp.} \quad L^2\subset [E_{att},H^{k/2,-\theta}]_{1-\delta} .
\ena
We also define $E_{cont}$ the (abstract) space of data that can be controled towards zero as the completion of the space 
$$\{u_0 \in L^2, \text{there exists } g \in L^2((0,T)\times \omega) , \text{ s.t. the associated solution $u$ to \eqref{e:control-heat} satisfies } u(T)=0 \}$$
 for the norm
\bna
\nor{u_0}{E_{cont}}=\inf \left\{\nor{g}{L^2((0,T)\times \omega)} , g \in L^2((0,T)\times \omega) \text{ s.t. the associated solution $u$ to \eqref{e:control-heat} satisfies } u(T)=0 \right\} .
\ena
From~\eqref{e:duality}, we obtain for all $u_1 \in E_{cont} \subset L^2$ and all $\y_0 \in L^2$
$$ 
\left| (\y(T) , u_1)_{L^2(\M)} \right| = \left| \int_0^T(\mathds{1}_{\omega}\y(T-t) , g )_{L^2(\M)} dt  \right| \leq \|\y_0\|_{F_{obs}} \|u_1\|_{E_{cont}}. 
$$
Similarly, the scalar product $\left\langle \y_0 , u_1\right\rangle= (e^{-T\L}\y_0 , u_1)_{L^2(\M)}$ extends uniquely as a duality product $\langle \y_0 , u_1 \rangle_{F_{obs} , E_{att}}$, allowing to identify $F_{obs}'$ with $E_{cont}$. With this same identification, we have as well $F_1' \approx L^2(\M)$ and $F_0' \approx e^{-T\L}L^2(\M)$ (resp. $F_1' \approx e^{-T\L}L^2(\M)$ and $F_0' \approx e^{-T\L}H^{k/2, -\theta}$). 
Similarly, \eqref{e:lebeau-dual} also yields 
\bna
L^2\subset [E_{cont},e^{-T \L}L^2]_{1-\delta}\\
e^{-T\L}L^2\subset [E_{cont},e^{-T\L}H^{k/2,-\theta}]_{1-\delta}. 
\ena
Note that when $k>2$, $e^{-T\L}H^{k/2,-\theta}\approx H^{k/2,-\theta}$ which is not a distributional set whatever the time $T$ is. Yet, if $k=2$, $e^{-T\L}H^{k/2,-\theta}\approx H^{1,-\theta+T}$, which is a space of very regular functions if $T>\theta$.
\end{remark}

The next remark concerns the results of Theorem~\ref{thm:parabolic} and Corollary~\ref{cor:parabolic}.
\begin{remark}[Large time approximate controllability with polynomial cost of ``critical anomalous diffusion'']
\label{rkdemichaleur}
The proof of Theorem~\ref{thm:parabolic} and Corollary~\ref{cor:parabolic} also applies to any positive selfadjoint operator satisfying spectral estimates (or similar estimates for the control of the heat equation for spectrally localized initial data) like
\bnan
\label{spectralhalf}
\nor{w}{L^2(\M)}\leq Ce^{c\lambda}\nor{w}{L^2(\omega)},\quad \text{for all} \quad w=\sum_{\lambda_j\leq \lambda}w_j \varphi_j.
\enan
This is in particular the case for the square root of the Laplacian $\sqrt{-\Delta_g}$ where $\Delta_g$ is the elliptic Laplace-Beltrami operator on a compact Riemannian manifold (even without the analyticity assumption and with Dirichlet boundary condition), see \cite{LR:95}. Therefore, for the associated evolution operator (so called ``critical anomalous diffusion'') $\d_t + \sqrt{-\Delta_g}$ (studied in~\cite{Miller:06}), the same approximate controllability result with a polynomial cost holds. Note also that even for the one dimensional case (namely the operator $\d_t + |\d_x|$ on the circle), it has been proved by Koenig \cite{Koe:17} that exact controllability in finite time $T>0$ never holds (as long as the control domain is not the whole circle). In particular, it suggests that approximate controllability at polynomial cost might be the best to obtain under general spectral assumptions like \eqref{spectralhalf}.
\end{remark}

\bigskip
\noindent
{\em Acknowledgements.} 
We wish to thank warmly Davide Barilari and Emmanuel Tr\'elat for very helpful discussions on sub-Riemannian geometry. In particular, Emmanuel Tr\'elat brought our attention to the references~\cite{RiffordTrelatMorse} and \cite{Derridj:71}, and the properties of Example~\ref{ex:lie-groups}, Remark~\ref{r:sub-riem-lap} and Appendix~\ref{app:davide} were explained to us by Davide Barilari.
The first author is partially supported by the Agence Nationale de la Recherche under grant  SRGI ANR-15-CE40-0018 and IPROBLEMS ANR-13-JS01-0006.
The second author is partially supported by the Agence Nationale de la Recherche under grant GERASIC ANR-13-BS01-0007-01.

\section{The quantitative Holmgren-John theorem of~\cite{LL:15}}
\label{s:Holmgren-john-LL15}

In this section, we briefly review some results obtained in~\cite{LL:15}, that will be at the core of the proof of the present paper. We shall only consider a very particular class of operators, namely second order operators with real principal symbol. Also, we shall only consider non-characteristic surfaces. This assumption can be also relaxed (see e.g.~\cite{LL:15} Definition~1.7 and Remark~1.9), even though we are not aware of any application of the refined result.

The interest of taking such operators and surfaces is that, in this context, several assumptions and formulations of the results in~\cite{LL:15} are simplified.
 In this section, we state results for an operator $P$ in $\R^n$, where, in the application in Section~\ref{s:hypo-wave} below, we shall mainly consider $P$ as a local version of $\d_t^2 + \L$ on $\R \times \R^d$ (recall that $\dim(\M)= d$), that is $n=d+1$.
For this reason (and as opposed to the notation of~\cite{LL:15}), we shall denote by $\z \in \R^n$ the running variable and $\zeta \in \R^n$ its cotangent variable. In the applications in the next section, we will have $\z=(t,x)$ and $\zeta = (\tau, \xi)$.

\subsection{A typical quantitative unique continuation result of~\cite{LL:15}}
\label{s:typical-result}

 A typical instance (in the situation describe above) of the main result of~\cite{LL:15} may be stated as follows (see~\cite[Theorem~1.11]{LL:15} together with~\cite[Remark~1.10]{LL:15}).

\bigskip

\begin{figure}[h!]
  \begin{center}
    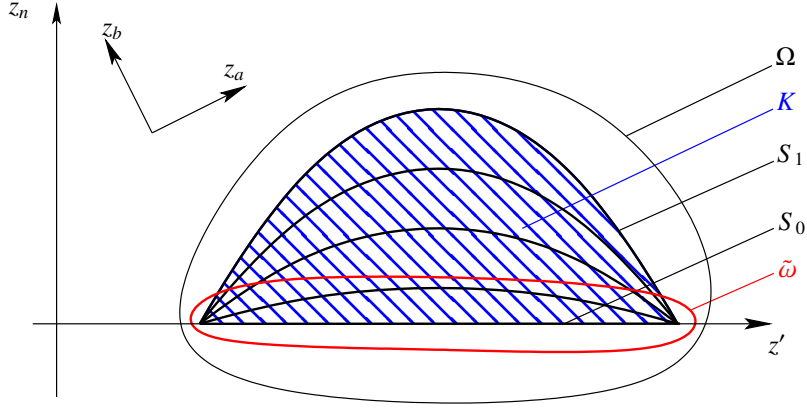 
    \caption{Geometric setting of Theorem~\ref{thmsemiglobal}}
    \label{f:geom-setting}
 \end{center}
\end{figure}

\textbf{Geometric setting:} (see Figure~\ref{f:geom-setting})
We first fix two splittings of $\R^n$:
\begin{itemize}
\item $\R^n=\R^{n_a} \times \R^{n_b}$, where $n_a + n_b= n$. We denote $\z=(\z_a,\z_b)$ the global variable and $\zeta=(\zeta_a,\zeta_b)$ the associated cotangent variable.
\item and $\R^n=\R_{\z'}^{n-1}\times \R_{\z_n}$,
\end{itemize} possibly in two different bases.
We let $D$ be a bounded open subset of $\R^{n-1}$ with smooth boundary and $G = G(z',\eps)$ a $C^2$ function defined in a neighborhood of $\overline{D} \times [0,1]$, such that  
\begin{enumerate}
\item \label{geom-set-1} For all $\eps \in (0,1]$, we have $\{\z' \in \R^{n-1} , G(\z' , \eps) \geq 0\} = \overline{D}$;
\item \label{geom-set-2} for all $\z' \in D$, the function $\eps \mapsto G(\z' , \eps)$ is strictly increasing;
\item \label{geom-set-3} for all $\eps \in (0,1]$, we have $\{\z' \in \R^{n-1} , G(\z' , \eps) = 0\} = \partial D$.
\end{enumerate}
We set $G(\z', 0) = 0$, $S_0 = \overline{D} \times \{0\}$ and,  for $\eps \in (0,1]$, 
\bna
&S_\eps = \{(\z',\z_n) \in \R^{n} , \z_n \geq 0 \text{ and } G(\z' , \eps) =\z_n \} =( \overline{D} \times \R )\cap\{(\z',\z_n) \in \R^{n} , G(\z' , \eps) =\z_n \} ;\\
&K = \{\z \in \R^n ,0\leq  \z_n \leq G(\z', 1)\}.
\ena

We recall that the local surface $S := \{\varphi=0\} \ni \z_0$, $d\varphi(\z_0)\neq 0$ is called {\em non-characteristic} at $\z_0$ for the differential operator $P$ with principal symbol $p$ if $p(\z_0 , d \varphi(\z_0)) \neq 0$, and that this is a property of the sole surface $S$ (together with the point $\z_0$ and the principal symbol of the operator $p$) and not its defining function $\varphi$.

\medskip
Note also that in the main part of the paper, the operators are analytic with respect to all variables, 
In this case, Theorem~\ref{thmsemiglobal} is a quantitative version of the Holmgren-John theorem (for second order operators, see~\cite{LL:15} in the general case), and may be seen as a generalization of~\cite{Leb:Analytic}, which concern the (analytic) wave operator.

\begin{theorem}
\label{thmsemiglobal}
In the above geometric setting, we moreover let $\Omega \subset \R^{n_a} \times \R^{n_b}$ be a bounded open neighborhood of $K$, and $P$ be a differential operator of order $2$ on $\Omega$ such that
\begin{itemize}
\item all coefficients of $P$ are smooth and depend analytically on the variable $\z_a$, 
\item the principal symbol of $P$, namely $p(\z,\zeta) = Q_\z(\zeta)$, is a $\z$-family of real quadratic forms such that $\zeta_b \mapsto Q_\z(0, \zeta_b)$ is definite on $\R^{n_b}$ for any $\z \in \Omega$. 
\end{itemize}

Assume also that, for any $\e\in [0,1+\eta)$, $\eta>0$, the surface $S_{\e}$ is {\em non-characteristic} for $P$ at each point of $S_\eps$.

Then, for any open neighborhood $\tilde{\omega} \subset \Omega$ of $S_0$, there exists a neighborhood $U$ of $K$, and constants $\kappa ,C ,\mu_0 >0$ such that for all $\mu\geq \mu_0$ and $u\in C^{\infty}_0(\R^n)$,  we have
\bnan
\label{e:estimate-hyp}
\nor{ u}{L^2(U)}\leq C e^{\kappa \mu}\left(\nor{ u}{H^{1}_b(\tilde{\omega})} + \nor{Pu}{L^2(\Omega)}\right)+\frac{C}{\mu}\nor{u}{H^{1}(\Omega)} ,
\enan
where we have denoted $\nor{u}{H^{1}_b(\tilde{\omega})}=\sum_{|\beta| \leq 1}\nor{ D_b^\beta u}{L^2(\tilde\omega)}$.
\end{theorem}

Unfortunately, in the present paper, this global result does not apply under this form. In order to state the refined (and more technical) version, used in the main part of the paper, we shall need some definitions taken from~\cite[Section~2.3]{LL:15}.

\subsection{Definitions and tools for propagating the information}
\label{s:def-mult}
We first define the following regularization process for functions $f$ defined on $\R^n$, by $f\mapsto f_\lambda$ with
$$
f_{\lambda}:=e^{-\frac{|D_a|^2}{\lambda}} f =  \mathcal{F}_a^{-1} \left( e^{-\frac{|\zeta_a|^2}{\lambda}} \mathcal{F}_a (f)(\zeta_a , \z_b)\right)(\z_a),
$$
where $\mathcal{F}_a$ denotes the Fourier transform in the variable $z_a$ only, 
or, equivalently
\bna
f_{\lambda}(\z_a , \z_b)=\left(\frac{\lambda}{4\pi}\right)^{\frac{n_a}{2}}\Big( e^{-\frac{\lambda}{4}|\cdot|^2} *_{\R^{n_a}} f(\cdot , \z_b) \Big)(\z_a)
 = \left(\frac{\lambda}{4\pi}\right)^{\frac{n_a}{2}}\int_{\R^{n_a}} f\left(y_a, \z_b\right)e^{-\frac{\lambda}{4} |\z_a- y_a|^2}~dy_a .
\ena
Then, we also need to introduce frequency localization functions, i.e. appropriately smoothed Fourier multipliers. 
Let $m(\zeta_a)$ be a smooth radial function (i.e. depending only on $|\zeta_a|$), compactly supported in $|\zeta_a| < 1 $ such that $m(\zeta_a)=1$ for $|\zeta_a| < 3/4$. We denote by $M^{\mu}$ the Fourier multiplier $M^{\mu}u= m\left(\frac{D_a}{\mu}\right)u$, that is 
$$
(M^{\mu} u) (\z_a , \z_b)=  \mathcal{F}_a^{-1} \left( m \left(\frac{\zeta_a}{\mu}\right) \mathcal{F}_a (u)(\zeta_a , \z_b)\right)(\z_a) .
$$
Given $\lambda , \mu >0$, we shall denote by $M^{\mu}_{\lambda}$  the Fourier multiplier of symbol $m^{\mu}_{\lambda}(\zeta_a)=m_{\lambda}\left(\frac{\zeta_a}{\mu}\right)$, i.e. $M^{\mu}_{\lambda}= m^{\mu}_{\lambda}(D_a) = m_{\lambda}\left(\frac{D_a}{\mu}\right)$ or
$$
(M^{\mu}_{\lambda} u) (\z_a , \z_b)=  \mathcal{F}_a^{-1} \left( m_{\lambda}\left(\frac{\zeta_a}{\mu}\right) \mathcal{F}_a (u)(\zeta_a , \z_b)\right)(\z_a) ,
$$
with, according to the above notation for the subscript $\lambda$, 
\bnan
\label {e:reg-mult}
m_{\lambda}(\zeta_a) 
 = \left(\frac{\lambda}{4\pi}\right)^{\frac{n_a}{2}}\int_{\R^{n_a}} m \left(\eta_a \right)e^{-\frac{\lambda}{4} |\zeta_a- \eta_a|^2}~d\eta_a .
\enan

Note that in this definition, the symbol is first regularized and then dilated. 
We stress the fact that these Fourier multipliers only act in the variable $\z_a$.

\bigskip
The typical local estimate of \cite[Theorem~3.1]{LL:15}, which is the building block for semiglobal statements like that of Theorem~\ref{thmsemiglobal}, reads as
\bna
\nor{M^{\beta\mu}_{\mu} \sigma_{\mu} u}{1}\leq C e^{\kappa \mu}\left(\nor{M^{\alpha \mu}_{\mu} \vartheta_{\mu} u}{1} + \nor{Pu}{L^2(B(0,R))}\right)+Ce^{-\kappa' \mu}\nor{u}{1} ,
\ena
for all $\mu\geq \mu_0$ and $u\in C^{\infty}_0(\R^n)$, where $\sigma$ is a cutoff function in a small ball $B(0,r)$, $r<R$, whereas $\vartheta$ is a cutoff in only one side (the one where the information is taken) of the hypersurface passing through zero (and being non-characteristic).

Here, and below, the norm $\nor{\cdot}{1}$ is the norm $\nor{\cdot}{H^1(\R^n)}$.

Such an estimate only provides information on the low frequency part of the function, through the frequency cutoff $M^{\beta\mu}_{\mu}$, with an exponentially small $Ce^{-\kappa' \mu}$ remainder (as opposed to the $1/\mu$ remainder term in~\eqref{e:estimate-hyp}).
 Iterating this result allows us to propagate the low frequency information. In this section, we recall some tools, used in~\cite[Section~4]{LL:15}, for this iterative procedure. 
They are aimed at describing how information on the low frequency part of the solution can be deduced from one subregion to another one.

The following is~\cite[Definition~4.4]{LL:15}, given here in the context of second order operators.
\begin{definition}
\label{defdependencestrong}
Fix $\Omega $ be an open set of $\R^n=\R^{n_a}\times \R^{n_b}$ and $P$ a differential operator of order $2$ defined in $\Omega$, and $(V_j)_{j\in J}$ and $(U_i)_{i\in I}$ two finite collections of bounded open sets of $\R^n$.
We say that $(V_j)_{j\in J}$ is \textbf{under the strong dependence} of $(U_i)_{i\in I}$, denoted by 
\bna
(V_j)_{j\in J} \lhd (U_i)_{i\in I},
\ena
 if there exists $W_i\Subset U_i$ such that 
 for any $\vartheta_i\in C^{\infty}_0(\R^n)$ such that $\vartheta_i(\z)=1$ on a neighborhood of $\overline{W_i}$, for any $\widetilde{\vartheta}_j\in C^{\infty}_0(V_j)$
and for all $\kappa, \alpha>0$, there exist $C, \kappa', \beta,\mu_0 >0$ such that for all $(\mu,u) \in [\mu_0, + \infty)\times C^{\infty}_0(\R^{n})$, we have 
\bna
\sum_{j\in J}\nor{M^{\beta\mu}_{\mu} \widetilde{\vartheta}_{j,\mu}  u}{1}\leq C e^{\kappa \mu}\left(\sum_{i\in I} \nor{M^{\alpha\mu}_{\mu} \vartheta_{i,\mu} u}{1} + \nor{Pu}{L^2(\Omega)}\right)+Ce^{-\kappa' \mu}\nor{u}{1} .
\ena
If the cardinal of $I$ is one, writing $U$ the single set of the family $(U_i)_{i\in I}$, we simply denote $(V_j)_{j\in J} \unlhd U$. We use the same convention for $V$ in case the cardinal of $J$ is one.
The norm $\nor{\cdot}{1}$ is taken in $\R^n$.
\end{definition}

We summarize the properties of this relation in the following proposition~\cite[Proposition~4.5]{LL:15}.
\begin{proposition}
\label{propstrong}
We have the following properties
\begin{enumerate}
\item \label{proptotostrong} If $(V_j)_{j\in J} \lhd (U_i)_{i\in I}$ with $U_i=U$ for all $i\in I$, then $(V_j)_{j\in J} \lhd U$.
\item \label{propincludestrong}If $V_i\Subset U_i$ for any $i\in I$, then, $(V_i)_{i\in I} \lhd  (U_i)_{i\in I}$.
\item \label{propunionstrong} If $V_i\Subset U_i$ for any $i\in I$, then $\bigcup_{i\in I} V_i \lhd (U_i)_{i\in I}$.
\item \label{propproduitstrong} If for any $i\in I$, $V_i\lhd U_i$, then $(V_i)_{i\in I}\lhd (U_i)_{i\in I}$. In particular, if for any $i\in I$, $U_i\lhd U$, then $(U_i)_{i\in I}\lhd U$. 
\item \label{proptransstrong} The relation $\lhd$ is transitive, that is \bna
\left[(V_j)_{j\in J} \lhd (U_i)_{i\in I} \textnormal{ and } (U_i)_{i\in I}\lhd (W_k)_{k\in K } \right]\Longrightarrow (V_j)_{j\in J} \lhd (W_k)_{k\in K }.
\ena
\end{enumerate}
\end{proposition}
Note that we do not always have $U \lhd U$.
\begin{remark}
\label{rkinvariance}
We stress the fact that the definition of $\lhd$ actually depends on the set $\Omega$, the splitting $\R^n=\R^{n_a}\times \R^{n_b}$ and the operator $P$.
The dependence of $\lhd$ upon these objects will be mentioned when needed.
 For the applications, it is important that the function $u$ is not necessarily supported in $\Omega$. 

In the following, we will only need to use this relation $\lhd$ in some appropriate coordinate charts. However, it will not be a problem for what we want to prove, even on a compact manifold. Indeed, we will fix some coordinate chart on an open set $\Omega\subset \R^n$ close to a point or close to a trajectory. Then, we will use the relation $\lhd$ related to $\Omega$ to finally obtain some estimates which will be invariant by changes of coordinates. 
\end{remark}

We will also use the following proposition, \cite[Proposition~4.9]{LL:15}, which allows to iterate local propagation results towards global ones.
\begin{proposition}
\label{propiterationabstrait}
Assume that there exists some open sets $U_0$, $U_{i,j}$, $\omega_{i,j}$, $V_{i,j}$, with $j\in \llbracket 1,N\rrbracket$ and $i\in I_j$ ($I_j$ finite) such that we have
\bna
& U_{i,j}\lhd V_{i,j} \quad \text{ and } \quad \omega_{i,j}\Subset U_{i,j}, \quad \text{ for all } j\in \llbracket 1,N\rrbracket \text{ and }i\in I_j ; & \\ 
& V_{m,l+1} \Subset \left[U_0 \cup\bigcup_{j\in \llbracket 1,l\rrbracket}\bigcup_{i\in I_{j}} \omega_{i,j}\right], \quad \text{ for all } l\in \llbracket 0,N -1\rrbracket,  \text{ and }m\in I_{l+1} ,&
\ena 
where we consider the union $\bigcup_{j\in \llbracket 1,l\rrbracket}$ empty if $l=0$.
Then, we have $\left[U_0 \cup\bigcup_{j\in \llbracket 1,N\rrbracket}\bigcup_{i\in I_{j}} \omega_{i,j}\right]\lhd V_0$ for any open set $V_0$ such that $U_0\Subset V_0$.
\end{proposition} 
In this proposition, the local propagation results is $U_{i,j}\lhd V_{i,j}$ but the iteration is made by packets. Roughly speaking, we use all sets corresponding to indices $i\in I_j$, $j\leq l$ to deduce the information on the sets with indices $i\in I_{l+1}$.
\subsection{Semiglobal estimates along foliation by hypersurfaces}
\label{subsectsemiglobal}

Now, we formulate the results of~\cite{LL:15} in the form they will be used in the next section, which is different from Theorem~\ref{thmsemiglobal} with two respects:
\begin{itemize}
\item First, we keep the formulation with $\lhd$; this means that we keep a frequency cutoff in both handsides of the estimate, as well as an exponentially small remainder (it is a low frequency estimate only).  This allows to patch estimates together (which is no longer the case when the high frequencies have been taken into account, i.e. when estimates take the form of~\eqref{e:estimate-hyp}). The high frequencies will then be taken into account to close the estimates with two different methods in Section~\ref{s:energy-estimates}.
\item Second, we allow the linear change of variables between the two splittings (namely $(\z_a, \z_b)$ for the analytic dependence and $(\z', \z_n)$ for the geometry) to be replaced by a diffeomorphism, which shall be very useful in the following.
\end{itemize}

We give a first statement that is a low frequency formulation of Theorem~\ref{thmsemiglobal}, using the notation $\lhd$ (see~\cite[Theorem~4.7]{LL:15}).
\begin{theorem}
\label{thmsemiglobaldep}
Under the assumptions of Theorem \ref{thmsemiglobal}, there exists an open neighborhood $U$ of $K$ such that
\bna
U\lhd \widetilde{\omega}.
\ena
\end{theorem}
This essentially means that Estimate~\eqref{e:estimate-hyp} may be replaced by the following: for all $\chi\in C^{\infty}_0(U)$ and $\varphi \in C^{\infty}_0(\Omega)$ such that $\varphi=1$ on a neighborhood of $\widetilde{\omega}$, we have
\bnan
\label{estimaMmu}
\nor{M^{\beta\mu}_{\mu} \chi_{\mu}  u}{1}\leq C e^{\kappa \mu}\left(\nor{M^{\mu}_{\mu} \varphi_{\mu} u}{1} + \nor{Pu}{L^2(\Omega)}\right)+Ce^{-\kappa' \mu}\nor{u}{1} ,
\enan
 (for any $\kappa>0$, there exist $C, \beta, \kappa', \mu_0 >0$ such that for $\mu\geq \mu_0$) i.e. keep the frequency cutoff and the exponentially small remainder.

A remaining drawback of this statement, given by the geometric framework of Theorems~\ref{thmsemiglobal}, is that the hypersurfaces are described by graphs in some coordinates (namely $(\z',\z_n)$). This choice of description is mainly convenient to make the foliation more effective and order the hypersurfaces more easily, but is too rigid for the application in the present paper. 
Now, we give a slight variant of Theorem \ref{thmsemiglobaldep}, more adapted to possible changes of variables.
\begin{theorem}
\label{thmsemiglobaldepchgt}
Let $\Omega\subset \R^n = \R^{n_a} \times \R^{n_b}$ and $P$ be a differential operator of order $2$ on $\Omega$ such that
\begin{itemize}
\item all coefficients of $P$ are smooth and depend analytically on the variable $\z_a$, 
\item the principal symbol of $P$ namely $p(\z,\zeta) = Q_\z(\zeta)$ is a $\z$-family of real quadratic forms, such that $\zeta_b \mapsto Q_\z(0, \zeta_b)$ is definite on $\R^{n_b}$ for any $\z \in \Omega$. 
\end{itemize}
Let $\Phi$ be a diffeomorphism of class $C^2$ from $\Omega$ to $\widetilde{\Omega}=\Phi(\Omega)$. Assume that the Geometric Setting of Theorem~\ref{thmsemiglobal} is satisfied for some $D$, $G$, $K$, $S_{\e}$ on $\widetilde{\Omega}$ (and not on $\Omega$). 
Assume further that for any $\eps \in [0, 1+\eta)$, $\eta>0$, the surface $\Phi^{-1}(S_{\e})$ (well defined on $\Omega$) is {\em non-characteristic} for $P$ (at every point of $\Phi^{-1}(S_{\e})$).

Then, for all neighborhood  $\omega$ of $\Phi^{-1}(S_{0})$, there exists an open neighborhood $U\subset \Omega$ of $\Phi^{-1}(K)$ such that
\bna
U\lhd \omega ,
\ena
where $\lhd=\lhd_{\Omega,P}$ is related to the operator $P$ defined on $\Omega$ (see Remark \ref{rkinvariance}).
\end{theorem}
In this result, the operator $P$ has the appropriate form in $\Omega \subset \R^{n_a} \times \R^{n_b}$ whereas the geometry of the surfaces is defined in $\tilde{\Omega}$, both being linked by a diffeomorphism.

With this theorem in hand, we may now prove the results presented in Section~\ref{s:results}.

\section{The hypoelliptic wave equation, proof of Theorem~\ref{thmwavehypo}}
\label{s:hypo-wave}
The main goal of this section is to prove Theorem \ref{thmwavehypo}.
The proof is inspired by the case of the classical wave equation (see~\cite[Section~6.1]{LL:15}) with mainly two differences:
\begin{itemize}
\item We are not able to construct global coordinates near any path $\gamma$, as in the case of the wave equation. However, if $\gamma$ is a normal geodesic (see definition~\ref{def:normal} below), we are able to do this ocnstruction locally. Then, this local result needs to be iterated.
\item The $H^1$ norm is no longer equivalent to the energy norm for the hypoelliptic operator. We thus need to use hypoelliptic estimates instead.
\end{itemize}
Therefore, the proof is divided in two parts, the first of which concerning the geometric iteration process, and the second the energy estimates.

\bigskip
Let us start by introducing geometrical definitions and facts used all along the proofs.
First, denote by $\ell =\ell (x,\xi) \in C^\infty(T^*\M)$ the principal symbol of the operator $\L$, that is
\bnan
\label{e:symbol-p}
\ell (x, \xi) =   \sum_{i=1}^m \left< \xi , X_i(x) \right>^2 .
\enan
where $ \left< \xi , X(x)  \right>  = \left< \xi , X(x)  \right>_{T_x^*\M , T_x \M}$ is the duality bracket.

\begin{remark}
In view of unique continuation results, note that a local hypersurface $\{\varphi = 0\}$ at $x_0\in \M$ (where $\varphi : \M \to \R$ with $\varphi(x_0)= 0$ and $d\varphi(x_0) \neq 0$) is characteristic for the operator $\L$ if $\ell(x_0 , d\varphi(x_0))= 0$, that is, according to~\eqref{e:symbol-p}, if
$$
\left< d\varphi(x_0), X_i(x) \right> = 0 \quad  \text{ for all } i \in \{1 ,\cdots , m\} .
$$
\end{remark}

\begin{definition}
\label{def:Hamilton}
The Hamiltonian curve of the symbol $\ell$ issued from $\rho_0 \in T^*\M$ is the unique maximal solution $\rho(s)= (\gamma(s),\xi(s))$ of the ODE
\bnan
\label{e:Hamiltonian-flow}
 \dot{\rho}(s)=H_\ell ( \rho(s)) , \quad \rho(0) = \rho_0 ,
 \enan
 where $H_\ell$ is the Hamiltonian vector field associated to the Hamiltonian $\ell$. In local charts, this is
\bneq
\dot{\gamma}(s)&=& \frac{\partial}{\partial \xi}\ell (\gamma(s),\xi(s))  ,\\
\dot{\xi}(s)&=& - \frac{\partial}{\partial x}\ell (\gamma(s),\xi(s)) .
\eneq
\end{definition}
Such a curve is smooth (even real analytic since $\M$ and $\ell$ are).
Moreover, the first equation writes
\bnan
\label{e:hamilt-horizont}
\dot{\gamma} (s)=  \sum_{j=1}^m  2 \left< \xi (s), X_i(\gamma(s))  \right> X_i(\gamma(s)),
\enan
so that the projection on $\M$ of a Hamiltonian curve is a smooth horizontal curve (see Definition~\ref{d:horizontal-path}). Note also that a Hamiltonian curve, locally defined according to the Cauchy-Lipschitz theorem, is actually globally defined (see e.g.~\cite[Proposition~2.3.2]{Rifford:book}).
Note finally that the Hamiltonian $\ell$ is preserved along a Hamiltonian curve of $\ell$, i.e. $\ell(\rho(s)) = \ell(\rho_0)$ for every $\rho_0 \in T^*\M$, $s \in \R$, where $\rho$ is the solution to~\eqref{e:Hamiltonian-flow}.

Given a Hamiltonian curve $\rho = (\gamma , \xi) : [0, S]\to T^*\M$, one may compute the speed along the horizontal geodesic $\gamma : [0, S]\to \M$. Namely, we have (see Lemma~\ref{l:davide} for a proof of the first identity)
\bna
 g(\gamma(s) , \dot{\gamma}(s))    = 
   \sum_{j=1}^m  4 \left< \xi (s), X_i(\gamma(s))  \right>^2  
 = 4  \ell(   \gamma(s), \xi (s))  = 4    \ell(  \rho(0))  , \quad \text{for all }s \in [0,S].
\ena  
We can hence compute the length of the horizontal path $\gamma$, namely, 
\bna
\length(\gamma)  =  \int_0^{S}\sqrt{g(\gamma(s) , \dot{\gamma}(s))} ds  
 =2 S \sqrt{ \ell(  \rho(0))  } 
.
 \ena

This motivates the following definition:
\begin{definition}
\label{def:normal}
We say that a horizontal curve $\gamma : [0,L_0] \to \M$ is a normal geodesic if there exists $\xi(s)\in T^*_{\gamma(s)}\M$ such that $s\mapsto (\gamma(s), \xi(s))$ is a Hamiltonian curve of the symbol $\ell$ with $\ell(\gamma(s),\xi(s))=\frac14$.
\end{definition}
As a consequence of this definition, such a curve $\gamma : [0,L_0] \to \M$ has unit speed:
\bna
\sqrt{g(\gamma(s) , \dot{\gamma}(s))}   =  
2   \sqrt{ \ell(   \gamma(s), \xi (s)) } = 2  \sqrt{ \ell(  \rho(0))  } = 1, \quad \text{for all }s \in [0,L_0], 
\ena 
(it is hence parametrized by arclength) and length $L_0$.

\begin{definition}
\label{def:minimizing}
We say that a curve $\gamma : [0,L] \to \M$ is a minimizing geodesic path between $x$ and $y$ if $\gamma(0)=x$, $\gamma(L)=y$, if $\gamma$ is a horizontal curve and if we have $d_\L(x,y) = \length(\gamma)$ together with $g(\gamma(s), \dot{\gamma}(s))$ constant. 
\end{definition}

See~\cite[Chapter~2]{Rifford:book} or~\cite[Section~3.3]{ABB:EMS16}.
Note that all above definitions are intrinsic.
The following key result in our proofs is~\cite[Theorem~1.1]{RiffordTrelatMorse}. 

\begin{theorem}[Rifford-Tr\'elat~\cite{RiffordTrelatMorse}]
\label{th:riff-trel}
For all $x_1\in \M$, there exists a dense subset $N_{x_1} \subset \M$ such that for all $x \in N_{x_1}$, there is a (unique) minimizing geodesic path between $x_1$ and $x$. Moreover, this path is a normal geodesic path.
\end{theorem}

As a direct corollary, we obtain the following result, which is a key step in the proof of Theorem~\ref{thmwavehypo} (in particular for obtaining the minimal time).
\begin{corollary}
\label{cor:time-optim}
Let $\omega$ a nonempty open subset of $\M$ and $T>\sup_{x\in \M} d_\L(x,\omega)$. Then, for any $x_1\in \M$, there exists $x_0\in \omega$ and a normal geodesic path $\gamma : [0,L] \to \M$ of length $L\in(0,T)$, so that $\gamma(L)=x_1$, $\gamma(0)=x_0$, and $\gamma$ is also a minimizing geodesic path.
\end{corollary}
Note that this path being minimizing, it is in particular non self-intersecting.
\begin{proof}
According to the definition of $T$, the set $\O_{x_1} := \omega \cap \{x \in \M, d_\L(x,x_1)<T\}$ is open (the continuity of $d_\L$ is a consequence of the Chow-Rashevski Theorem~\ref{t:chow}) and nonempty. Hence, it intersects the dense set $N_{x_1}$ given by Theorem~\ref{th:riff-trel}. Taking any $x_0 \in \O_{x_1} \cap N_{x_1}$, there is a normal geodesic curve $\gamma : [0,L] \to \M$ joining $x_0$ and $x_1$ of length $L\in(0,T)$ which is also a minimizing geodesic. 
\end{proof}

Now, one of the main purposes of the present Section~\ref{s:hypo-wave} is to give a proof of the following proposition (which essentially amounts to~\eqref{e:intro-estim-partial-wave}). Indeed, Theorem~\ref{thmwavehypo-s} in the simple case $s=k$ will follow. The proof of the full range of $s$ in Theorem~\ref{thmwavehypo-s} will require more work.
\begin{proposition}
\label{prop:local}
Let $T>0$, $x_0 , x_1\in \M$, and assume that there is a normal geodesic path of length $L\in(0,T)$ between $x_0$ and $x_1$. Then, for any $\e>0$, there exists $\widetilde{\e}>0$, there is $C,\kappa,\mu_0>0$ such that for all $u \in H^{1}(]-T,T[\times \M)$ solution of $(\d_t^2+\L)u=0$ on $]-T,T[\times \M$, and for all $\mu \geq \mu_0$, we have
$$
\nor{ u}{L^2(]-\widetilde{\e},\widetilde{\e}[\times B(x_1,\widetilde{\e}))}\leq C e^{\kappa \mu}\nor{u}{L^2(]-T,T[\times B(x_0,\e))}+\frac{C}{\mu}\nor{u}{H^{1}(]-T,T[\times \M)} .
$$
\end{proposition}

Note here that the $H^{1}(]-T,T[\times \M)$ norm is the usual one issued from the structure of compact manifold. Also, since $\e$ and $\widetilde{\e}$ are arbitrary small, the balls could be defined according to any metric on $\M$ defining an equivalent topology (balls could equivalently be replaced by neighborhoods). Yet, we will use the distance induced by the sub-Riemannian geometry since it is the important one in other parts of the proof (for defining the distance $d_\L$ for instance) and to avoid any confusion.

Using Proposition~\ref{prop:local}, together with Corollary~\ref{cor:time-optim} and a compactness argument directly yields the following global estimate, which is the main result of this step.
\begin{corollary}
\label{cor:H1th}
Let $\omega$ a nonempty open set $\M$ and $T>\sup_{x\in \M} d_\L(x,\omega)$. Then, there exist $\widetilde{\e}, C,\kappa,\mu_0>0$ such that for all $u \in H^{1}(]-T,T[\times \M)$ solution of $(\d_t^2+\L)u=0$ on $]-T,T[\times \M$, and for all $\mu \geq \mu_0$, we have, 
$$
\nor{ u}{L^2(]-\widetilde{\e},\widetilde{\e}[\times \M)}\leq C e^{\kappa \mu}\nor{u}{L^2(]-T,T[\times \omega)}+\frac{C}{\mu}\nor{u}{H^{1}(]-T,T[\times \M)} .
$$
\end{corollary}
The last step towards the proof of Theorem~\ref{thmwavehypo-s} (in the case $s=k$; the general case requires a slightly more precise version of this result), performed in Section~\ref{s:energy-estimates},  will be to transfer the time-space information carried by this inequality into some Sobolev norm $\H^s_\L$ related to the operator $\L$. This will be the object of the next paragraph.

\bigskip
Sections~\ref{subsectstepgeometry} and~\ref{sectstepsmallness} are devoted to the proof of Proposition~\ref{prop:local}. 
For this, the main point is to apply the quantitative unique continuation result adapted to changes of variable, namely Theorem~\ref{thmsemiglobaldepchgt}. The main drawback of this result is that it only works in a subset $\Omega$ of $\R^n$, i.e. it is not invariant by diffeomorphism. More precisely, all estimates obtained from the results of Section~\ref{subsectsemiglobal} can only be patched together in straight coordinates.

As a consequence, in the present context, we need to introduce {\em global} coordinates near the trajectory between the points $x_0$ and $x_1$. Then, locally, we shall define hypersurfaces to match the geometric setting of Section~\ref{s:typical-result}. This will be done in another set of coordinates (not necessarily analytic), which is allowed by the precise formulation of Theorem~\ref{thmsemiglobaldepchgt}.

Finally, the energy estimates needed to conclude the proof of Theorem~\ref{thmwavehypo-s} are performed in Section~\ref{s:energy-estimates}.

\subsection{Step 1: Geometric setting and non-characteristic hypersurfaces}
\label{subsectstepgeometry}

We now consider $x_1, x_0 \in \M$ (recall that $\dim(\M) = d$), together with a normal geodesic $\gamma$ as in Corollary~\ref{cor:time-optim}. The curve $\gamma$ is non self-intersecting, so that there exists an open neighborhood $\mathcal{N}_\gamma$ of $\gamma$ in $\M$, an open set $\Omega_\gamma \subset \R^d$ and an {\em analytic} diffeomorphism 
\bnan
\label{phigamma}
\phi_\gamma : \mathcal{N}_\gamma \to \Omega_\gamma \subset \R^d
\enan
such that $\phi_\gamma(\gamma([0,L])) \subset \Omega_\gamma$. 
From now on, \textbf{we shall only work in the chart $(\Omega_\gamma , \phi_\gamma)$}.
For the sake of readability, we shall keep the same notation for all objects pulled back form $\mathcal{N}_\gamma \subset \M$ to $\Omega_\gamma \subset \R^d$.
For instance, we shall still denote $X_j$ instead of $(\phi_\gamma^{-1})^*X_j$, $\L$ instead of $(\phi_\gamma^{-1})^*\L \phi_\gamma^*$, $\gamma \subset \Omega_\gamma$ instead of $(\phi_\gamma^{-1})^*\gamma$ etc...
Recall that all above definitions (in particular Definitions~\ref{def:Hamilton}, \ref{def:normal} and~\ref{def:minimizing}) are intrinsic, so that, in particular, the whole sub-Riemannian structure may be transported to $\Omega_\gamma$, and the curve $\gamma$ is still a normal geodesic in $\Omega_\gamma$.

Now, we define other local coordinates in which to construct the (local) noncharacteristic surfaces in order to apply Theorem~\ref{thmsemiglobaldepchgt}. We first need the following lemma.
\begin{lemma}[local coordinates]
\label{lmlocalcoord}
Given $\gamma : [0,L] \to \Omega_\gamma$ a normal geodesic path, for any point $x_0  =\gamma(s_0)$ on this curve, there exists an open neighborhood $V_{x_0}$ of $x_0$, and appropriate coordinates, denoted $x= (\check{x},x_d) \in \R^{d-1} \times \R$ (with associated cotangent variables  $\xi = (\check{\xi},\xi_d) \in \R^{d-1} \times \R$) in which
\begin{itemize}
\item the symbol $\ell$ can be written
\bnan
\label{Pnormal}
\ell(x, \xi)=\xi_d^2+r(x,\check{\xi}) ,
\enan
where $r(x,\check{\xi})$ is a non negative quadratic form in $\check{\xi}$;
\item the point $x_0$ is sent to $(0,\cdots,0,s_0)$ and $\gamma(s)=(0,\cdots,0,s )$ for $s$ close enough to $s_0$.
\end{itemize}
\end{lemma}
\bnp
This is e.g. a consequence of~\cite[Appendix~C.5]{Hoermander:V3}. 
More precisely, we denote by $s \mapsto (\gamma(s), \xi(s)) \in T^*\Omega \setminus 0$ the Hamiltonian curve associated to the normal geodesic $\gamma$. We let $\varphi$ be a real-valued function defined locally in a neighborhood of $x_0$, such that $\varphi(x_0) = 0$ and $d\varphi(x_0) = \xi(s_0)$. Then, $\varphi$ is a non-characteristic function for $\ell$ at $x_0$ since $\ell(x_0, d\varphi(x_0)) = \ell (\gamma(s_0) , \xi(s_0)) = \frac14 \neq 0$, according to the definition of a normal geodesic path. According to~\cite[Corollary~C.5.3]{Hoermander:V3}, there are local coordinates $(\check{x},x_d) \in \R^{d-1} \times \R$, defined in a neighborhood of $0$ in which
\begin{itemize}
\item $x_0$ is sent to $0$,
\item the surface $\{\varphi = 0\}$ is given by $\{x_d=0\}$,
\item the first item of the lemma holds. 
\end{itemize}
We now just have to check that the second item of the lemma holds in these coordinates.
First remark that $d\varphi(x_0) = \xi(s_0)$ is sent to $(0, \xi_d)$ for some $\xi_d \in \R^*$, so that $\ell (\gamma(s_0) , \xi(s_0)) = \frac14$ implies, with the form of $\ell$ in \eqref{Pnormal}, that $\xi_d^2=\frac14$. Up to changing $x_d \mapsto -x_d$ (without changing any of the three properties described in the above items), we may further assume that $\xi_d >0$. Hence $d\varphi(x_0) = \xi(s_0)$ is sent to $(\check{\xi},\xi_d) = (0, \frac12)$.
Second, the form of $\ell$ in~\eqref{Pnormal} yields that the Hamiltonian curves of $\ell$ satisfy in these coordinates: 
\bnan
\label{e:hamilt}
\dot{\check{x}}  = \d_{\check{\xi}}r(x,\check{\xi}) , \quad \dot{x}_d = 2\xi_d , \quad \dot{\check{\xi}} =- \d_{\check{x}}r(x,\check{\xi}) , \quad \dot{\xi}_d = - \d_{x_d}r(x,\check{\xi})  .
\enan
The Hamiltonian curve associated to the normal geodesic $\gamma$ in these coordinates is the unique curve of~\eqref{e:hamilt} passing through $x=0$ and $(\check{\xi},\xi_d) = (0, \frac12)$. But the function $(\check{x}, x_d, \check{x} , \xi_d)(s) = (0, s- s_0, 0, \frac12)$ solves \eqref{e:hamilt} since $\d_{\check{x}}r$, $\d_{x_d}r$ are quadratic in $\check{\xi}$ and $\d_{\check{\xi}}r$ is linear in $\check{\xi}$ (and, in particular, all vanish at $\check{\xi}=0$). It also starts at time $s_0$ at $(0,0, 0, \frac12)$, so that $(0, s-s_0, 0, \frac12)$ is the sought Hamiltonian curve. As a consequence, the normal geodesic $\gamma$ is given by $(\check{x}, x_d)(s) = (0, s-s_0)$ in these coordinates.
This concludes the proof after the linear change of variable $(\check{x},x_d)  \mapsto (\check{x}, x_d + s_0)$.
\enp
\begin{lemma}
[Construction of non characteristic hypersurfaces in normal coordinates]
\label{lmconstruction}
Assume that for some $r_0, l_0>0$, the symbol $\ell$ is given by~\eqref{Pnormal} in coordinates $(\check{x}, x_d) \in \widetilde{\Omega} :=B(0,r_0) \times ]-l_0,2l_0[ \subset \R^{d-1}_{\check{x}} \times \R_{x_d}$.
Then, for any $t_0>l_0$ and $0<r_1<r_0$, there exists $D$, $G$, $K$, $S_{\e}$ satisfying items~\ref{geom-set-1}-\ref{geom-set-2}-\ref{geom-set-3} of the Geometric Setting of Section~\ref{s:typical-result} in $\R^n = \R^{d+1}$ in the coordinates
\bnan
\label{e:coord-x-x}
(\z' , \z_n) =  (t, \check{x} , x_d), \quad  \text{ with } \quad  \z'  =  (t, \check{x}) \quad \text{ and } \quad \z_n = x_d ,
\enan
together with 
\begin{enumerate}
  \setcounter{enumi}{3}
\item \label{item1geom} $D \subset [-t_0,t_0] \times \ovl{B}(0,r_1)$, that is $S_0\subset [-t_0,t_0]_{t}\times \ovl{B}(0,r_1)_{\check{x}} \times \{0\}_{x_d} \subset \R_t \times \subset \R^{d-1}_{\check{x}} \times  \{0\}_{x_d}$; 
\item \label{item2geom}$\{0\}_{t}\times \{0\}_{\check{x}}\times [0,l_0]_{x_d}\subset K$;
\item \label{item3geom} for any $\e\in [0,1+\eta)$, the surface $S_{\e} $ is {\em non-characteristic} for $P=\partial_t^2+\L$ at each point of $S_\eps$.
\end{enumerate}
\end{lemma}
\bnp
The principal symbol of the operator $P=\partial_t^2+\L$ in the coordinates of Lemma~\ref{lmlocalcoord} is given by 
\bnan
\label{fromdiagm}
p(\check{x},x_d, \tau, \check{\xi}, \xi_d) =- \tau^2 + \ell(\check{x},x_d, \tau, \check{\xi}, \xi_d)= - \tau^2  + \xi_d^2 + r(x,\check{\xi}) ,  \quad\xi = (\check{\xi}, \xi_d), 
\enan
To match the geometric setting of Section~\ref{s:typical-result}, we define the coordinates $(\z',\z_n)$ according to~\eqref{e:coord-x-x}, as well as
 \bna
 D=\left\{(t,\check{x})\left| \Big(\frac{\check{x}}{r_1}\Big)^2+\Big(\frac{t}{t_0}\Big)^2 < 1\right.\right\}, \qquad 
G(t,\check{x},\e)=  \e l_0\psi \left(\sqrt{\Big(\frac{\check{x}}{r_1}\Big)^2+\Big(\frac{t}{t_0}\Big)^2} \right) , \\
 \quad 
 \phi_\eps(t,\check{x},x_d): =  G(t,\check{x} , \eps) - x_d , \qquad S_\eps = \{ \phi_\eps = 0\},  \quad  \e \in [0,1 + \eta) ,
\ena with $r_1, \eta >0$ small to be fixed, 
where $\psi$ is such that, for some $\eta_0,\eta_1>0$,
 \bna
& \psi : [-1-\eta_0,1+\eta_0]\to [-\eta_1,1] , \quad \text{smooth and  even,}\quad  \psi(\pm1) =0, \quad \psi(0)=1,&\\
&\psi(s)\geq 0,  \textnormal{ if and only if }s\in [-1,1], \quad\text{ and } |\psi'| \leq \alpha \text{ on }[-1-\eta_0,1+\eta_0], &
 \ena
with $1<\alpha <\frac{t_0}{l_0}$. This is possible since $\frac{t_0}{l_0}>1$.

Note first that Item~\ref{item1geom} is satisfied according to the definition of $D$.
Note also that the point $(t,\check{x},x_d)=(0,0,l_0)$ belongs to $S_1= \left\{\phi_{1}=0\right\}$. Hence, Item \ref{item2geom} is satisfied since $\{0\}_{t}\times \{0\}_{\check{x}}\in D$ and $G(0_{t,\check{x}},1)=l_0$, so that we have $0\leq x_d\leq G(0_{t,\check{x}},1)=l_0$ if $x_d\in [0,l_0]$.

Let us now check Item~\ref{item3geom}. We have
  $$
d \phi_\eps(t,\check{x},x_d) 
=   \e l_0  \left(\Big(\frac{\check{x}}{r_1}\Big)^2+\Big(\frac{t}{t_0}\Big)^2\right)^{-1/2} \psi' \left(\sqrt{\Big(\frac{\check{x}}{r_1}\Big)^2+\Big(\frac{t}{t_0}\Big)^2} \right) \left( \frac{tdt}{t_0^2}  + \frac{\check{x} d\check{x}}{r_1^2}  \right) - d x_d .
 $$
Given the form of the principal symbol of the operator $P$ in these coordinates (see \eqref{fromdiagm}), we obtain
\bna
- p( \check{x},x_d, d \phi_\eps(t,\check{x},x_d) ) =  \e^2 l_0^2\frac{t^2}{t_0^4}    \left(\Big(\frac{\check{x}}{r_1}\Big)^2+\Big(\frac{t}{t_0}\Big)^2 \right)^{-1}|\psi'|^2  
- l_0^2\frac{\e^2}{r_1^4}  r(x,\check{x}) \left(\Big(\frac{\check{x}}{r_1}\Big)^2+\Big(\frac{t}{t_0}\Big)^2 \right)^{-1} |\psi'|^2  
- 1
\ena
where $|\psi'|^2 $ is taken at the point $\left(\sqrt{\Big(\frac{\check{x}}{r_1}\Big)^2+\Big(\frac{t}{t_0}\Big)^2} \right)$. Since $r$ is non negative, we get
\bna
-p( \check{x},x_d, d \phi_\eps(t,\check{x},x_d) ) \leq   \e^2 l_0^2\frac{t^2}{t_0^4}    \left(\Big(\frac{\check{x}}{r_1}\Big)^2+\Big(\frac{t}{t_0}\Big)^2 \right)^{-1}|\psi'|^2  
- 1
\ena
Since $|\psi'|\leq \alpha $ and   $\e\in [0,1+ \eta]$, we obtain for any $(t,\check{x},x_d)\in \overline{D}\times [0,l_0]$,
\bna
- p( \check{x},x_d, d \phi_\eps(t,\check{x},x_d) )& \leq&   \frac{\e^2}{t_0^2}l_0^2   \Big(\frac{t}{t_0}\Big)^2\left(\frac{t}{t_0} \right)^{-2}  \alpha^2  
-  1  \\
& \leq &(1+\eta)^2 \frac{l_0^2}{t_0^2} \alpha^2-1 <0,
 \ena
the last constant being negative for $\eta$ small enough because $\alpha<\frac{t_0}{l_0}$.
Therefore, the surface $S_\e = \{\phi_\eps =0\}$ is noncharacteristic for any $\eps \in [0,1+\eta ]$, which concludes the proof of Item~\ref{item3geom}, and hence of the lemma.
\enp

\subsection{Step 2: Propagation of smallness}
\label{sectstepsmallness}
We recall that $\z = (t,x)$ and $n=d+1$, and introduce the notation 
$$
\vois(K,r): = \bigcup_{\z \in K} B(\z, r) \quad \text{ for } K \subset \R^n .
$$

\begin{figure}[h!]
  \begin{center}
    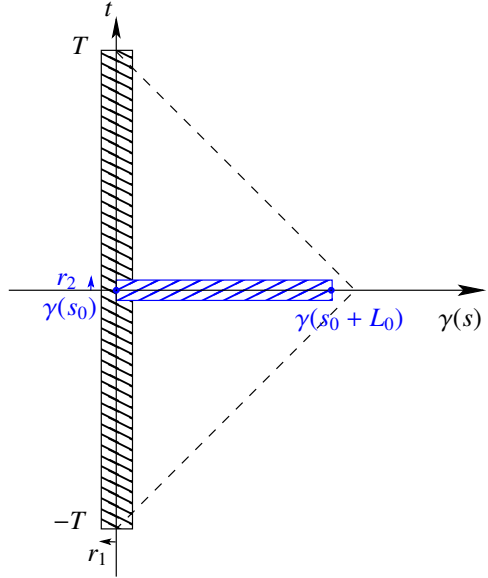 
    \caption{Schematic representation of the sets involved in~ Lemma~\ref{lemhypolocal}.}
    Recall that $\length(\gamma([s_0,s_0+L_0]))= L_0 < T$; the dashed line represents the boundary of the wave cone.
    \label{f:stick1}
 \end{center}
\end{figure}

\begin{lemma}[local version near a piece of a normal geodesic]
\label{lemhypolocal}
Given $\gamma : [0,L] \to \Omega_\gamma$ a normal geodesic path and fix $\tilde{s}\in [0,L]$. Then, there exists $\tilde{L}_s, \tilde{r}_s$ small such that for any  $s_0\in [0,L]$, $L_0>0$ so that $[s_0,s_0+L_0]\subset ]\tilde{s}-\tilde{L}_s,\tilde{s}+\tilde{L}_s[$, for all $T>L_0$ and $0<r_1<\tilde{r}_s$, there exists $r_2>0$ so that 
\bnan
\label{deptriloc}
]-r_2,r_2[ \times \vois(\gamma([s_0,s_0+L_0]),r_2)\lhd ]-T,T[ \times \vois(\gamma(s_0),r_1) ,
\enan
where $\lhd$ is related to the operator $P=\d_t^2+\L$ in the set $]-T-\e,T+\e[\times \Omega_\gamma$.
\end{lemma}
See Figure~\ref{f:stick1} for a picture of the sets involved in~\eqref{deptriloc}.

The proof is almost the same as in Theorem 6.3 of \cite{LL:15}. The only difference is that the coordinates where we have a nice diagonal form for the operator $P$ are not global and are not those where we want to apply the local result. Note that this would not have been a problem if we had proved that the relation $\lhd$ is invariant by change of coordinates. Now, we perform the following steps:
\begin{itemize}
\item use Lemma \ref{lmlocalcoord} to obtain nice coordinates in a neighborhood of $\gamma(s)$; 
\item construct the non characteristic hypersurfaces in these coordinates according to Lemma~\ref{lmconstruction};
\item apply Theorem \ref{thmsemiglobaldepchgt} in these coordinates, w.r.t. those surfaces, keeping in mind that the fact to be non characteristic is invariant by changes of coordinates.
\end{itemize}
Note also that the presence of $\tilde{s}$ and $s_0$ in the statement may seem strange and it would look simpler to consider only $s_0$ and intervals $[s_0,s_0+L_0]$. Yet, this will be useful later in compactness and covering arguments where we have to substract and cut some intervals. 
\bnp
By Lemma \ref{lmlocalcoord}, we can find a diffeomorphism from a neighborhood of $\gamma(\widetilde{s})$ in $\Omega_\gamma$ to $B_{\R^{d-1}}(0,r_0) \times ]\widetilde{s}-4\tilde{L}_s,\tilde{s}+4\tilde{L}_s[ \ni (\check{x} , x_d)$ in which the symbol $\ell$ is as in~\eqref{Pnormal} and $\gamma(s)$ is given by $(\check{x} , x_d) = (0,\widetilde{s})$. By translation and using $[s_0,s_0+L_0]\subset ]\tilde{s}-\tilde{L}_s,\tilde{s}+\tilde{L}_s[$, we have a diffeomorphism $\Phi$ from a neighborhood $V$ of $\gamma(s_0)$ in $\Omega_\gamma$ onto $B_{\R^{d-1}}(0,r_0)\times ]-L_0,2L_0[  \ni (\check{x} , x_d)$ where $\Phi(\gamma(s_0))=0$. 

For later purposes, fix $r_3$ so that 
\bnan
\label{inclusomegacompact}
\Phi^{-1}(B_{\R^{d}}(0,r_3))\Subset B_{\R^{d}}(\gamma(s_0),r_1).
\enan 

We will keep denoting $\Phi$ the same diffeomorphism acting on $]-T-\e,T+\e[\times V$ leaving the $t$ variable unchanged, and set $\widetilde{\Omega}: = ]-T-\e,T+\e[ \times B_{\R^{d-1}}(0,r_0)\times ]-L_0,2L_0[$. 

Now, construct $D$, $G$, $K$, $S_{\e}$ according to Lemma \ref{lmconstruction}. In particular, all surfaces $S_\eps$ are non characteristic for $P$ (or, more precisely, for $(\Phi^*)^{-1}P \Phi^*$) in $\widetilde{\Omega}$). By change of coordinates, the surface $\Phi^{-1}(S_{\e})$ is non characteristic for $P$ at each point of $\Phi^{-1}(S_{\e})$.

For $0<L_0<t_0<T$, take $\omega=\Phi^{-1}(]-t_0,t_0[ \times B(0,r_3))$ so that Item \ref{item1geom} of Lemma \ref{lmconstruction} implies that $\omega$ is a neighborhood of $\Phi^{-1}(S_{0})$. 

The assumptions of Theorem \ref{thmsemiglobaldepchgt} are fulfilled (with $n_a=n=d+1$, i.e. in the Holmgren-John case), so there exists an open neighborhood $U\subset \Omega$ of $\Phi^{-1}(K)$ such that $U\lhd \omega$. Note here that the strict application of Theorem \ref{thmsemiglobaldepchgt} yields this result for the relation $\lhd_{]-T-\e,T+\e[\times V,P}$, but the latter then implies the same property for the relation $\lhd_{]-T-\e,T+\e[\times \Omega_{\gamma},P}$. We will keep the notation~$\lhd$.

Moreover, Item \ref{item2geom} of Lemma \ref{lmconstruction} implies, after having applied $\Phi^{-1}$, that
\bna
\Phi^{-1}\left(\{0\}_{t}\times \{0\}_{\check{x}}\times [0,l_0]_{x_d}\right)\subset \Phi^{-1}(K)\subset U .
\ena
Using the form of $P$ on $\widetilde{\Omega}$ and that $\gamma$ is a normal geodesic, we obtain that $\Phi^{-1}\left(\{0\}_{t}\times \{0\}_{\check{x}}\times [0,L_0]_{x_d}\right)=\{0\}_{t}\times\gamma([s_0,s_0+L_0])$. In particular since $U$ is open and the previous set is compact, we can find $r_2>0$ so that $]-r_2,r_2[\times \vois(\gamma([s_0-r_2,s_0+L_0]),r_2) \Subset U$.
Items \ref{propincludestrong} and \ref{proptransstrong} of Proposition \ref{propstrong} imply 
\bna
]-r_2,r_2[_{\R_t}\times \vois(\gamma([s_0-r_2,s_0+L_0]),r_2)\lhd \omega.
\ena
Finally, the definition of $\omega$, $T>t_0$ and \eqref{inclusomegacompact} imply
$\omega \Subset ]-T,T[ \times \vois(\gamma(s_0),r_1)$. This gives the final result by applying again Item \ref{propincludestrong} and \ref{proptransstrong} of Proposition \ref{propstrong}.
\enp
We can iterate the previous local result to get a more global one, which will be the main step for Proposition~\ref{prop:local}.
\begin{proposition}[global version near a normal geodesic]
\label{propdepgeodes}
Given $\gamma : [0,L] \to \Omega_\gamma$ a normal geodesic path, and let $0<L<T$. 
Then, there exists $r_0$ small, such that for any $0<r_1<r_0$, there exists $r_2>0$ such that 
\bnan
\label{deptriglob}
]-r_2,r_2[\times \vois(\gamma([0,L]),r_2)\lhd ]-T,T[\times \vois(\gamma(0),r_1).
\enan
\end{proposition}
\bnp[Proof of Proposition~\ref{propdepgeodes}]
We prove the result for another $\widetilde{T}>T$. It gives the result since it is arbitrary. 

For any $\tilde{s}\in[0,L]$, Lemma \ref{lemhypolocal} provides $\tilde{L}_s$ and $\tilde{r}_s$ and an interval $]\tilde{s}-\tilde{L}_s,\tilde{s}+\tilde{L}_s[$ with the appropriate conclusion. By compactness of $\gamma([0,L])$, we can extract a finite covering such that  $\gamma([0,L])\subset \bigcup_{j\in \llbracket 1, N\rrbracket} \gamma(]\tilde{s}_j-\tilde{L}_j,\tilde{s}_j+\tilde{L}_j[)$. Then, the issue is that $\sum_{j=1}^N2\tilde{L}_j$ may be very large with respect to $2T$. To overcome this difficulty, starting from this covering, we can always obtain (for this, suppress some intervals and cut them when they overlap too much) a finite number of intervals $[s_j,s_j+L_j]$ and associate times $T_j$ that satisfy the following properties:
\begin{itemize}
\item $[0,L]\subset \cup_{j=1}^N]s_j,s_j+L_j[$,
\item $[s_j,s_j+L_j]\subset ]\tilde{s}_j-\tilde{L}_j,\tilde{s}_j+\tilde{L}_j[$, for all $j\in \llbracket 1, N\rrbracket$,
\item $s_1=0$,
\item $s_{j+1}\in ]s_j,s_j+L_j[$, for all $j\in \llbracket 1, N-1\rrbracket$,
\item $L<\sum_{j=1}^N L_j< T$,
\item $L_j<T_j$ and $\sum_{j=1}^N T_j< T$.
\end{itemize}
So, for any $j\in \llbracket 1, N\rrbracket$, since  $[s_j,s_j+L_j]\subset ]\tilde{s}_j-\tilde{L}_j,\tilde{s}_j+\tilde{L}_j[$, Lemma \ref{lemhypolocal} can be applied to the path $\gamma([s_j,s_j+L_j])$ and gives the existence of $\tilde{r}_s$ associated to $\tilde{s}_j$, which we here denote $r_0^j$. We also denote by $r_0$ the minimum of all $r_0^j$, $j\in \llbracket 1, N\rrbracket$, so that the conclusion of Lemma \ref{lemhypolocal} remains true with any choice of $r_1^j<r_0$. We next define $r_1^j$ and $r_2^j$ recursively in the following way:
\begin{itemize}
\item $r_1^1=\min(r_0, r_1)/2$ and $r_2^1$ is given by the Lemma \ref{lemhypolocal} for the interval $[s_1,s_1+L_1]=[0,L_1]$.
\item We choose $r_1^{j+1}=\min(r_0, r_2^{j})/4$ and $r_2^{j+1}$ is given by Lemma \ref{lemhypolocal} applied to the path $[s_j,s_j+L_j]\subset ]\tilde{s}_j-\tilde{L}_j,\tilde{s}_j+\tilde{L}_j[$ and the time $T_j>L_j$.
\end{itemize}
The conclusion of Lemma \ref{lemhypolocal} is then
\bnan
\label{concluLM}
]-r_2^j,r_2^j[_{\R^t}\times \vois(\gamma([s_j,s_j+L_j]),r_2^j)\lhd ]-T_j,T_j[_{\R^t}\times \vois(\gamma(s_j),r_1^j).
\enan
Now, for any $l\in \llbracket 1, N\rrbracket$, consider a sequence of time $(t_i^l)_{i\in I_l}$ such that $\left(]t_i^l-r_2^{l},t_i^l+r_2^{l}[\right)_{i\in I_l}$ is a finite covering of $]-T+\sum_{j=1}^l T_j,T-\sum_{j=1}^l T_j[$. One can also impose $t_i^l\in ]-T+\sum_{j=1}^l T_j,T-\sum_{j=1}^l T_j[$.

We want to apply Proposition \ref{propiterationabstrait} with the following definitions for $j\in \llbracket 1, N\rrbracket$, $i\in I_j$
\begin{itemize}
\item $U_{i,j}=]t_i^j-r_2^j,t_i^j+r_2^j[\times \vois(\gamma([s_j,s_j+L_j]),r_2^j)$,
\item $\omega_{i,j}=]t_i^j-r_2^j/2,t_i^j+r_2^j/2[\times \vois(\gamma([s_j,s_j+L_j]),r_2^j/2)$,
\item $V_{i,j}=]t_i^j-T_j,t_i^j+T_j[\times \vois(\gamma(s_j),r_1^j)$,
\item $U_0=]-T,T[\times \vois(\gamma(0),2r_1^1)$.
\end{itemize}

Since Lemma \ref{lemhypolocal} is invariant by translation in time, \eqref{concluLM} and the choices of $r_1^j$, $ r_2^j$ give $U_{i,j}\lhd V_{i,j}$. We also have $\omega_{i,j} \Subset U_{i,j}$. So, the main point to check is 
\bnan
\label{inclus}
V_{m,l+1} \Subset \left[U_0 \cup\bigcup_{j\in \llbracket 1,l\rrbracket}\bigcup_{i\in I_{j}} \omega_{i,j}\right] ,\quad \text{for all }m \in I_{l+1} , \quad \text{and } l \in \llbracket 1, N-1\rrbracket .
\enan
We first check the degenerate case $l=0$, which amounts to proving that $V_{m,1} \Subset U_0$ for $m \in I_1$. Since $t_m^1\in ]-T+T_1,T- T_1[$ (by choice), we have $]t_m^1-T_1,t_m^1+T_1[\Subset ]-T,T[$. Moreover, since $s_1=0$, we have $\vois(\gamma(s_1),r_1^j)\Subset \vois(\gamma(0),2r_1^1)$, and $V_{m,1} \Subset U_0$ by definition.

\medskip
Concerning the case $l\in \llbracket 1,N-1\rrbracket$, we prove the stronger property
\bnan
\label{stronger}
V_{m,l+1} \Subset \bigcup_{i\in I_{l}} \omega_{i,l} ,\quad \text{for all }m \in I_{l+1} .
\enan
First, we have by definition
\bna
\bigcup_{i\in I_{l}} \omega_{i,l}= \left[\bigcup_{i\in I_{l}} ]t_i^l-r_2^l/2,t_i^l+r_2^l/2[\right]\times \vois(\gamma([s_l,s_l+L_l]),r_2^l/2).
\ena
Since $\left(]t_i^l-r_2^{l},t_i^l+r_2^{l}[\right)_{i\in I_l}$ is a finite covering of $]-T+\sum_{j=1}^l T_j,T-\sum_{j=1}^l T_j[$, we have 
\bnan
\label{inclusomega}
]-T+\sum_{j=1}^l T_j,T-\sum_{j=1}^l T_j[\times \vois(\gamma([s_l,s_l+L_l]),r_2^l/2) \subset \bigcup_{i\in I_{l}} \omega_{i,l} .
\enan
We also have $t_i^{l+1}\in ]-T+\sum_{j=1}^{l+1} T_j,T-\sum_{j=1}^{l+1} T_j[$, so that 
$$
]t_i^{l+1}-T_{l+1},t_i^{l+1}+T_{l+1}[\Subset ]-T-T_{l+1}+\sum_{j=1}^{l+1} T_j,T+T_{l+1}-\sum_{j=1}^{l+1} T_j[,
$$
 that is
\bnan
\label{inclusint}
]t_i^{l+1}-T_{l+1},t_i^{l+1}+T_{l+1}[\Subset]-T+\sum_{j=1}^{l} T_j,T-\sum_{j=1}^{l} T_j[.
\enan
Moreover, $\gamma(s_{l+1})\in \gamma(]s_l,s_l+L_l[)$ and $r_1^{l+1}<r_2^l/2$ by construction, so 
\bnan
\label{inclusVois}
\vois(\gamma(s_{l+1}),r_1^{l+1})\Subset \vois(\gamma([s_l,s_l+L_l]),r_2^l/2).
\enan
Combining the definition of $V_{i,l+1}=]t_i^{l+1}-T_{l+1},t_i^{l+1}+T_{l+1}[\times \vois(\gamma(s_{l+1}),r_1^{l+1})$, \eqref{inclusint}, \eqref{inclusVois} and \eqref{inclusomega}, we obtain $V_{i,l+1}\Subset \bigcup_{i\in I_{l}} \omega_{i,l}$. This finishes the proof of \eqref{stronger} and therefore \eqref{inclus}, so that all assumptions of Proposition~\ref{propiterationabstrait} are satisfied.

The conclusion of this proposition can then be written as
\bnan
\label{conclusionprop}
\left[U_0 \cup\bigcup_{j\in \llbracket 1,N\rrbracket}\bigcup_{i\in I_{j}} \omega_{i,j}\right]\lhd V_0 \textnormal{ for any }V_0 \text{ such that } U_0\Subset V_0.
\enan
Now, we pick $r_2< \min(T- \sum_{j=1}^N T_j,r_2^l/2)$. Using \eqref{inclusomega} and then the covering property $[0,L] \subset \bigcup_{j \in  \llbracket 1, N\rrbracket} ]s_j,s_j+L_j [$, we obtain
\bna
]-r_2,r_2[\times \vois(\gamma([s_j,s_j+L_j]),r_2/2) \subset \bigcup_{i\in I_{j}} \omega_{i,j} , \\
]-r_2,r_2[\times \vois(\gamma([0,L]),r_2/2) \subset \bigcup_{j\in \llbracket 1,N\rrbracket} \bigcup_{i\in I_{j}} \omega_{i,j}.
\ena
In particular, we have
\bna
]-r_2/4,r_2/4[\times \vois(\gamma([0,L]),r_2/4) \Subset \bigcup_{j\in \llbracket 1,N\rrbracket} \bigcup_{i\in I_{j}} \omega_{i,j}.
\ena
Now, since $\widetilde{T}>T$ and $r_1^1\leq r_1/2$, we have $U_0\Subset ]-\widetilde{T},\widetilde{T}[\times \vois(\gamma(0),2r_1):=V_0$. Combining this together with \eqref{conclusionprop} and the use of Proposition \ref{propstrong} (several times), we finally obtain
\bna
]-r_2/4,r_2/4[\times \vois(\gamma([0,L]),r_2/4) \lhd ]-\widetilde{T},\widetilde{T}[\times \vois(\gamma(0),2r_1).
\ena
This concludes the proof of Proposition \ref{propdepgeodes}.
\enp

From Proposition~\ref{propdepgeodes}, let us now briefly explain the proof of Proposition~\ref{prop:local}. We proceed exactly as in~\cite[Section~4.2]{LL:15}, in the proof that Theorem~4.7 implies Theorem~1.11 (which, in Section~\ref{s:Holmgren-john-LL15} of the present paper, corresponds to the proof that Theorem \ref{thmsemiglobaldep} implies Theorem \ref{thmsemiglobal}).
Note that it only consists in getting rid of the frequency cutoffs appearing in the definition of $\lhd$, i.e. considering all frequencies, at the cost of replacing the $e^{-\kappa \mu}$ exponentially small remainder by a $\frac{1}{\mu}$. 
This concludes the proof of Proposition~\ref{prop:local}.

\bigskip
For an application in the context of Section~\ref{s:non-analytic} (partially analytic operators with a boundary), we need to relax the condition that $\gamma$ is globally a normal geodesic.

\begin{remark}
\label{rem:propvariantbypiece}
The proof of Step 2 took the following structure.
\begin{itemize}
\item Lemma \ref{lemhypolocal} proves some relations of dependence \eqref{deptriloc} for some local region around some small part of a normal geodesic.
\item Proposition \ref{propdepgeodes} uses the fact that the relations of dependence that we obtained in Lemma \ref{lemhypolocal} can be iterated to be around a global normal geodesic to get some relation of dependence that has the same form but globally, namely \eqref{deptriglob}.
\end{itemize}
Therefore, with exactly the same iteration process as Proposition \ref{propdepgeodes}, except that we invoke Proposition \ref{propdepgeodes} itself instead of Lemma \ref{lemhypolocal}, we can obtain that the same result as Proposition \ref{propdepgeodes} is true if $\gamma$ is only geodesic by piece. This is the following Proposition. 
\end{remark}

\begin{proposition}[global version near a piecewise normal geodesic]
\label{propvariantbypiece}
The same result as Proposition \ref{propdepgeodes} is true if $\gamma$ is only normal geodesic by piece.
\end{proposition}

\subsection{Step 3: Energy estimates}
\label{s:energy-estimates}
\subsubsection{Simple energy estimates concluding the proof of Theorem~\ref{thmwavehypo-s} with $s=k$}
\label{s:simple-case}
As precised earlier in the introduction, Theorem~\ref{thmwavehypo-s} is easier to prove in the specific case $s=k$.
To conclude this proof from Corollary~\ref{cor:H1th}, it only remains to prove for solutions of~\eqref{hypoelliptic-wave} the two estimates:
\bnan
\nor{(u_0,u_1)}{L^2\times \H^{-1}_\L} \leq C_{\widetilde{\e}} \nor{ u}{L^2(]-\widetilde{\e},\widetilde{\e}[\times \M)} , \label{lowerL2} \\
\nor{u}{H^{1}(]-T,T[\times \M)} \leq C_T \nor{(u_0,u_1)}{\H^{k}_\L\times \H^{k-1}_\L} \label{upperH1}.
\enan
On the one hand, \eqref{lowerL2} is a "straightforward observability estimate", the proof of which is exactly the same as inequality~(6.9) in the proof of Theorem~6.1 in~\cite{LL:15} for the classical wave equation.

On the other hand, the proof of~\eqref{upperH1} relies on the (optimal) subelliptic estimates stated in Theorem~\ref{thmhypoestim} and Corollary~\ref{corhypoestim}.
First define the energies
\bnan
\label{def-energy}
\E_s(u)= \frac12 \| \d_t u\|^2_{\H^{s-1}_{\L}} + \frac12 \|u\|^2_{\H^{s}_{\L}} =  \frac12 \nor{(u, \d_t u)}{\H^{s}_\L\times \H^{s-1}_\L}^2 , \quad s \in \R .
\enan
Then, rewriting~\eqref{hypoelliptic-wave} as $(\d_t^2 + \L + 1)u = u$ and taking the inner product with $(\L+1)^{s-1}\d_t u$ yields
\bna
\frac{d}{dt}\E_s(u)(t) = \left((\L+1)^{\frac{s-1}{2}}u, (\L+1)^{\frac{s-1}{2}}\d_t u\right)_{L^2}, \quad \text{and hence } - \E_s(u)(t) \leq \frac{d}{dt}\E_s(u)(t)\leq \E_s(u)(t) ,
\ena
so that we have 
\bnan
\label{e:gronwall}
\E_s(u)(t) \leq C_T\E_s(u)(0) \quad \text{for } t \in [-T,T]. 
\enan
Also, according to Corollary~\ref{cor:HsHsL}, we have $\nor{u}{H^{1}(\M)}\leq C \nor{u}{\H_\L^{k}}$. 
This estimate yields 
 \bna
\nor{u}{H^{1}(]-T,T[\times M)}^2 = \int_{-T}^T \left(\nor{\partial_t u}{L^2(\M)}^2+\nor{u}{H^1(\M)}^2 \right) dt 
\leq C \int_{-T}^T \left( \E_1(u) +  \nor{u}{\H_\L^{k}}^2  \right)dt .
\ena
Since $\E_s(u) \leq \E_\sigma(u)$ for $s \leq \sigma$, we obtain
 \bna
\nor{u}{H^{1}(]-T,T[\times M)} ^2 
\leq C \int_{-T}^T  \E_k(u) (t)dt \leq C_T  \E_k(u) (0) =2C_T  \nor{(u_0,u_1)}{\H^{k}_\L\times \H^{k-1}_\L}^2 .
\ena
This proves~\eqref{upperH1}, which, combined with the estimate of Corollary~\ref{cor:H1th} and~\eqref{lowerL2} implies
$$
\nor{(u_0,u_1)}{L^2\times \H^{-1}_\L} \leq C_s e^{\kappa \mu}\nor{u}{L^2(]-T,T[\times \omega)} +\frac{C}{\mu}\nor{(u_0,u_1)}{\H^{k}_\L\times \H^{k-1}_\L}.
$$
This finally proves the estimate~\eqref{e:wavehypo-s} with $s=k$, up to changing $\mu$ by $\mu^k/C$, and $\mu_0$ and $\kappa$ accordingly. 
Estimate~\eqref{th-estimate-Lambda} is then a direct consequence of~\cite[Lemma~A.3]{LL:15} and the inequality $\nor{u}{L^2(]-T,T[\times \omega)}\leq C \nor{(u_0,u_1)}{\H^{k}_\L\times \H^{k-1}_\L}$. 

\begin{remark}
\label{rkNRJpartial-easy}
Note that the previous energy estimates do not use any analyticity property of the solution and are equally true in the partially analytic case.
\end{remark}
\begin{remark}
Until this point, all proofs work as well if $\L$ is replaced by $\L + V$ where $V$ is a time-dependent nonnegative complex-valued analytic potential. 
Beware that in Section~\ref{s:energy-general-case} below, the spectral theory used restricts the discussion to time-independent real-valued analytic potentials. 
\end{remark}

\subsubsection{Interlude: eigenfunction tunneling, a proof of Theorem~\ref{t:spec-ineq}}
\label{s:interlude}
Denoting $v(t,x)= \cos(\sqrt{\lambda_j }t)\varphi_j$, we remark that $v$ is solution to
\bneq
\partial_t^2 v+\L v&=&0\\
(v,\partial_t v)_{t=0}&=&(\varphi_j,0).
\eneq
Therefore, thanks to Theorem~\ref{thmwavehypo-s} with $s=k$, we have the estimate 
\bna
\nor{\varphi_j}{\H^1_\L}\leq Ce^{c\Lambda} \nor{v}{L^{2}(]-T,T[\times\omega)}=C_Te^{c\Lambda}\nor{\varphi_j}{L^2(\omega)} ,
\ena
with $\Lambda=\frac{\nor{\varphi_j}{\H^{k}_\L}}{\nor{\varphi_j}{L^2}} =  (\lambda_j+1)^{k/2}$. This proves Theorem~\ref{t:spec-ineq} from Theorem~\ref{thmwavehypo-s} with $s=k$ (or any given $s>0$).

\subsubsection{End of the proof of Theorem~\ref{thmwavehypo-s} in the general case}
\label{s:energy-general-case}
We now prove the appropriate energy estimates to conclude the proof of Theorem~\ref{thmwavehypo-s} for any $s>0$. Recall again that $\z = (t,x)$ and that all variables are analytic, that is $z=z_a$. Hence, we have $\zeta = \zeta_a=(\tau, \xi)$ together with $D=D_a=D_{t,x}$. 
\begin{remark}
The main idea of this section is that since we are only dealing with low frequency estimates, the remainder can be taken in any "arbitrary weak" norm; this gives then a lot of flexibility for the the norms that we can finally take.

Yet, in our estimates, the remainder term is actually in $H^{1}$. It would certainly be possible to obtain a weaker norm in the general Theorem of Quantitative Unique Continuation of~\cite{LL:15}, but it would require to revisit all proofs of~\cite{LL:15}, and even the Carleman estimate itself. Instead, we could try to only apply our general estimate to a function with low frequency, for example $u=m^{\mu}(D_{t,x})v$ where $v$ is solution of $Pv=0$. But in this case, to obtain good estimates of $Pu$, we need some nice exponential estimate for the commutator $[P,m^{\mu}(D_{t,x})]$. This would certainly be possible but quite lengthy and would require some analytic regularity properties for $P$. Instead, for the specific case of the wave type equations,  regularity in time morally implies regularity in space. That is why it is more convenient to consider a regularizer $m^{\mu}(D_{t})$ which commutes exactly with the wave operator. Remark that the method may extend to other {\em evolution} equations.

\end{remark}
\begin{remark}
In this section, we are in the case where all variables are analytic, that is $z_a =z = (t,x)$. Therefore, the proof below does not a priori apply to the partially analytic case described in Section~\ref{s:without-analyticity}. This explains why the statement of Theorem \ref{t:partially-anal} is slightly less general. The same results might be true in the partially analytic case but would certainly require more work and additional arguments.
\end{remark}

\bigskip
In the relation $\lhd$, the remainder terms are always measured in the norm $H^1$. In this section, we describe how, in the case of solutions of the wave equation, the remainder term can be chosen in any weak Sobolev norm. 

\bigskip
 The starting point of the proof is Proposition~\ref{propdepgeodes}, which implies the following statement (without the use of the notation $\lhd$).
Let $\gamma : [0,L] \to \Omega_\gamma$ a normal geodesic path such that $\gamma(0) \in \omega$, and let $0<L<T$. 
Then, there exists $\eps>0$ such that for any $\vartheta \in C^\infty_0(\R^{d+1})$ equal to one in a neighborhood of $]-L,L[\times \{\gamma(0)\}$ and $\widetilde{\vartheta} \in C^\infty_0(\R^{d+1})$, supported in $]-2\e,2\e[\times B(\gamma(L),2\e)$, we have: for all $\kappa, \alpha>0$, there exist $C, \kappa', \beta,\mu_0 >0$ such that for all $(\mu,v) \in [\mu_0, + \infty)\times C^{\infty}_0(\R^{d+1})$, we have 
 \bnan
\label{inegtriang}
\nor{M^{\beta\mu}_\mu \widetilde{\vartheta}_{\mu}  v}{1}\leq C e^{\kappa \mu}\left(\nor{M^{\alpha\mu}_\mu \vartheta_{\mu} v}{1} + \nor{Pv}{L^2(\widetilde{\Omega}_{\gamma})}\right)+Ce^{-\kappa' \mu}\nor{v}{1} .
\enan
where $ \widetilde{\Omega}_{\gamma}=]-T,T[\times \Omega_{\gamma}$.
In particular, we may (and we shall) take $\vartheta$ supported in $]-T,T[\times \omega$ and $\widetilde{\vartheta}=1$ in a neighborhood of $]-\e,\e[\times B(x_0,\e)$. Notice also that we can impose $\beta<\alpha<1$; indeed, the estimate for a smaller $\beta$ is actually worse (than that for a larger $\beta$) up to an error term of the form $Ce^{-c \mu}\nor{v}{1}$ (see e.g. Lemma~\ref{Lemma23}, Item 3 below).

Up to choosing the coordinate charts $\Omega_{\gamma}$ smaller, we can still select some other open set $\Omega\subset \R^d$ with $\Omega_{\gamma}\Subset \Omega$ so that there is an analytic diffeomorphism from an open neighborhood of $\gamma$ to $\Omega$ and such that this diffeomorphism coincides with the $\phi_\gamma $ defined in \eqref{phigamma} on $\Omega_{\gamma}$.

\bigskip

From this starting point, the plan of the proof of Theorem~\ref{thmwavehypo-s} is to apply \eqref{inegtriang} to the function 
\bnan
\label{e:def-u-from-v}
v=\chi^0(t)\chi_{\mu}^1(\frac{D_t}{\mu})\chi_{\mu}^2(t)\chi^3(x)u, 
\enan
where $u$ is a solution of 
\bneqn
\label{hypoelliptic-wave-v}
Pu = (\partial_t^2 + \L) u&=&0\\
(u,\partial_t u)_{t=0}=(u_0,u_1) & = &U_0 ,
\eneqn
pulled back to the local coordinate chart $\Omega$ (we however keep the notation $u$ to lighten the notation). 

The time, frequency, and spatial cut off functions $\chi^i$ are chosen as follows:
\begin{itemize}
\item $\chi^2 \in C^{\infty}_0(]-T-\eps, T+\eps[)$ so that $\chi^2(t)=1$ on $]-T,T[$; we also write $\chi^2_\mu = (\chi^2)_\mu$,
\item $\chi^0 \in C^{\infty}_0(]-T-\eps, T+\eps[)$ so that $\chi^0(t)=1$ in an neighborhood of $\supp(\chi^2)$,
\item $\chi^3 \in C^{\infty}_0(\Omega)$ so that $\chi^3(x)=1$ for $x \in \Omega_{\gamma}$, 
\item $\chi^1 \in C^{\infty}_0(\R)$ is supported in $]-2,2[$ and such that $\chi^1(\tau)=1$ for $\tau \in ]-1,1[$; we also write $\chi^1_\mu = (\chi^1)_\mu$.
\end{itemize}
\begin{remark}
Note that since $\chi^3\in C^{\infty}_0(\Omega)$ and $\chi^0 \in C^{\infty}_0(\R)$, $v$ is a well defined function in $C^{\infty}_0(\R^{d+1})$ if $U_0$ and hence $u$ is smooth. So, \eqref{inegtriang} is applicable. Moreover, $v$ is not a local function in terms of $u$ in the time variable; indeed, it depends on all values of $u(t)$ for $t\in \R$. However, since $u$ is a solution to Equation~\eqref{hypoelliptic-wave-v}, $u$ may rather be seen as a function of the data $U_0$ only.
\end{remark}

Our task is now to estimate each term in~\eqref{inegtriang} in terms of the observation $\nor{u}{L^2(]-T,T[\times \omega)}$ and appropriate (weak) norms on $U_0$ (this is done respectively in Lemmata~\ref{lmtermobs}, \ref{lmtermP} and~\ref{lmtermreste} below). Finally, to conclude the proof, it will remain to estimate also the high frequency part of the solution (this is done in Lemma~\ref{lmhighweak}). Several technical lemmata and estimates used in this proof are postponed to the next section for the sake of readability.

All along this section, we shall use the following natural product norms for $s\in \R$,
$$
\nor{U_0}{s,\times} = \nor{U_0}{\H_\L^{s}\times\H_\L^{s-1}} = \sqrt{2 \E_{s}(U_0)} .
$$
\begin{lemma}
For all $s\geq 0$, there is $C, N, c, \mu_0 >0$ such that for all $\mu \geq \mu_0$ and for all $u$ solution to~\eqref{hypoelliptic-wave-v} and $v$ defined accordingly in~\eqref{e:def-u-from-v}, we have
\label{lmtermobs}
\bna
\nor{M^{\alpha\mu}_\mu \vartheta_{\mu} v}{1}\leq C\mu^{N}\nor{u}{L^2(]-T,T[\times \omega)}+Ce^{-c\mu}\nor{U_0}{-s,\times} .
\ena
\end{lemma}
\bnp[Proof of Lemma~\ref{lmtermobs}]
First note that, according to~\eqref{Estimate29} below, we have (for some $N\in \N$, changing from line to line),
$$
m_\mu(\eta) \leq C  \mu^N \langle \eta \rangle^N e^{-\frac{\mu}{4} \dist(\eta , \supp(m))^2} .
$$
Hence, recalling the definition of $m$ in Section~\ref{s:def-mult}, we obtain, for all $s\geq 0$, for some $C, N, c >0$,
\bna
|\eta|^{s+1} m_\mu(\eta) & = &|\eta|^{s+1} m_\mu(\eta)\mathds{1}_{|\eta|\leq 3/2}  +|\eta|^{s+1} m_\mu(\eta)\mathds{1}_{|\eta|\geq 3/2} \\
 & \leq & C \mathds{1}_{|\eta|\leq 3/2}  +C \mu^N  \langle \eta \rangle^N e^{-\frac{\mu}{2}(|\eta|-1)^2} \mathds{1}_{|\eta|\geq 3/2} \\
& \leq &C \mathds{1}_{|\eta|\leq 3/2} + C \mu^N \langle \eta \rangle^N e^{-\frac{\mu}{8} |\eta|^2} \mathds{1}_{|\eta|\geq 3/2} \\
 &\leq &C \mathds{1}_{|\eta|\leq 3/2} + C e^{-c \mu} .
\ena
As a consequence, since $m=1$ on a neighborhood of $[0,3/4]$, we obtain, for all $s\geq 0$, the existence of $C, N, c >0$ such that for $\mu \geq \mu_0$, we have
\bna
\langle \zeta \rangle m_\mu(\frac{\zeta}{\alpha \mu}) & \leq & C \mu^N\langle \zeta \rangle^{-s} \mathds{1}_{\frac{|\zeta|}{\alpha \mu}\leq 3/2} + C \langle \zeta \rangle^{-s}  e^{-c \mu} \\
& \leq & C \mu^N\langle \zeta \rangle^{-s} m_\mu(\frac{\zeta}{2\alpha \mu}) + C \langle \zeta \rangle^{-s}  e^{-c \mu} .
\ena
This implies
\bnan
\label{e:two-terms-obs}
\nor{M^{\alpha\mu}_\mu \vartheta_{\mu} v}{1}\leq C \mu^N \nor{M^{2\alpha\mu}_\mu  \vartheta_{\mu} v}{-s}  + C e^{-c \mu} \nor{ \vartheta_{\mu} v}{-s} .
\enan
Concerning the last term, we have 
$$
\nor{\vartheta_{\mu}v}{-s} \leq C \nor{v}{-s} =  C\nor{ \chi^0(t)\chi_{\mu}^1(\frac{D_t}{\mu})\chi_{\mu}^2(t)\chi^3(x)u }{-s} \leq C \nor{\chi_{\mu}^2(t)\chi^3(x)u }{-s} \leq C \nor{U_0}{-s,\times} ,
$$
where, in the last inequality, we used~\eqref{hypoellidual}. 

Concerning the first term in the right hand-side of~\eqref{e:two-terms-obs}, we write
\bnan
\label{e:splitintwo}
\nor{M^{2\alpha\mu}_\mu  \vartheta_{\mu} v}{-s}& = &
\nor{ M^{2\alpha\mu}_\mu  \vartheta_{\mu} \chi^0(t)\chi_{\mu}^1(\frac{D_t}{\mu})\chi_{\mu}^2(t)\chi^3(x)u}{-s} \nonumber \\
 &\leq &\nor{M^{2\alpha\mu}_\mu  \vartheta_{\mu} \chi^0(t)(1-\chi_{\mu}^1(\frac{D_t}{\mu}))\chi_{\mu}^2(t)\chi^3(x)u}{-s}+\nor{ M^{2\alpha\mu}_\mu  \vartheta_{\mu} \chi^0(t)\chi_{\mu}^2(t)\chi^3(x)u}{-s} .
\enan
Concerning the first term in~\eqref{e:splitintwo}, we have 
\bna
&&\nor{M^{2\alpha\mu}_\mu \vartheta_{\mu} \chi^0(t)(1-\chi_{\mu}^1(\frac{D_t}{\mu}))\chi_{\mu}^2(t)\chi^3(x)u}{-s}\\
&\leq& \nor{M^{2\alpha\mu}_\mu \vartheta_{\mu} (1-\chi_{\mu}^1(\frac{D_t}{\mu}))\chi_{\mu}^2(t)\chi^3(x)u}{-s}+\nor{M^{2\alpha\mu}_\mu \vartheta_{\mu} (1-\chi^0(t))(1-\chi_{\mu}^1(\frac{D_t}{\mu}))\chi_{\mu}^2(t)\chi^3(x)u}{-s}\\
&\leq &Ce^{-c\mu}\nor{\chi_{\mu}^2(t)\chi^3(x)u}{-s} ,
\ena
where we have used Lemma \ref{lmmixspcetime} for the first term and then  that $\chi^0(t)=1$ on $\supp \vartheta$ and Lemma \ref{lmintersect} for the second. That this term is bounded by $C  e^{-c\mu}\nor{U_0}{-s,\times}$ then follows from to \eqref{hypoellidual}.

The second term in~\eqref{e:splitintwo} is simpler to handle. Consider $\theta_{\omega}$ a smooth cutoff function supported in $]-T,T[\times \omega$ and equal to one on a neighborhood of $\supp \vartheta$. We have
\bna
\nor{M^{\alpha\mu}_\mu \vartheta_{\mu} \chi^0(t)\chi_{\mu}^2(t)\chi^3(x)u}{-s}
&\leq &\nor{\vartheta_{\mu} \theta_{\omega}\chi^0(t)\chi_{\mu}^2(t)\chi^3(x)u}{-s}+\nor{\vartheta_{\mu}(1-\theta_{\omega}) \chi^0(t)\chi_{\mu}^2(t)\chi^3(x)u}{-s}\\
&\leq&\nor{u}{L^2(]-T,T[\times \omega)}+Ce^{-c\mu }\nor{ \chi^2_\mu(t)\chi^3(x)u}{-s} \\
&\leq & \nor{u}{L^2(]-T,T[\times \omega)}+Ce^{-c\mu }\nor{ U_0}{-s,\times},
\ena
where we have used Lemma \ref{lmintersect} and then again Estimate \eqref{hypoellidual} in the last step.
The last two estimates combined with~\eqref{e:two-terms-obs} and~\eqref{e:splitintwo} conclude the proof of the lemma.
\enp

The proofs of the following two lemmata are based on the spectral representation~\eqref{e:spectral-representation} of the solution $u$ as
\bnan
\label{e:spectral-representation}
u(t,x)=\sum_{j \in \N} \left(a_j^+e^{i \sqrt{\lambda_j}t}+a_j^-e^{-i\sqrt{\lambda_j}t} \right)\varphi_j(x), \quad (t,x) \in \R \times \M .
\enan
This explicit expression allows to prove that a time-frequency cutoff $\chi (\frac{D_t}{\mu})$ truncates also space-frequencies (see in particular the use of Lemma \ref{lmestimint} below).

\begin{lemma}
For all $s \geq 0$, there is $C, c >0$ such that for all $\mu \geq 1$ and for all $u$ solution to~\eqref{hypoelliptic-wave-v} and $v$ defined accordingly in~\eqref{e:def-u-from-v}, we have
\label{lmtermP}
\bna
\nor{Pv}{L^2(\widetilde{\Omega}_{\gamma})}\leq Ce^{-c\mu}\nor{U_0}{-s,\times} .
\ena
\end{lemma}
\bnp[Proof of Lemma~\ref{lmtermP}]
Since $Pu=0$, $\chi^3(x)=1$ on $\Omega_{\gamma}$ and $\chi^0(t)=1$ on $]-T,T[$, we have on $\widetilde{\Omega}_{\gamma}= ]-T,T[ \times \Omega_{\gamma}$
\bnan
\label{e:comm-Pv}
Pv=(\d_t^2+\L)\chi^0(t)\chi_{\mu}^1(\frac{D_t}{\mu})\chi_{\mu}^2(t)\chi^3(x)u=\chi^0(t)\chi_{\mu}^1(\frac{D_t}{\mu})[\d_t^2,\chi_{\mu}^2(t)]\chi^3(x)u , 
\enan
with $[\d_t^2,\chi_{\mu}^2(t)]=(\d_t^2\chi^2)_{\mu}(t) + 2(\partial_t\chi^2)_{\mu}\partial_t$ (the terms $[\chi^0, \d_t^2]$ and $[\chi^3 , \L]$ being supported outside of $\widetilde{\Omega}_{\gamma}$). We only treat the second term in this commutator, the first one being simpler to handle.

We split $u$ given in~\eqref{e:spectral-representation} into high and low frequencies as $u=u_{\leq }+u_{> }$ with 
\bnan
\label{e:def-vleq}
u_{\leq } :=\mathds{1}_{\sqrt{\L}\leq 8 \mu}u = \sum_{\sqrt{\lambda_j}\leq 8 \mu} \cdots ,  \quad \quad \quad 
 u_{> } := \mathds{1}_{\sqrt{\L}> 8 \mu}u = \sum_{\sqrt{\lambda_j} > 8 \mu} \cdots.
\enan
We also write $f(t)=\partial_t\chi^2(t)$ and $f_\mu(t)=(\partial_t\chi^2)_\mu(t)$.
We first estimate the low frequencies as
\bna
\nor{\chi^0(t)\chi_{\mu}^1(\frac{D_t}{\mu}) f_{\mu}\chi^3(x)\partial_tu_{\leq }}{L^2(\widetilde{\Omega}_{\gamma})}
&\leq& C\nor{\mathds{1}_{]-T,T[}\chi_{\mu}^1(\frac{D_t}{\mu}) f_{\mu} \chi^3(x)\partial_tu_{\leq }}{L^2(\R_t\times \Omega_\gamma)}\\
&\leq& C\nor{\mathds{1}_{]-T,T[}\chi_{\mu}^1(\frac{D_t}{\mu}) f_{\mu} \chi^0_\mu(t) \chi^3(x)\partial_tu_{\leq }}{0}\\
&& + C\nor{\mathds{1}_{]-T,T[}\chi_{\mu}^1(\frac{D_t}{\mu}) f_{\mu} (1-\chi^0_\mu(t))\chi^3(x)\partial_tu_{\leq }}{0} .
\ena
Concerning the fist term in the right hand-side, we have
\bna
\nor{\mathds{1}_{]-T,T[}\chi_{\mu}^1(\frac{D_t}{\mu}) f_{\mu} \chi^0_\mu(t) \chi^3(x)\partial_tu_{\leq }}{L^2(\R_t\times \Omega_\gamma)}
&\leq & \nor{\mathds{1}_{]-T,T[}\chi_{\mu}^1(\frac{D_t}{\mu}) f_{\mu} }{L^2(\R_t)\mapsto L^2(\R_t)} \nor{\chi^0_\mu(t) \chi^3(x)\partial_tu_{\leq}}{0} \\
&\leq &Ce^{-c \mu} \nor{(u_{\leq}, \d_t u_{\leq})|_{t=0}}{1,\times} \leq Ce^{-c \mu} \mu^{s}\nor{U_0}{-s,\times} ,
\ena
after having used Lemma \ref{lmtrucQUCP} together with~\eqref{hypoellidual-bis}.
Concerning the second term in the right hand-side, we write
$$\nor{\mathds{1}_{]-T,T[}\chi_{\mu}^1(\frac{D_t}{\mu}) f_{\mu}(1- \chi^0_\mu(t))\chi^3(x)\partial_tu_{\leq }}{L^2(\R_t\times \Omega_\gamma)} 
\leq  \nor{\mathds{1}_{]-T,T[}\chi_{\mu}^1(\frac{D_t}{\mu})(1- \chi^0_\mu(t))}{L^2(\R_t)\mapsto L^2(\R_t)} \nor{f_\mu(t) \chi^3(x)\partial_tu_{\leq}}{0} 
$$and we conclude with the same arguments. Hence, we have the low frequency estimate
\bnan
\label{e:low-frq-j}
\nor{\chi^0(t)\chi_{\mu}^1(\frac{D_t}{\mu}) f_{\mu}\chi^3(x)\partial_tu_{\leq }}{L^2(\widetilde{\Omega}_{\gamma})} \leq Ce^{-c \mu} \nor{U_0}{-s,\times} .
\enan
Concerning now the high-frequency part, still denoting $f=\partial_t\chi^2$, we have
\bna
\nor{\chi^0(t)\chi_{\mu}^1(\frac{D_t}{\mu})f_{\mu}(t)\chi^3(x)\partial_tu_{>}}{L^2(\widetilde{\Omega}_{\gamma})}\leq \nor{\chi_{\mu}^1(\frac{D_t}{\mu})f_{\mu}(t)\chi^3(x)\partial_tu_{>}}{L^2(\R_t\times \M)} ,
\ena
where 
\bna
\chi_{\mu}^1(\frac{D_t}{\mu})f_{\mu}(t)\chi^3(x)\partial_tu_{>} = 
\sum_{\sqrt{\lambda_j} > 8 \mu} i \sqrt{\lambda_j} \left[\chi_{\mu}^1(\frac{D_t}{\mu})f_{\mu}(t) \left(a_j^+e^{i \sqrt{\lambda_j}t} -a_j^-e^{-i\sqrt{\lambda_j}t} \right) \right]  \left[\chi^3 \varphi_j \right](x) .
\ena
As a consequence of the triangular inequality, the Plancherel theorem and the fact that $(\varphi_j)$ is an orthonormal family in $L^2(\M)$, we obtain
\bna
&& \nor{\chi^0(t)\chi_{\mu}^1(\frac{D_t}{\mu})f_{\mu}(t)\chi^3(x)\partial_tu_{>}}{L^2(\widetilde{\Omega}_{\gamma})} \\
& & \leq \sum_{\sqrt{\lambda_j}\geq 8 \mu} \sqrt{\lambda_j}\nor{\chi_{\mu}^1(\frac{D_t}{\mu})f_{\mu}(t) \left(a_j^+e^{i \sqrt{\lambda_j}t} -a_j^-e^{-i\sqrt{\lambda_j}t} \right)}{L^2(\R_t)}\nor{\chi^3(x)\varphi_j}{L^2(\M)}\\
&& \leq \sum_{\sqrt{\lambda_j}\geq 8 \mu} \sqrt{\lambda_j}\nor{ \chi_{\mu}^1\left(\frac{\tau}{\mu}\right) \left(\widehat{f_{\mu}}(\tau -\sqrt{\lambda_j})a_j^+ -\widehat{f_{\mu}}(\tau +\sqrt{\lambda_j})a_j^-\right)}{L^2(\R_\tau)} \\
&& \leq \sum_{\sqrt{\lambda_j}\geq 8\mu} \sqrt{\lambda_j} |a_j^+| \nor{ \chi_{\mu}^1\left(\frac{\tau}{\mu}\right) \widehat{f_{\mu}}(\tau -\sqrt{\lambda_j}) }{L^2(\R_\tau)}
+ \sqrt{\lambda_j} |a_j^-| \nor{  \chi_{\mu}^1\left(\frac{\tau}{\mu}\right) \widehat{f_{\mu}}(\tau +\sqrt{\lambda_j})}{L^2(\R_\tau)} .
\ena 
Using now Lemma \ref{lmestimint}, this yields
\bna
 \nor{\chi^0(t)\chi_{\mu}^1(\frac{D_t}{\mu})f_{\mu}(t)\chi^3(x)\partial_tu_{>}}{L^2(\widetilde{\Omega}_{\gamma})} 
   & \leq & \sum_{\sqrt{\lambda_j}\geq 8 \mu} \sqrt{\lambda_j} C e^{-c \sqrt{\lambda_j}} (|a_j^+| +  |a_j^-| )\\
    &  \leq & C \Big( \sum_{\sqrt{\lambda_j}\geq 8 \mu} e^{-c \sqrt{\lambda_j}} (|a_j^+|^2 +  |a_j^-|^2 ) \Big)^{1/2}
\leq Ce^{-c\mu}\nor{U_0}{-s,\times} ,
\ena
after having used the Cauchy-Schwarz inequality in $\ell^2(\N)$. This together with~\eqref{e:low-frq-j} and~\eqref{e:comm-Pv} concludes the proof of the lemma.
\enp

The following Lemma will be used to estimate the last term in~\eqref{inegtriang}. 
\begin{lemma}
\label{lmtermreste}
For all $s \geq 0$, there is $C, N >0$ such that for all $\mu \geq 1$ and for all $u$ solution to~\eqref{hypoelliptic-wave-v} and $v$ defined accordingly in~\eqref{e:def-u-from-v}, we have
\bna
\nor{v}{1}\leq C\mu^{N}\nor{U_0}{-s,\times} .
\ena
\end{lemma}
\bnp
We first write
$$
\nor{v}{1} \leq \nor{\chi_{\mu}^1(\frac{D_t}{\mu})\chi^2_{\mu}(t)\chi^3(x)u}{H^1(\R_t\times \M)} ,
$$ and decompose $u$ in~\eqref{e:spectral-representation} as  $u=u_{\leq}+u_{>}$ with $u_{\leq} = \mathds{1}_{\sqrt{\L}\leq 8 \mu}u$ and $u_{>} =\mathds{1}_{\sqrt{\L}> 8\mu}u$ being defined as in~\eqref{e:def-vleq}. We have
\bna
\nor{\chi_{\mu}^1(\frac{D_t}{\mu})\chi^2_{\mu}(t)\chi^3(x)u}{H^1(\R_t\times \M)}
\leq \nor{\chi_{\mu}^1(\frac{D_t}{\mu})\chi^2_{\mu}(t)\chi^3(x)u_{\leq}}{H^1(\R_t\times \M)}
+\nor{\chi_{\mu}^1(\frac{D_t}{\mu})\chi^2_{\mu}(t)\chi^3(x)u_{>}}{H^1(\R_t\times \M)} .
\ena
Concerning the first term (low-frequencies), we simply use \eqref{hypoellidual-ter} to write
\bna
\nor{\chi_{\mu}^1(\frac{D_t}{\mu})\chi^2_{\mu}(t)\chi^3(x)u_{\leq}}{H^1(\R_t\times \M)}
\leq C  \nor{\chi^2_{\mu}(t)\chi^3(x)u_{\leq}}{H^1(\R_t\times \M)}
\leq C \nor{(u_{\leq}, \d_t u_{\leq})|_{t=0}}{k,\times} \leq C \mu^{k+s}\nor{U_0}{-s,\times} .
\ena
Concerning the second term (high-frequencies), we use $\|f(t)g(x)\|_{H^1(\R\times \M)} \leq C \|f\|_{H^1(\R)} \|g\|_{H^1(\M)}$ and proceed as in the above proof of Lemma~\ref{lmtermP}. This yields
\bna
&& \nor{\chi_{\mu}^1(\frac{D_t}{\mu})\chi^2_{\mu}(t)\chi^3(x)u_{>}}{H^1(\R_t\times \M)} \\
&\leq & C \sum_{\sqrt{\lambda_j}\geq 8 \mu} \nor{\left\langle \tau\right\rangle  \chi_{\mu}^1\left(\frac{\tau}{\mu}\right) \left(\widehat{\chi^2_{\mu}}(\tau -\sqrt{\lambda_j})a_j^+ +\widehat{\chi_{\mu}^2}(\tau +\sqrt{\lambda_j})a_j^-\right)}{L^2(\R_\tau)}\nor{\chi^3(x)\varphi_j}{H^1(\M)}\\
&\leq& C \sum_{\sqrt{\lambda_j}\geq 8 \mu} \lambda_j^{N(k,s)}\nor{\chi_{\mu}^1\left(\frac{\tau}{\mu}\right) \left(\widehat{\chi^2_{\mu}}(\tau -\sqrt{\lambda_j})a_j^+ +\widehat{\chi^2_{\mu}}(\tau +\sqrt{\lambda_j})a_j^-\right)}{L^2(\R_\tau)}\\
&\leq & C \sum_{\sqrt{\lambda_j}\geq 8 \mu}Ce^{- c\sqrt{\lambda_j}}\left(|a_j^+|+|a_j^-|\right)\leq C \nor{U_0}{-s,\times} ,
\ena
where we have used $\| \varphi_j\|_{H^s(\M)}\leq C (\lambda_j+1)^{\frac{ks}{2}}$ (which follows directly from Corollary~\ref{cor:HsHsL}), Lemma \ref{lmestimint} and the Cauchy-Schwarz inequality in $\ell^2(\N)$.
\enp
Before concluding the proof of Theorem \ref{thmwavehypo-s}, we need to explain how to estimate the high-frequency part of the solution (recall indeed that our starting point, Estimate~\eqref{inegtriang}, is a low-frequency estimate only). This is the aim of the following lemma, which proof is close to the proof that Theorem~4.7 implies Theorem~1.11 in~\cite{LL:15}.
\begin{lemma}
\label{lmhighweak}
For all $s\in [0,k]$, there is $C, \mu_0 >0$ such that for all $\mu \geq \mu_0$ and for all $u$ solution to~\eqref{hypoelliptic-wave-v}, we have
\bna
\nor{(1-M_{\mu}^{\beta\mu}) \widetilde{\vartheta}_{\mu} \chi^0(t) \chi_{\mu}^2(t)\chi^3(x)u}{L^2(\R_t\times \R^d_x)}\leq \frac{C}{\mu^{s/k}}\nor{U_0}{s,\times} .
\ena
\end{lemma}
Note that this estimate is almost 
\bna
\nor{(1-M^{\beta\mu}_{\mu})\widetilde{\vartheta}_{\mu}v}{L^2(\R_t\times \R^d_x)}\leq \frac{C}{\mu^{s/k}}\nor{U_0}{s,\times} .
\ena
(which is also true, but not used here), the difference being that $v$ contains an additional time-frequency cutoff $\chi_{\mu}^1(\frac{D_t}{\mu})$ (which does not play any role in the estimates below).

The proof below only gives the endpoint case $s=k$, the intermediate situations being deduced by interpolation. A direct proof of intermediate estimates would follow the same lines, yet being slightly longer.
\bnp
We first notice that it is enough to treat the case $s=k$. Indeed, $s=0$ is direct by standard energy estimates for $v$ (see for instance \eqref{e:gronwall}) and uniform bound on $\chi_{\mu}^2$ and $\widetilde{\vartheta}_{\mu}$. Hence, since all operators involved are linear, the result for $s\in [0,k]$ follows by interpolation, see for instance \cite[Chapter 23]{TartarBookinterp}.

Concerning the case $s=k$, we write $w =  \chi^0(t) \chi_{\mu}^2(t)\chi^3(x)u$, together with
\bna
\nor{(1-M^{\beta\mu}_{\mu}) \widetilde{\vartheta}_{\mu}w}{0}\leq C \sup_{\eta \in \R^{d+1}}\left|\left\langle \eta\right\rangle^{-1}(1-m_{\mu})(\frac{\eta}{\beta\mu})\right| \nor{\widetilde{\vartheta}_{\mu}w}{1}.
\ena
In the range $|\eta|\geq \beta\mu/2$ with $\mu\geq \mu_0$, we have the loose estimate $\left|\left\langle \eta\right\rangle^{-1}(1-m_{\mu})(\frac{\eta}{\beta\mu})\right| \leq \frac{C}{\mu}$
whereas in the range $|\eta|\leq \beta\mu/2$, using $\dist\left( \supp (1-m(\frac{\cdot}{\beta})), \left\{|\eta|\leq \beta/2\right\}\right)>0$, we have $\left|(1-m_{\mu})(\frac{\zeta}{\beta\mu})\right|\leq Ce^{-c\mu}$ according to estimate \eqref{Estimate29} below. 
This implies, for $\mu\geq \mu_0$ that 
\bnan
\label{refref-1}
\nor{(1-M^{\beta\mu}_{\mu})\widetilde{\vartheta}_{\mu}w}{L^2(\R_t\times \R^d_x)}\leq \frac{C}{\mu}\nor{\widetilde{\vartheta}_{\mu}w}{H^{1}(\R_t\times \R^d_x)} .
\enan
Next, we have that 
\bnan
\label{refref-2}
\nor{\widetilde{\vartheta}_{\mu}w}{H^{1}(\R_t\times \R^d_x)}\leq  \nor{w}{H^{1}(\R_t\times \R^d_x)} \leq  \nor{\chi_{\mu}^2(t)\chi^3(x)u}{H^{1}(\R_t\times \R^d_x)} .
\enan
uniformly for $\mu \geq 1$, since all derivatives of $\widetilde{\vartheta}_{\mu}$ are uniformly bounded for $\mu \geq 1$.
We write for $\tilde{w} = \chi_{\mu}^2(t)\chi^3(x)u$
\bna
\nor{\tilde{w}}{H^{1}(\R_t\times \R^d_x)}\leq \nor{\tilde{w}}{H^{1}(\R_t; L^2(\R^d_x))}+ \nor{\tilde{w}}{L^2(\R_t; H^{1}(\R^d_x))} ,
\ena
and estimate each term separately.
Concerning the first term in this estimate, we have
\bna
\nor{\tilde{w}}{H^{1}(\R_t; L^2(\R^d_x))} & \leq & \nor{\tilde{w}}{L^2(\R_t\times \R^d_x)}+\nor{\partial_t \tilde{w}}{L^2(\R_t\times \R^d_x)} \\
&\leq& \nor{\chi_{\mu}^2(t)\chi^3(x)u}{L^2(\R_t\times \R^d_x)} + \nor{(\d_t\chi^2)_{\mu}(t)\chi^3(x)u}{L^2(\R_t\times \R^d_x)} + \nor{\chi_{\mu}^2(t)\chi^3(x) \d_t u}{L^2(\R_t\times \R^d_x)} \\
&\leq&   C\nor{U_0}{k,\times}.
\ena
after having used~\eqref{hypoellidual}, ~\eqref{hypoellidual-bis} and $k\geq 1$.
Similarly, the second term is estimated as
\bna
\nor{\tilde{w}}{L^2(\R_t; H^{1}(\R^d_x))} 
=  \nor{\chi_{\mu}^2(t)\chi^3(x)u}{L^2(\R_t; H^{1}(\R^d_x))} 
\leq C\nor{U_0}{k,\times} ,
\ena
as a direct consequence of~\eqref{hypoellidual-ter}. The above three estimates together with \eqref{refref-1} and~\eqref{refref-2} conclude the proof of the lemma.
\enp

With the above four lemmata in hand, we can now conclude the proof of Theorem \ref{thmwavehypo-s}.
\bnp[Proof of Theorem \ref{thmwavehypo-s}] We first prove the result for $s\in ]0,k]$, the conclusion for all $s>0$ being then a consequence of Remark~\ref{rem:Hsspaces}.

Starting form Estimate~\eqref{inegtriang}, combined with Lemmata \ref{lmtermobs}, \ref{lmtermP} and \ref{lmtermreste} to bound the terms in the right hand side, we first obtain the intermediate estimate 
\bnan
\label{e:rhs-done}
\nor{M_{\mu}^{\beta\mu} \widetilde{\vartheta}_{\mu}  v}{1}\leq Ce^{c\mu}\nor{u}{L^2(]-T,T[\times \omega)}+ Ce^{-c'\mu}\nor{U_0}{-s,\times}.
\enan
Note that in order to obtain this inequality, we have chosen $\kappa$ in \eqref{inegtriang} to be small enough compared to the constant $c$ appearing in Lemmata \ref{lmtermobs} and \ref{lmtermP}.
Now, recalling that $v=  \chi^0(t)\chi_{\mu}^1(\frac{D_t}{\mu})\chi_{\mu}^2(t)\chi^3(x)u$, we decompose 
\bnan 
\label{decomp-finale}
\widetilde{\vartheta}_{\mu} \chi^0(t) \chi_{\mu}^2(t)\chi^3(x)u
&= &M_{\mu}^{\beta\mu} \widetilde{\vartheta}_{\mu}  v 
+ M_{\mu}^{\beta\mu} \widetilde{\vartheta}_{\mu} \chi^0(t)\left(1-\chi_{\mu}^1(\frac{D_t}{\mu})\right)\chi_{\mu}^2(t)\chi^3(x)u \nonumber \\
&&+ (1-M_{\mu}^{\beta\mu}) \widetilde{\vartheta}_{\mu} \chi^0(t) \chi_{\mu}^2(t)\chi^3(x)u .
\enan
The first term is estimated in~\eqref{e:rhs-done}, the last one is estimated in Lemma \ref{lmhighweak} so that it only remains to estimate the second one. We have
\bna
M_{\mu}^{\beta\mu} \widetilde{\vartheta}_{\mu} \chi^0(t)\left(1-\chi_{\mu}^1(\frac{D_t}{\mu})\right)\chi_{\mu}^2(t)\chi^3(x)u
& = &M_{\mu}^{\beta\mu} \widetilde{\vartheta}_{\mu}\left(1-\chi_{\mu}^1(\frac{D_t}{\mu})\right)\chi_{\mu}^2(t)\chi^3(x)u \\
&&+ M_{\mu}^{\beta\mu} \widetilde{\vartheta}_{\mu}(1- \chi^0(t))\left(1-\chi_{\mu}^1(\frac{D_t}{\mu})\right)\chi_{\mu}^2(t)\chi^3(x)u .
\ena
According to Lemma~\ref{lmmixspcetime} and~\eqref{hypoellidual} below, the first of these two terms satisfies
\bna
\nor{M_{\mu}^{\beta\mu} \widetilde{\vartheta}_{\mu}\left(1-\chi_{\mu}^1(\frac{D_t}{\mu})\right)\chi_{\mu}^2(t)\chi^3(x)u }{0}
 \leq Ce^{-c\mu} \nor{\chi_{\mu}^2(t)\chi^3(x)u }{0} \leq  Ce^{-c\mu} \nor{U_0}{0,\times} ,
\ena
whereas, according to Item~2 of Lemma \ref{Lemma23}, we have $\nor{\widetilde{\vartheta}_{\mu}(1- \chi^0(t))}{L^2(\R^n)\to L^2(\R^n)} \leq C e^{-c\mu}$ so that the second one is bounded as well as
\bna
\nor{M_{\mu}^{\beta\mu} \widetilde{\vartheta}_{\mu}(1- \chi^0(t))\left(1-\chi_{\mu}^1(\frac{D_t}{\mu})\right)\chi_{\mu}^2(t)\chi^3(x)u}{0}
 \leq Ce^{-c\mu} \nor{\chi_{\mu}^2(t)\chi^3(x)u }{0} \leq  Ce^{-c\mu} \nor{U_0}{0,\times} .
\ena
Coming back to the decomposition~\eqref{decomp-finale}, using the above two estimates together with~\eqref{e:rhs-done} and Lemma \ref{lmhighweak}, we now have (for $s \in ]0,k]$)
\bna
\nor{\widetilde{\vartheta}_{\mu} \chi^0(t) \chi_{\mu}^2(t)\chi^3(x)u}{L^2(\R_t\times \R^d_x)}\leq Ce^{\kappa\mu}\nor{u}{L^2(]-T,T[\times \omega)}+\frac{C}{\mu^{s/k}}\nor{U_0}{s,\times}.
\ena
Now, since $\widetilde{\vartheta}=1$ in a neighborhood of $]-\e,\e[\times B(x_0,\e)$, for $\mu$ large enough, we have $\widetilde{\vartheta}_{\mu}\geq 1/2$ on $]-\e,\e[\times B(x_0,\e)$, and the same also holds for $ \chi_{\mu}^2(t)$. Hence, we obtain
\bna
\nor{u}{L^2(]-\e,\e[\times B(x_0,\e))}\leq Ce^{\kappa\mu}\nor{u}{L^2(]-T,T[\times \omega)}+\frac{C}{\mu^{s/k}}\nor{U_0}{s,\times}.
\ena
Since $x_0\in \M$ is arbitrary, a compactness argument allows to obtain other constants still denoted $C, \kappa, \mu_0,\e>0$ such that for $\mu \geq \mu_0$, 
\bna
\nor{u}{L^2(]-\e,\e[\times \M)}\leq Ce^{\kappa\mu}\nor{u}{L^2(]-T,T[\times \omega)}+\frac{C}{\mu^{s/k}}\nor{U_0}{s,\times}.
\ena
We conclude the proof of Theorem \ref{thmwavehypo-s}, in the case $s\in ]0,k]$, by using Estimate \eqref{lowerL2} and changing $\mu$ into $\mu^{k}$.
In the case $s>k$, the proof follows from the estimate with $s=k$ and an interpolation argument, as explained in Remark~\ref{rem:Hsspaces}.
\enp

\subsubsection{Technical lemmata used in the previous section (only)}
In this section, we collect some technical results used in the above Section~\ref{s:energy-general-case} (and in that section only). We first state results that are either directly taken from~\cite{LL:15}, or direct consequences of these. Second, we prove three lemmata using these results.

\paragraph{Estimates taken from \cite{LL:15}.}
The following estimate is \cite[Equation~(2.9)]{LL:15} in a simple situation: given a continuous function $f$ on $\R^{n_a}$, we have, for all $\zeta \in \R^{n_a}$
\bnan
\label{Estimate29}
\left|f_{\lambda}(\zeta)\right|
\leq  C \left\langle \lambda\right\rangle^{\frac{n_a}{2}}
\nor{f}{L^\infty} 
\left\langle \dist(\zeta , \supp(f) ) \right\rangle^{n_a-1}
e^{-\frac{\lambda}{4} \dist(\zeta , \supp(f) )^2} 
\enan

\begin{lemma}[Lemma 2.3 of \cite{LL:15}]
\label{Lemma23}
The following three statements holds true.
\begin{enumerate}
\item
For any $d>0$, there exist $C,c>0$ such that for any $f_1, f_2\in L^{\infty}(\R^n)$ such that $\dist(\supp(f_1),\supp(f_2))\geq d$ and all $\lambda\geq 0$, we have 
\bna
\nor{f_{1,\lambda}f_2 }{L^{\infty}} \leq C e^{-c \lambda}\nor{f_1}{L^{\infty}} \nor{f_2}{L^{\infty}} ,
\qquad 
\nor{f_{1,\lambda}f_{2,\lambda}}{L^{\infty}}\leq C e^{-c \lambda}\nor{f_1}{L^{\infty}} \nor{f_2}{L^{\infty}} .
\ena
\item If moreover $f_1, f_2\in C^{\infty}(\R^n)$ have bounded derivatives, then for all $k \in \N$, there exist $C,c>0$ such that for all $\lambda\geq 1$, we have
\bna
\nor{f_{1,\lambda}f_2 }{H^k(\R^n) \to H^k(\R^n) } \leq C e^{-c \lambda} .
\ena
\item 
Let $f_1, f_2\in L^{\infty}(\R^{n_a})$ such that $\dist(\supp(f_1),\supp(f_2))>0$ . Then there exist $C,c>0$ such that for all $\lambda\geq 1$, for all $k \in \N$, for all $\mu \geq 1$,  we have 
\bna
\nor{f_{1,\lambda}(D_a /\mu )f_2 (D_a/\mu)}{H^k(\R^n) \to H^k(\R^n) } \leq C e^{-c \lambda} , \\
\nor{f_{1,\lambda}(D_a /\mu )f_{2,\lambda} (D_a/\mu)}{H^k(\R^n) \to H^k(\R^n) } \leq C e^{-c \lambda} .
\ena
\end{enumerate}
\end{lemma}

\begin{lemma}[Lemma~2.9 of \cite{LL:15}]
\label{Lemma29}
Let $k\in \N$ and $f \in C^\infty_0(\R^n)$. 
Then, there exist $C,c$ such that, for any $\lambda, \mu>0$, we have
\bna
\nor{M^{\mu}_{\lambda}f_{\lambda}(1-M^{2\mu}_{\lambda})}{H^k(\R^n) \to H^k(\R^n)}  \leq C e^{-c \frac{\mu^2}{\lambda}} +C e^{-c \lambda} ; \\
\nor{(1-M^{2\mu}_{\lambda})f_{\lambda}M^{\mu}_{\lambda}}{H^k(\R^n) \to H^k(\R^n)}  \leq C e^{-c \frac{\mu^2}{\lambda}} +C e^{-c \lambda} .
\ena
\end{lemma}

\begin{lemma}
\label{lmintersect}
Let $f_1 , f_2 \in C^\infty(\R)$ with all derivatives bounded and such that $\supp( f_1)\cap \supp (f_2)=\emptyset$. Then, for any $s \in \N$, there is $C, c >0$ such that for all $w \in H^{-s}(\R)$, we have
\bna
\nor{f_1 f_{2,\mu}w}{-s}\leq Ce^{-c\mu }\nor{w}{-s} .
\ena
\end{lemma}
Lemma~\ref{lmintersect} is obtained by duality from Item~2 of Lemma \ref{Lemma23}.

\begin{lemma}
\label{lmtrucQUCP}
Let $f \in C^\infty(\R)$ bounded such that $\supp(f) \cap [-T,T] = \emptyset$. Then, there is $C,c>0$ such that 
$$\nor{\mathds{1}_{]-T,T[} \chi_{\mu}^1(D_t/\mu) f_{\mu}}{L^2(\R_t)\mapsto L^2(\R_t)}\leq Ce^{-c\mu} . $$
\end{lemma}
Lemma~\ref{lmtrucQUCP} is a particular case of \cite[Lemma~2.10]{LL:15}.

\paragraph{A few estimates using the above lemmata.}
Recalling the definition of $\E_s$ in~\eqref{def-energy} we now refine the rough Estimate~\eqref{e:gronwall}. Indeed, on account to the spectral theory of $\L$, if we denote $\Pi_0$ the spectral projector on $\ker(\L)=  \vect_{L^2}(1)$, we have, for all $u \in \H^s_\L$,
$$
\nor{u}{\H_\L^s} \simeq \nor{\L^{\frac{s}{2}}u}{L^2(\M)} + |\Pi_0 u | .
$$
Notice now that the energy $\nor{\L^{\frac{s}{2}} u}{L^2}^2 +\nor{\L^{\frac{s-1}{2}}\d_t u}{L^2}^2$ is preserved by the equation~\eqref{hypoelliptic-wave-v}, and that the equation for the zero frequency is 
$\d_t^2 \Pi_0 u= 0$ hence growing at most linearly. As a consequence, we finally obtain that for all $s\in \R$, there is $C>0$ such that for all solution $u$ of~\eqref{hypoelliptic-wave-v}, we have
\bnan
\label{e:gronwall-refined}
 \E_{s}(u)(t) dt\leq C(1+|t|) \E_{s}(u)(0) , \quad \text{ for all }t\in \R .
\enan
This estimate is now used to bound some integrals of $\E_{s}(u)$.
\begin{lemma}
Let $\eta \in C^\infty_0(\R)$, denote $\chi^3(x)$ as above and fix $s \geq 0$. Then, there is $C>0$ such that for all $\mu \geq 1$ and  for all solution $u$ of~\eqref{hypoelliptic-wave-v}, we have
\bnan
\label{hypoellidual}
\nor{ \eta_\mu (t)\chi^3(x)u}{-s}\leq C\nor{U_0}{-s,\times} , \\
\label{hypoellidual-bis}
 \nor{ \eta_\mu (t)\chi^3(x)\d_t u}{0}\leq C\nor{U_0}{1,\times} , \\
\label{hypoellidual-ter}
\nor{ \eta_\mu (t)\chi^3(x)u}{1} 
\leq C\nor{U_0}{k,\times} , \\
\label{hypoellidual-quad}
\nor{ \eta_\mu (t)\chi^3(x)u}{L^2(\R_t; H^{\frac{s}{k}}(\R^d_x))}
\leq C\nor{U_0}{s,\times} .
\enan

\end{lemma}
\bnp
We first remark (see e.g.~\eqref{Estimate29}) that there is $C,c>0$ such that for $\mu \geq 1$, we have $0\leq \eta_\mu (t) \leq C e^{-c|t|}$ for all $t \in \R$.
To prove~\eqref{hypoellidual}, we now simply write 
\bna
\nor{ \eta_\mu (t)\chi^3(x)u}{-s}^2 &\leq & \int_\R C e^{-c|t|} \nor{\chi^3(x)u}{H^{-s}(\Omega_\gamma)}^2 dt =  \int_\R C e^{-c|t|} \nor{\chi^3(x)u}{H^{-s}(\M)}^2 dt \\
 &\leq & \int_\R C e^{-c|t|} \nor{u}{H^{-s}(\M)}^2 dt \leq   \int_\R C e^{-c|t|} \nor{u}{\H_\L^{-s}}^2 dt ,
\ena
where we used~\eqref{e:sobolev-simple-dual} and $s\geq 0$ in the last inequality. Recalling the definition of $\E_s$ in~\eqref{def-energy} together with estimate~\eqref{e:gronwall-refined}, we now have
\bna
\nor{ \eta_\mu(t)\chi^3(x)u}{-s}^2  \leq \int_\R C e^{-c|t|} \E_{-s}(u)(t) dt\leq C \left(\int_\R e^{-c|t|}(1+|t|) dt \right)\E_{-s}(u)(0) = C \nor{U_0}{\H_\L^{-s}\times\H_\L^{-s-1}}^2 ,
\ena
which concludes the proof of~\eqref{hypoellidual}. 
The proof of~\eqref{hypoellidual-bis} is the same, except that we use $\nor{\d_t u}{L^2(\M)}^2 \leq  2\E_{1}(u)$ instead of~\eqref{e:sobolev-simple-dual}.
The proof of~\eqref{hypoellidual-ter} is similar: after using the chain rule, each term is either of the form $\eta_\mu (t)\chi^3(x) u$ (treated in~\eqref{hypoellidual}) or of the form $\eta_\mu (t)\chi^3(x) \d_t u$ (treated in~\eqref{hypoellidual-bis}) or of the form $\eta_\mu (t)\chi^3(x) \d_x u$, for which the proof is the same using $\nor{ u}{H^{\frac{s}{k}}(\M)}^2 \leq  2\E_{s}(u)$, consequence of Corollary~\eqref{cor:HsHsL}, instead of~\eqref{e:sobolev-simple-dual}. The proof of~\eqref{hypoellidual-quad} is the same, still using Corollary~\eqref{cor:HsHsL}.
\enp

\begin{lemma}
\label{lmmixspcetime}
Fix $\alpha<1$. Let $\chi^1 \in C^{\infty}_0(\R)$ such that $\chi^1(\tau)=1$ for $\tau \in ]-1,1[$ and $\vartheta \in C^{\infty}_0(\R^{1+d})$.
Then, for any $s \in \N$, there is $C, c >0$ such that for all $w \in H^{-s}(\R^{1+d})$, we have
\bna
\nor{M^{\alpha\mu}_\mu \vartheta_{\mu} (1-\chi_{\mu}^1(\frac{D_t}{\mu}))w}{-s}\leq Ce^{-c\mu }\nor{w}{-s} .
\ena
\end{lemma}
Note that in the proofs above, the parameter $\alpha$ of Lemma~\ref{lmmixspcetime} is both taken to be the parameters $\alpha$ or $\beta$ appearing estimate~\eqref{inegtriang}. This is the reason why we assumed $\alpha, \beta <1$ there.
\bnp
This lemma is very close to (and a consequence of) Lemma \ref{Lemma29} except that $\chi_{\mu}^1(D_t)$ is a Fourier cutoff in $D_t$ only whereas $M^{\alpha\mu}_\mu$ are Fourier cutoffs in the whole $D_{t,x}$ (and that the Sobolev orders are negative).
Recall that $M^{\alpha\mu}_\mu  = m_{\mu}\left( \frac{D}{\alpha \mu}\right)$, $D=D_{t,x}$, where $m$ is compactly supported in $|\xi| < 1$ and $m(\xi)=1$ for $|\xi| < 3/4$ (see the beginning of Section~\ref{s:def-mult}).

Let $\tilde{m}$ be a radial smooth function on $\R^{1+d}$ such that $\tilde{m}(\xi) = 1$ in a neighborhood of $|\xi| \leq \alpha$ and $\tilde{m}(\xi) = 0$ in a neighborhood of $|\xi| \geq 1$. Then we have $1-\chi^1=0$ on the support of $\tilde{m}$.
For $s\in \N$, we write
\bna
\nor{(1-\chi_{\mu}^1(D_t)) \vartheta_{\mu}  M^{\alpha\mu}_\mu w}{s}
& \leq &  \nor{(1-\chi_{\mu}^1(D_t/\mu))(1-\tilde{m}_{\mu}(D/\mu)) \vartheta_{\mu}  M^{\alpha\mu}_\mu w}{s}\\
&& +\nor{(1-\chi_{\mu}^1(D_t))\tilde{m}_{\mu}(D/\mu) \vartheta_{\mu}  M^{\alpha\mu}_\mu w}{s} \\
& \leq & C \nor{(1-\tilde{m}_{\mu}(D/\mu)) \vartheta_{\mu}  M^{\alpha\mu}_\mu w}{s}
+C\nor{(1-\chi_{\mu}^1(D_t))\tilde{m}_{\mu}(D/\mu) \vartheta_{\mu}  M^{\alpha\mu}_\mu w}{s} .
\ena
According to Lemma \ref{Lemma29} and the respective supports of $\tilde{m}$ and $m(\frac{\cdot}{\alpha})$, we have
$$
\nor{(1-\tilde{m}_{\mu}(D/\mu)) \vartheta_{\mu}  M^{\alpha\mu}_\mu w}{s} \leq C e^{-c\mu }\nor{w}{s} .
$$
Also, according to Item~3 of \ref{Lemma23},  and the respective supports of $\tilde{m}$ and $\chi^1$, we have
$$
\nor{(1-\chi_{\mu}^1(D_t/\mu))\tilde{m}_{\mu}(D/\mu)}{H^s \to H^s} \leq  C e^{-c\mu },
$$
and hence 
$$
\nor{(1-\chi_{\mu}^1(D_t/\mu))\tilde{m}_{\mu}(D/\mu) \vartheta_{\mu}  M^{\alpha\mu}_\mu w}{s} \leq  C e^{-c\mu }\nor{w}{s}
$$
This finally yields for $s \in \N$
$$
\nor{(1-\chi_{\mu}^1(D_t/\mu)) \vartheta_{\mu}  M^{\alpha\mu}_\mu w}{s} \leq  C e^{-c\mu }\nor{w}{s},
$$
and the sought estimate by a duality argument.
\enp

\begin{lemma}
\label{lmestimint}
Let $\chi\in C^{\infty}_0(\R)$ and $m\in C^{\infty}_0(]-1,1[)$, and define
$$f_{\mu,\lambda}(\tau)=m_{\mu}\left(\frac{\tau}{\mu}\right) \widehat{\chi_\mu}(\tau-\lambda) .
$$
Then, for all $\sigma \in \R$, there is $C,c>0$ so that we have 
\bnan
\label{estimflarge}
\nor{\langle \tau \rangle^\sigma f_{\mu,\lambda}(\tau)}{L^2} \leq  Ce^{-c|\lambda|} , \quad  \textnormal{ for all } \lambda\in \R , \ \mu\geq 0 \ \text{ such that }|\lambda|\geq 4\mu  .
\enan
\end{lemma}
\bnp
We decompose $f_{\mu,\lambda}=f_{\mu,\lambda}^1+f_{\mu,\lambda}^2$ with $f_{\mu,\lambda}^1(\tau)=f_{\mu,\lambda}(\tau)\mathds{1}_{|\tau|\leq 2\mu}$ and $f_{\mu,\lambda}^2(\tau) = f_{\mu,\lambda}(\tau)\mathds{1}_{|\tau|> 2 \mu}$.
Using that $m_{\mu}$ is uniformly bounded, we have
\bna
\nor{\langle \tau \rangle^\sigma f^1_{\mu,\lambda}}{L^2}^2&\leq& C\langle 2 \mu \rangle^\sigma \int_{-2\mu}^{2\mu}e^{-\frac{|\tau-\lambda|^2}{\mu}}|\widehat{\chi}(\tau-\lambda)|^2~d\tau
=C \langle 2 \mu \rangle^\sigma \int_{-2\mu-\lambda}^{2\mu-\lambda}e^{-\frac{|\tau|^2}{\mu}}|\widehat{\chi}(\tau)|^2~d\tau\\
&\leq &C\langle \lambda  \rangle^\sigma\int_{|\tau|\geq |\lambda|/2}e^{-\frac{|\tau|^2}{\mu}}|\widehat{\chi}(\tau)|^2~d\tau
\leq C\langle \lambda  \rangle^\sigma e^{-\frac{|\lambda|^2}{4\mu}}\int_{|\tau|\geq |\lambda|/2}|\widehat{\chi}(\tau)|^2\leq Ce^{-|\lambda|} ,
\ena
where we have used that $|\lambda|\geq 4\mu$ implies $|\lambda|/2\leq |\lambda|-2\mu$ (and in particular $\tau\in [-2\mu-\lambda,2\mu-\lambda]$ implies $|\tau|\geq |\lambda|/2$) and $\frac{|\lambda|^2}{4\mu}\geq |\lambda|$.

Concerning now $f_{\mu,\lambda}^2$, remark that $|s|\geq 2$ implies $\dist(s,[-1,1])\geq |s|/2$. Hence, using~\eqref{Estimate29} together with the support of $m$, we have uniformly
\bna
\left|m_{\mu}(s) \mathds{1}_{|s|\geq 2}\right|\leq C\left\langle \mu\right\rangle^{1/2}e^{-\frac{\mu s^2}{16}} \mathds{1}_{|s|\geq 2} .
\ena
Using this with $s=\tau/\mu$, we obtain
\bna
\nor{ \langle \tau \rangle^\sigma f^2_{\mu,\lambda}}{L^2}^2
\leq C\left\langle \mu\right\rangle^{1/2} \int_{|\tau|\geq 2\mu}\langle \tau \rangle^\sigma e^{-\frac{\tau^2}{16\mu }}e^{-\frac{|\tau-\lambda|^2}{\mu}}|\widehat{\chi}(\tau-\lambda)|^2~d\tau
\leq C \left\langle \mu\right\rangle^{1/2}  e^{-\frac{\lambda^2}{5\mu}} \int_\R \langle \tau \rangle^\sigma |\widehat{\chi}(\tau-\lambda)|^2~d\tau ,
\ena
where we have used the estimate $\frac{\tau^2}{16}+|\tau-\lambda|^2\geq \lambda^2 \min\left\{\left.s^2/16+(s-1)^2\right|s\in\R\right\}\geq c\lambda^2$ with $c>0$. 
Using now that $|\lambda|\geq 4\mu$, we have 
\bna
\nor{ \langle \tau \rangle^\sigma f^2_{\mu,\lambda}}{L^2}^2
 \leq C \left\langle \mu\right\rangle^{1/2}  e^{-c\frac{\lambda^2}{\mu}} \langle \lambda \rangle^\sigma
  \int_\R \langle \tau \rangle^\sigma |\widehat{\chi}(\tau)|^2~d\tau 
\leq C \langle \lambda \rangle^{\sigma +1/2}e^{-4c|\lambda|} ,
\ena
which concludes the proof of the lemma.
\enp

\section{The hypoelliptic heat equation}
\label{s:hypo-heat}
This section is devoted to the proofs of Theorems~\ref{t:approx-control-heat},~\ref{thm:para-gevrey} and~\ref{thm:parabolic}, which all rely on the methods of ~\cite[Propositions 1 and 2]{EZ:11} (proved in~\cite[Section~3]{EZ:11s}). We summarize these results in the next proposition for readibility.
 \begin{proposition}[\cite{EZ:11,EZ:11s}]
\label{propoEZ}
Let $T,S>0$ and $\alpha > 2 S^2$. Then, there exists some kernel function $k_T(t,s)$ such that
\begin{itemize}
\item
 if $\y$ is solution of the heat equation~\eqref{abstractheat}, then $\w(s)=\int_0^T k_T(t,s)\y(t)dt$ is solution of 
\bneqn
\label{wave-w}
\partial_s^2\w +\L \w&=&0, \quad \text{ for } s \in ]-S,S[ ,  \\
(\w,\partial_s\w)|_{s=0}&=&\left(0,\int_0^T \d_s k_T (t,0) \y(t)dt\right) = \left(0,\int_0^T e^{-\alpha \left(\frac{1}{t}+\frac{1}{T-t}\right)}\y(t)dt\right) ;
\eneqn
\item for all $\delta \in ]0,1[$, there is $C>0$ such that for all $(t, s) \in ]0,T[ \times ]-S,S[$, $k_T$ satisfies
\bnan
\label{e:estim-kT}
|k_T(t,s)|\leq  C|s|\exp\left( \frac{1}{\min \left\{ t,T-t\right\}}\left(\frac{s^2}{\delta}-\frac{\alpha}{(1+\delta)}\right)\right) .
\enan
\end{itemize}
\end{proposition}
Note that this last estimate is most useful for $\delta$ sufficiently close to one so that $\alpha \geq S^2(1+\frac{1}{\delta})$.

The proof of Theorems~\ref{thm:parabolic} and~\ref{thm:para-gevrey} then follows the Lebeau-Robbiano transmutation method, as implemented in~\cite{EZ:11}, splitting high and low frequencies. The proof of Theorems~\ref{t:approx-control-heat} is slightly different and does not rely on this splitting.
For the purposes of Theorems~\ref{thm:parabolic} and~\ref{thm:para-gevrey}, we define 
$$
E_{\lambda}=\vect \left\{\varphi_j, \lambda_j\leq \lambda\right\} ,
$$ 
where $(\lambda_j ,  \varphi_j)$ are the spectral elements of $\L$, defined in~\eqref{e:spectral-elts}. The first step of the proofs of Theorems~\ref{thm:parabolic} and~\ref{thm:para-gevrey} is to show, using the above transmutation technique, that we can transfer estimates obtained for solutions of the wave equation to solutions of the heat equation. More precisely, we first prove the following low-frequency observability estimate, with a precise estimation of the observability constant with respect to the cutoff frequency.
\begin{lemma}
\label{l:estim-heat-BF}
There exist $C, \gamma >0$ such that for any $T>0 , \lambda \geq0$, for every $\y_0 \in E_{\lambda}$ and associated solution $\y$ to~\eqref{abstractheat}, we have
\bnan
\label{estimLFheat}
\nor{\y(T)}{L^2}^2\leq \frac{C}{T}e^{\left( 2\gamma \lambda^{k/2} +\frac{C}{T}\right)}\int_0^T \int_{\omega}\left| \y(t,x)\right|^2~dt~dx .
\enan
Moreover, there exists $c_0>0$ such that for any $T>0$ there exists $C=C_{T} >0$ such that for any $\lambda \geq 0$, any $\y_0 \in E_{\lambda}$ and associated solution $\y$ to~\eqref{abstractheat}, we have
\bnan
\label{estimLFheat-y0}
\nor{\y_0}{L^2}^2\leq  C e^{2 c_0 \lambda^{k/2}}\int_0^T \int_{\omega}\left| \y(t,x)\right|^2~dt~dx .
\enan
\end{lemma}
\begin{remark}
\label{rem-ctes}
\begin{itemize} 
\item The constant $\gamma$ appearing in the exponent in \eqref{estimLFheat} may exactly be taken as $\gamma = \kappa +\eps$ for any $\eps>0$ where $\kappa$ appears in the exponent in Estimate~\eqref{th-estimate-k}, Theorem \ref{thmwavehypo} for some $S> \sup_{x\in \M} d_\L(x,\omega)$. In this case, the constant $C>0$ in front of the exponential also depend on $\eps$.
\item The constant $c_0$ appearing in the exponent in \eqref{estimLFheat-y0} may also be taken as $c_0 = \kappa +\eps$ for any $\eps>0$ (where $\kappa$ appears in the exponent in Estimate~\eqref{th-estimate-k}, Theorem \ref{thmwavehypo} for some $S> \sup_{x\in \M} d_\L(x,\omega)$) in the case where $k\geq 2$, but only $c_0 = \kappa + 2 \sqrt{\alpha}+\eps$ for any $\eps>0$ in the case $k=1$, which is the classical (elliptic) heat equation (where $\alpha$ is any constant $>\sqrt{2}\sup_{x\in \M} d_\L(x,\omega)$).
\item This is exactly the cost of controlling low frequencies, following~\cite{LR:95}.
For instance, \eqref{estimLFheat} implies that for all $\y_0\in \Pi_\lambda  L^2(\M) = E_\lambda$ ($\Pi_\lambda$ being the orthogonal projector associated to the spectral space of $\L$ with eigenvalues lower that $\lambda$), there exists $f \in L^2((0,T);L^2(\omega))$ with $\nor{f}{L^2((0,T);L^2(\omega))}^2\leq \frac{C}{T}e^{\left( 2\gamma \lambda^{k/2} +\frac{C}{T}\right)}\nor{\y_0}{L^2}^2$ such that the solution to 
\bneqn
\label{abstractheat-BF-Control}
\partial_t \y +\L \y&=&\Pi_\lambda \mathds{1}_\omega f\\
\y(0)&=&\y_0
\eneqn
satisfies $\y(T)=0$. Note that this finite dimensional observablity/controllabilty property is interesting in itself. For the time being and to the authors' knowledge, it is now understood in few situations, i.e. essentially in case $\L$ is an elliptic selfadjoint second order operator~\cite{LR:95}, the bi-Laplace operator~\cite{LRR:15}, the Stoke operator~\cite{CL:16}, and in case of some lower order perturbation of such operators~\cite{Lea:10}.

Again, the situation of Example~\ref{ex:Grushin++}, the exponent $\lambda^{k/2}$ is optimal in general, as can be seen when testing on eigenfunctions and using Proposition \ref{Prop:BCG}.
\end{itemize}
\end{remark}

In the proofs of Estimate~\eqref{estimLFheat-y0} and Theorem~\ref{t:approx-control-heat}, we shall moreover need the following definitions:
\bnan
\label{def:I-int}
\mathcal{I}(T,\lambda) = \int_0^T e^{-\alpha \left(\frac{1}{t}+\frac{1}{T-t}\right)} e^{-\lambda t }dt  ,
\enan 
and 
\bnan
\label{def:I}
\mathcal{I}(T,\L) u =  \sum_{j \in \N} \mathcal{I}(T,\lambda_j ) a_j \varphi_j = \sum_{j \in \N} \left( \int_0^T e^{-\alpha \left(\frac{1}{t}+\frac{1}{T-t}\right)} e^{-\lambda_j t }dt \right) a_j \varphi_j  , \quad \text{for } u = \sum_{j \in \N} a_j \varphi_j .
\enan

\bnp[Proof of Lemma~\ref{l:estim-heat-BF}]
The proofs of~\eqref{estimLFheat} and~\eqref{estimLFheat-y0} are similar. Let us start with that of~\eqref{estimLFheat}.
We start by using Theorem \ref{thmwavehypo-s} in the simpler case $s=k$. Now, we fix any $S> \sup_{x\in \M} d_\L(x,\omega)$: Estimate~\eqref{th-estimate-k} yields the existence of $C, \kappa, \mu_0>0$ such that  for all $\W_0 =(\w, \d_t\w)|_{s=0}$ (note that all constants then depend on these, and hence on the chosen $S>0$), the associated solution to~\eqref{wave-w} satisfies 
\bnan
\label{th-estimate-k-used}
\nor{\W_0}{L^2\times \H^{-1}_\L} \leq C e^{\kappa \mu}\nor{ \w}{L^2(]-S,S[\times \omega)} +\frac{1}{\mu}\nor{\W_0}{\H^{k}_\L\times \H^{k-1}_\L} , \quad \mu \geq \mu_0 .
\enan
Note that~\eqref{th-estimate-k-used} implies the same estimate for all $\mu>0$, in which case $\kappa$ has to be replaced by a bigger constant.

Assume now that $\w(s)$ is associated to $\y$ as $\w(s)=\int_0^T k_T(t,s)\y(t)dt$, where $\y$ is the solution to~\eqref{abstractheat} with initial datum $\y_0 \in E_{\lambda}$. Then, in~\eqref{wave-w}, $\W_0$ is of the particular form $\W_0 =  \left(0,\int_0^T e^{-\alpha \left(\frac{1}{t}+\frac{1}{T-t}\right)}\y (t)dt\right)$, so that a calculation (see~\cite[Equation~(3.3)]{EZ:11}) yields
\bnan
\label{e:brutale-est}
\nor{\W_0}{L^2\times \H^{-1}_\L}^2 & \geq & (1+ \lambda)^{-1}\nor{\W_0}{\H^1_\L \times L^2}^2 
= (1+ \lambda)^{-1}  \nor{\int_0^T e^{-\alpha \left(\frac{1}{t}+\frac{1}{T-t}\right)}\y (t)dt}{L^2}^2 \nonumber \\
& \geq & \frac{(1+ \lambda)^{-1}T^2}{9}e^{-\frac{9\alpha}{T}}\nor{\y(T)}{L^2}^2 .
\enan
Moreover, we have $\W_0 \in E_{\lambda} \times E_\lambda$ so that 
\bnan
\label{e:fct-freqW}
\frac{\nor{\W_0}{\H^{k}_\L \times \H^{k-1}_\L}}{\nor{\W_0}{L^2\times \H^{-1}_\L}} \leq (1+\lambda)^{\frac{k}{2}}.
\enan
As a consequence,~\eqref{th-estimate-k-used} implies
\bna
\left(1 - \frac{(1+\lambda)^{\frac{k}{2}}}{\mu} \right) \nor{\W_0}{L^2\times \H^{-1}_\L} \leq C e^{\kappa \mu}\nor{ \w}{L^2(]-S,S[\times \omega)}  , \quad \mu \geq \mu_0 ,
\ena
and hence, choosing $\mu = (1+\lambda)^{\frac{k}{2}}(1+\eps)$ for $\e \in (0,1)$, this is 
\bna
 \eps \nor{\W_0}{L^2\times \H^{-1}_\L} \leq C (1+\e)e^{\kappa  (1+\lambda)^{\frac{k}{2}}(1+\eps)}\nor{ \w}{L^2(]-S,S[\times \omega)}  , \quad \lambda \geq \lambda_0 = \mu_0^{\frac{2}{k}} ,
\ena
and  $\nor{\W_0}{L^2\times \H^{-1}_\L} \leq C_\eps e^{(\kappa +\eps)\lambda^{\frac{k}{2}}}\nor{ \w}{L^2(]-S,S[\times \omega)}$ for all $\eps >0$ (different from that in the previous line).
Using Cauchy-Schwarz inequality together with~\eqref{e:estim-kT} (with $\delta$ sufficiently close to one) we obtain that, for some $C>0$ depending only on $S,\alpha, \delta$, but not on $T$, we have
\begin{align}
\label{estimobserv}
\nor{\w}{L^2(]-S,S[\times \omega)}^2 &\leq \left( \int_{]0,T[\times ]-S,S[} k_T(t,s)^2 dt~ds \right) \int_{0}^T \int_{\omega}\left|\y (t,x)\right|^2 dx~dt \nonumber \\
& \leq C T e^{\frac{C}{T}}  \int_{0}^T \int_{\omega}\left|\y(t,x)\right|^2 dx~dt,
\end{align}
which then gives~\eqref{estimLFheat} for $\lambda \geq \lambda_0$. The estimate for $\lambda \in [0,\lambda_0]$ remains valid up to changing the constant $C$.

To prove~\eqref{estimLFheat-y0}, we follow the same lines, except for the lower bound~\eqref{e:brutale-est}, which we replace by an estimate of Corollary~\ref{cor:I} below. Namely, with the notation~\eqref{def:I}, we have $\W_0 = (0, \mathcal{I}(T,\L)\y_0)$,
so that, according to Corollary~\ref{cor:I}, we have for all $T>0$ and $s \in \R$, the existence of $C=C_{\alpha, T, s}>0$ such that 
\bnan
\label{e:refinedsqrt}
 \nor{\W_0}{L^2\times \H^{-1}_\L} = \nor{\mathcal{I}(T,\L) \y_0}{\H^{-1}_\L} \geq C \nor{\y_0}{\frac12,-2 \sqrt{\alpha}, -1-\frac32} 
 =C \nor{ (\L+1)^{-\frac{5}{4}} e^{-2 \sqrt{\alpha \L}} \y_0}{L^2} .
\enan
Recalling that $\y_0 \in E_\lambda$, this implies
\bna
 \nor{\W_0}{L^2\times \H^{-1}_\L} 
 \geq C \nor{ (\L+1)^{-\frac{5}{4}} e^{-2 \sqrt{\alpha \L}} \y_0}{L^2} \geq C \nor{ e^{-(2+\eps) \sqrt{\alpha \L}} \y_0}{L^2} \geq C e^{-(2+\eps) \sqrt{\alpha \lambda}} \nor{ \y_0}{L^2} .
 \ena
Applying then~\eqref{th-estimate-k-used} (for any $\mu>0$) with this lower bound, together with~\eqref{estimobserv} and~\eqref{e:fct-freqW} as above and the fact that $e^{(2+\eps) \sqrt{\alpha \lambda}} \leq C_\eps e^{\eps \lambda^{\frac{k}{2}}}$ for any $k\geq 2$ (in case $k=1$, the constant $(2+\eps) \sqrt{\alpha}$ has to be taken into account), concludes the proof of~\eqref{estimLFheat-y0}, and hence of the lemma.
\enp
Note that in the proof of~\eqref{estimLFheat-y0}, and in the case $k>2$, we could simply replace~\eqref{e:refinedsqrt} by the rough estimate
$$
\nor{\y_{\lambda}(0)}{L^2}^2\leq e^{2\lambda T}\nor{\y_{\lambda}(T)}{L^2}^2 ,
$$
which would be enough for the purpose of Estimate~\eqref{estimLFheat-y0}. This is not possible at all in case $k=1$, and in case $k=2$, would require $c_0$ to depend (linearly) on $T$.

\subsection{Approximate controllability with polynomial cost in large time: Proof of Theorem~\ref{thm:parabolic}}
\label{s:proof-thm:parabolic}
From the low-frequency Lemma~\ref{l:estim-heat-BF}, Estimate~\eqref{estimLFheat}, the proof of the theorem follows the spirit of~\cite{LR:95,Miller:10,EZ:11} but is simpler. It combines the cost of controllability of low frequencies, of order $e^{\gamma \lambda^{\frac{k}{2}}}=e^{\gamma \lambda}$ ($k=2$ in this part) and the dissipation of the heat at high frequency, of order $e^{-t\lambda}$. However, here, we do not perform the usual iterative procedure since it does not seem to improve the estimates.

\bnp[Proof of Theorem~\ref{thm:parabolic}]
For $\y \in L^2(\M)$, we decompose $\y=\y_{\lambda}+r_{\lambda}$ with $\y_\lambda \in E_\lambda$ and $r_\lambda \in E_\lambda^\perp$. 

On the one hand, using Lemma~\ref{l:estim-heat-BF} Estimate~\eqref{estimLFheat} for $\y_\lambda$ on the time interval $(T-\eta , T)$ (the problem being time invariant), we obtain, uniformly with respect to $T >0, \eta \in ]0,T[, \lambda >0 $, 
\bnan
\label{estimLFheat-bis}
\nor{\y_\lambda (T)}{L^2}^2\leq C e^{\left( 2\gamma \lambda^{k/2} +\frac{C}{\eta} \right)} \int_{(T-\eta)}^T \int_{\omega}\left| \y_\lambda(t,x)\right|^2~dt~dx .
\enan
On the other hand, we have 
\bnan
\nor{r_{\lambda}(t)}{L^2}\leq  e^{-\lambda t}\nor{r_{\lambda}(0)}{L^2}\leq e^{-\lambda t}\nor{\y(0)}{L^2}  ,  \label{e:rest-lam}\\
\int_{T-\eta}^T \int_{\omega}\left|r_{\lambda}(t,x)\right|^2~dt~dx\leq \frac{1}{2\lambda} e^{-2 \lambda (T-\eta)}\nor{\y(0)}{L^2}^2 \nonumber .
\enan
The last estimate gives
\bna
\int_{(T-\eta)}^T \int_{\omega}\left|\y_{\lambda}(t,x)\right|^2~dt~dx\leq 2\int_{(T-\eta)}^T \int_{\omega}\left|\y(t,x)\right|^2~dt~dx+\frac{1}{\lambda} e^{-2\lambda (T-\eta)}\nor{\y(0)}{L^2}^2.
\ena
So, using successively~\eqref{e:rest-lam}, \eqref{estimLFheat-bis} and the last estimate, we finally obtain for $T >0, \eta \in ]0,T[, \lambda >0$,
\bna
\nor{\y(T)}{L^2}^2&= & \nor{\y_{\lambda}(T)}{L^2}^2+ \nor{r_{\lambda}(T)}{L^2}^2\\
&\leq &\nor{\y_{\lambda}(T)}{L^2}^2+e^{-2\lambda T}\nor{\y(0)}{L^2}^2\\
&\leq & C e^{\left( 2\gamma \lambda^{k/2} +\frac{C}{\eta} \right)}  \int_{T-\eta}^T \int_{\omega}\left|\y(t,x)\right|^2~dt~dx
+C\left(e^{-2\lambda T}+ e^{\left( 2\gamma \lambda^{k/2} +\frac{C}{\eta} -2\lambda (T-\eta) \right)} \right)\nor{\y(0)}{L^2}^2 \\
&\leq & C e^{\left( 2\gamma \lambda^{k/2} +\frac{C}{\eta} \right)}  \int_{T-\eta}^T \int_{\omega}\left|\y(t,x)\right|^2~dt~dx
+2 C e^{\left( 2\gamma \lambda^{k/2} +\frac{C}{\eta} -2\lambda (T-\eta) \right)} \nor{\y(0)}{L^2}^2 .
\ena
Now, we recall that we assume $k=2$ (for $k>2$, the diffusion cannot compete with the cost of controllability of low frequencies).
Consequently, we obtain, for all $T >0, \eta \in ]0,T[, \lambda >0$,
\bnan
\label{e:heat-approx-lambda}
\nor{\y(T)}{L^2}^2 \leq  Ce^{C/\eta} \left( e^{2\gamma \lambda}  \int_{T-\eta}^T \int_{\omega}\left|\y(t,x)\right|^2~dt~dx
+ e^{2\lambda( \gamma + \eta - T)} \nor{\y(0)}{L^2}^2  \right).
\enan
This now provides information if $T$ is sufficiently large. Namely, setting $\eps = e^{-2\lambda}$, we obtain the existence of $C>0$ such that for all $\eta>0$,
$T\geq\gamma + \eta$, all  $\eps \in ]0, 1[$, we have 
\bna
\nor{\y(T)}{L^2}^2
\leq Ce^{C/\eta} \left( \frac{1}{\eps^{\gamma} } \int_{T-\eta}^T \int_{\omega}\left|\y(t,x)\right|^2~dt~dx+\eps^{T-(\gamma +\eta)} \nor{\y(0)}{L^2}^2 \right).
\ena
Changing $\eps^{T-(\gamma +\eta)}$ into $\eps$, this implies the existence of $C>0$ such that for all $\eta >0$, all
$T > \gamma + \eta$ and all  $\eps \in ]0, 1[$, we have 
\bna
\frac{e^{-C/\eta}}{C} \nor{\y(T)}{L^2}^2
\leq \frac{1}{\eps^{\frac{\gamma}{T-(\gamma +\eta)}}}
\int_{T-\eta}^T \int_{\omega}\left|\y(t,x)\right|^2~dt~dx+\eps \nor{\y(0)}{L^2}^2.
\ena
This concludes the proof of Theorem~\ref{thm:parabolic} with $T_0:=\gamma$ after having remarked that the parabolic dissipation yields $\nor{\y(T)}{L^2}^2\leq \nor{\y(0)}{L^2}^2$, and hence the case $\eps \geq 1$.
\enp

\subsection{Approximate controllability in Gevrey-type spaces: Proof of Theorem~\ref{thm:para-gevrey}}
\label{s:proof-thm:para-gevrey}
The proof of Theorem~\ref{thm:para-gevrey} follows the same lines as Theorem~\ref{thm:parabolic}, decomposing into low and high frequencies, but uses Estimate~\eqref{estimLFheat-y0} instead of~\eqref{estimLFheat}.

\bnp[Proof of Theorem~\ref{thm:para-gevrey}]
For $\y \in L^2(\M)$ arbitrary, we again write the decomposition $\y=\y_{\lambda}+r_{\lambda}$ with $\y_\lambda \in E_\lambda$ and $r_\lambda \in E_\lambda^\perp$. 
Note that, using the fact that $\y_{\lambda}$ is solution of the heat equation in $E_{\lambda}$, we obtain from Lemma~\ref{l:estim-heat-BF} that
\bna
\nor{\y_\lambda (0)}{L^2}^2\leq  C_T e^{2 c_0 \lambda^{k/2}}\int_0^T \int_{\omega}\left| \y_\lambda (t,x)\right|^2~dt~dx .
\ena
Moreover, we have (recall that the norm $\nor{\cdot}{\alpha, \theta}$ is defined in~\eqref{e:gevrey-norms}; $\alpha$ will be eventually taken equal to $k/2$)
\bna
\nor{r_{\lambda}(t)}{L^2}\leq e^{-\lambda t}\nor{r_{\lambda}(0)}{L^2}
\leq e^{-\lambda t} e^{-\theta \lambda^{\alpha }} \nor{r_{\lambda}(0)}{\alpha,\theta} \leq e^{-\lambda t -\theta \lambda^{\alpha }} \nor{\y(0)}{\alpha,\theta} \\
\int_{T-\eta}^T \int_{\omega}\left|r_{\lambda}(t,x)\right|^2~dt~dx\leq \frac{1}{2\lambda} e^{-2 \lambda (T-\eta)} e^{-2\theta \lambda^{\alpha }}\nor{\y(0)}{\alpha,\theta}^2 \nonumber. 
\ena
 From this last estimate, we obtain 
\bna
\int_{0}^T \int_{\omega}\left|\y_{\lambda}(t,x)\right|^2~dt~dx
& \leq & 2\int_{0}^T \int_{\omega}\left|\y(t,x)\right|^2~dt~dx+ 2\int_{0}^T \int_{\omega}\left|r_{\lambda}(t,x)\right|^2~dt~dx  \\
& \leq & 2\int_{0}^T \int_{\omega}\left|\y(t,x)\right|^2~dt~dx+  \lambda^{-1} e^{-2\theta \lambda^{\alpha }}\nor{\y(0)}{\alpha,\theta}^2.
\ena
So, combining all these estimates,  we finally obtain for $\lambda>0$
\bnan
\label{e:k<=2proof}
\nor{\y(0)}{L^2}^2& = &\nor{\y_{\lambda}(0)}{L^2}^2+ \nor{r_{\lambda}(0)}{L^2}^2 \nonumber \\
&\leq & \nor{\y_{\lambda}(0)}{L^2}^2+e^{-2\theta \lambda^{\alpha }} \nor{\y(0)}{\alpha,\theta}^2  \nonumber \\
& \leq &   C_T e^{2c_0 \lambda^{k/2}} \left( 2\int_{0}^T \int_{\omega}\left|\y(t,x)\right|^2~dt~dx
+   e^{-2\theta \lambda^{\alpha }}\nor{\y(0)}{\alpha,\theta}^2 \right) 
+e^{-2\theta \lambda^{\alpha }} \nor{\y(0)}{\alpha,\theta}^2  .
\enan
Now, for $\alpha=k/2$, we find for all $\lambda>0$ that
\bna
\nor{\y(0)}{L^2}^2\leq C e^{2c_0\lambda^{k/2}}\int_{0}^T \int_{\omega}\left|\y(t,x)\right|^2~dt~dx+Ce^{-2(\theta-c_0) \lambda^{k/2}}\nor{\y(0)}{k/2,\theta}^2.
\ena
Assuming $\theta > c_0$, and setting $\eps = Ce^{-2(\theta-c_0) \lambda^{k/2}} \in ]0,1[$, this is precisely~\eqref{estimobserheathypohigh} with $\theta_{0} = c_0$. The full range of $\eps>0$ follows from the simple estimate $\nor{\y(0)}{L^2}^2 \leq \nor{\y(0)}{k/2,\theta}^2$.
\enp

\subsection{Approximate controllability in natural spaces with exponential cost: Proof of Theorem~\ref{t:approx-control-heat}}
\label{s:proof-t:approx-control-heat}
Let us now proceed to the proof of Theorem~\ref{t:approx-control-heat}. It does not rely on frequency cutoff (we do not distinguish between low and high frequencies), and hence on Lemma~\ref{l:estim-heat-BF}. Instead, we directly apply the transmutation result of Proposition~\ref{propoEZ} to the full solution and use precise properties of the operator $\mathcal{I}(T,\L)$ defined in~\eqref{def:I} (which we aready used in the proof of Estimate~\eqref{estimLFheat-y0}), proved in the next section. Note also that here, as opposed to the above two sections, we need to use the strong version of Theorem~\ref{thmwavehypo}.

\bnp[Proof of Theorem~\ref{t:approx-control-heat}]
We apply directly the transmutation kernel to the solution. Using Theorem~\ref{thmwavehypo}, we obtain
\bnan
\label{th-estimate-k-H1}
\nor{\W_0}{L^2\times \H^{-1}_\L} \leq C e^{\kappa \mu^k}\nor{ \w}{L^2(]-S,S[\times \omega)} +\frac{1}{\mu}\nor{\W_0}{\H^{1}_\L\times L^2} , \quad \mu >0 ,
\enan
see e.g.~\cite[Lemma~A.3]{LL:15} to obtain the range $\mu \in [0,\mu_0]$ (however deteriorating the constant $\kappa$).
Then, we recall that, with the notation~\eqref{def:I}, we have 
$$
\W_0 = (0, \mathcal{I}(T,\L)\y_0) ,
$$
so that, according to Corollary~\ref{cor:I}, we have for all $T>0$ and $s \in \R$, the existence of $C=C_{\alpha, T, s}>0$ such that 
\bna
C^{-1} \nor{\y_0}{\frac12,-2 \sqrt{\alpha}, s-\frac32} \leq \nor{\W_0}{\H^{s+1}_\L\times \H^{s}_\L} = \nor{\mathcal{I}(T,\L) \y_0}{\H^{s}_\L} \leq C \nor{\y_0}{\frac12,-2 \sqrt{\alpha}, s-\frac32} , 
\ena
where $\nor{\y_0}{\frac12,-2 \sqrt{\alpha},s} = \nor{ (\L+1)^{\frac{s}{2}} e^{-2 \sqrt{\alpha \L}} \y_0}{L^2}$  (see~\eqref{def:anal-spaces} for the definition of the norms). In particular, this implies
\bna
\frac{\nor{\W_0}{\H^{1}_\L\times L^2}}{\nor{\W_0}{L^2\times \H^{-1}_\L}}
\leq C \frac{ \nor{\y_0}{\frac12,-2 \sqrt{\alpha},-\frac32}}{ \nor{\y_0}{\frac12,-2 \sqrt{\alpha}, -1 -\frac32}}
=  C\Lambda_1(H(\L) \y_0)
 \leq 2 C\Lambda_1(\y_0), \quad \Lambda_1(\y_0) = \frac{\nor{\y_0}{\H^1_\L}}{\nor{\y_0}{L^2}} , 
\ena
where $H(\lambda) = (\lambda+1)^{-\frac{5}{4}} e^{-2 \sqrt{\alpha \lambda}}$, which is positive, decreasing to zero, and we have thus used Corollary~\ref{cor:croissancefrequency} in the last inequality.

When combined with~\eqref{estimobserv} (still valid in this context), we now obtain, for all $\mu >0$,
\bna
  \nor{\W_0}{L^2\times \H^{-1}_\L}^2 \leq C e^{2\kappa \mu^k}  \int_{0}^T \int_{\omega}\left|\y(t,x)\right|^2 dx~dt+\frac{C}{\mu^2} \Lambda_1(\y_0)^2 \nor{\W_0}{L^2\times \H^{-1}_\L}^2.
\ena
Writing $\Lambda =\Lambda_1(\y_0)$, taking $\mu=\sqrt{2C}\Lambda$, and recalling that $C^{-1} \nor{\y_0}{\frac12,-2 \sqrt{\alpha}, -1-\frac32} 
 \leq \nor{\W_0}{L^2\times \H^{-1}_\L}$, this gives after absorption
\bnan
\label{almost-ultimate-heat}
 \nor{\y_0}{\frac12,-2 \sqrt{\alpha}, -1-\frac32}^2  \leq Ce^{c\Lambda^k}  \int_{0}^T \int_{\omega}\left|\y(t,x)\right|^2 dx~dt,\qquad \Lambda = \frac{\nor{\y_0}{\H^1_\L}}{\nor{\y_0}{L^2}} .
\enan
To conclude, we recall that $ \nor{\y_0}{\frac12,-2 \sqrt{\alpha}, -1-\frac32} = \nor{ (\L+1)^{-\frac{5}{4}} e^{-2 \sqrt{\alpha \L}} \y_0}{L^2} \geq C \nor{ e^{-3 \sqrt{\alpha (\L+1)}} \y_0}{L^2}$ and we use Lemma \ref{l:jensen-norm} with $F(s)=s+1$ and $G(s)= e^{-3 \sqrt{\alpha s}}$ which is convex, to finally obtain 
$$
\nor{\y_0}{L^2}\leq Ce^{c\Lambda_1(\y_0)}  \nor{ e^{-3 \sqrt{\alpha (\L+1)}} \y_0}{L^2} \leq Ce^{c\Lambda_1(\y_0)} \nor{\y_0}{\frac12,-2 \sqrt{\alpha}, -1-\frac32} .
$$ 
Together with~\eqref{almost-ultimate-heat}, this conludes the proof of~\eqref{e:approx-control-heat1}. Now, to prove \eqref{e:approx-control-heat2}, take any $\mu >0$. Either $\Lambda_1(\y_0) = \frac{\nor{\y_0}{\H^1_\L}}{\nor{\y_0}{L^2}} \geq \mu$, and \eqref{e:approx-control-heat2} holds (without the observation term), or else $\Lambda_1(\y_0) \leq \mu$, and \eqref{e:approx-control-heat1} yields \eqref{e:approx-control-heat2} (without additional term on the right hand-side). This concludes the proof of the Theorem.
\enp
\subsection{Technical lemmata used for the heat equation}
In this section, we collect three technical lemmata that we used in the proofs of Theorems~\ref{t:approx-control-heat},~\ref{thm:para-gevrey} and~\ref{thm:parabolic} above. 

First, we need an asymptotic expansion of the integral $\mathcal{I}(T,\lambda)$ defined in~\eqref{def:I-int} as $\lambda \to +\infty$.
\begin{lemma}
\label{lemme-laplace}
For all $\alpha>0$ and $T>0$, there exists $C_T, \lambda_0 >0$ such that for all $\lambda \geq \lambda_0$, there is $R(T, \lambda) \in \R$ such that we have 
\bna
\mathcal{I}(T,\lambda) =  \sqrt{\pi} \frac{\alpha^\frac14}{\lambda^\frac34} 
 e^{-\frac{\alpha}{T}} e^{-2\sqrt{\alpha \lambda}} \left( 1+ \frac{R(T,\lambda)}{\lambda^{\frac14}} \right), \qquad  | R(T,\lambda) | \leq C_T .
\ena
\end{lemma}
Next, this lemma allows us to link the operator $\mathcal{I}(T,\L)$ defined in~\ref{def:I}
and the norms 
\bnan
\label{def:anal-spaces}
\nor{u}{\delta,\theta,\sigma}^2
= \nor{ (\L+1)^{\frac{\sigma}{2}} e^{\theta \L^{\delta}}u}{L^2}^2
= \sum_{j \in \N} (\lambda_j+1)^\sigma e^{2\theta \lambda_j^{\delta}} |a_j|^2 
,  \quad \text{for } u = \sum_{j \in \N} a_j \varphi_j .
\enan
\begin{corollary}
\label{cor:I}
For all $s\in \R$, for all $T, \alpha>0$, there exists $C>1$ such that we have
\bna
C^{-1} \nor{u}{\frac12,-2 \sqrt{\alpha}, s-\frac32} \leq \nor{\mathcal{I}(T,\L) u}{\H^s_\L} \leq C \nor{u}{\frac12,-2 \sqrt{\alpha}, s-\frac32} , 
\ena
\end{corollary}
The proof of the corollary only consists in remarking that, once $\alpha ,T$ are fixed, we have, according to Lemma~\ref{lemme-laplace}, that
 $$ 
 0< \mathcal{I}(T,\lambda) \left(  \frac{e^{-2\sqrt{\alpha \lambda}} }{(1+\lambda)^\frac34}   \right)^{-1}  \to 
  \sqrt{\pi} \alpha^\frac14  e^{-\frac{\alpha}{T}}  >0 \quad \text{as } \lambda \to + \infty ,
 $$
 and this quantity does not vanish on $\R^+$ so that there is $C>1$ such that for all $\lambda \geq 0$, we have
  $$ 
C^{-1}\leq \mathcal{I}(T,\lambda) \left(  \frac{e^{-2\sqrt{\alpha \lambda}} }{(1+\lambda)^\frac34}   \right)^{-1}  \leq C ,
 $$
 which yields the result.
 
\bnp[Proof of Lemma~\ref{lemme-laplace}]
Note first that, given $\eps \in (0,1)$, we may assume that $\lambda_0$ is chosen such that $T \geq \sqrt{\frac{\alpha}{\lambda}}(1+\eps)$ for $\lambda \geq \lambda_0$. 
We first change variables in $\mathcal{I}(T,\lambda)$, denoting $\omega = \sqrt{\alpha \lambda}$ (new large parameter) and setting $t=\sqrt{\frac{\alpha}{\lambda}}s =\frac{\alpha}{\omega} s$, we have
\bna
\mathcal{I}(T,\lambda) = \frac{\alpha}{\omega}\int_0^{\frac{\omega T}{\alpha}} f_{\omega} (s)e^{-\omega \left( \frac{1}{s} + s\right)} ds, \qquad f_{\omega}(s) = \exp \left(- \frac{\alpha}{T - \frac{\alpha s}{\omega}} \right).
\ena
The phase $h(s) := \frac{1}{s} + s$ admits a single global strict (nondegenerate) minimum at the point $s=1$, with $h(1)=2$. Note also that $0\leq f_\omega \leq 1$. Hence, using that $\frac{\omega T}{\alpha} \geq 1+\eps$ by assumption, we have for  $\eps \in (0,1)$, the estimate
\bna
\mathcal{I}(T,\lambda) = \frac{\alpha}{\omega}\int_{1-\eps}^{1+\eps} f_{\omega} (s)e^{-\omega \left( \frac{1}{s} + s\right)} ds + O_T (e^{-(h(1) +c_\eps) \omega}) , \quad   c_\eps >0 .
\ena
Let $\phi : (1-\eps , 1+\eps) \to (-\eps_1 , \eps_2)$ be a (Morse) diffeomorphism for some $\eps_1, \eps_2>0$, such that $\phi(1)=0$, and with $u = \phi(s)$, we have
$$
h(s) = h(1) + h''(1) \frac{u^2}{2} = 2 +  u^2  , \quad \sgn(u) = \sgn(s-1) .
$$
Note that it is actually explicit, namely $\phi(s)=(s-1)/\sqrt{s}$. We change variable, setting $u = \phi(s)$, and obtain
\bna
\mathcal{I}(T,\lambda) = \frac{\alpha}{\omega}\int_{-\eps_1}^{\eps_2} e^{-\omega(2 +  u^2)}  f_{\omega} \circ \phi^{-1}(u) |(\phi^{-1})'(u)|du + O_T (e^{-(2 +c_\eps) \omega}) , \quad   c_\eps >0 ,
\ena
where $(\phi^{-1})'(0)=1$.

Moreover, for $s \in [1-\eps , 1+\eps] \subset [0, \frac{\omega T}{\alpha}]$, we write
$f_{\omega}(s)  =  f_{\omega}(1) + R_\omega(s)$ with  $f_{\omega}(1) = \exp \left(- \frac{\alpha}{T - \frac{\alpha}{\omega}} \right)$ and 
\bnan
\label{e:Romegataylor}
| R_\omega(s) | \leq |s-1|  \sup_{[1-\eps , 1+\eps] }|f_\omega'| \leq  \frac{c}{\omega}  |s-1| ,
\enan
where we used $f_\omega'(s) = - \frac{1}{\omega} X^{-2}e^{-\frac{1}{X}}$ with $X = (T-\frac{\alpha}{\omega}s) \alpha^{-1} \in [0,T/\alpha]$.

As a consequence, we obtain
\bna
\mathcal{I}(T,\lambda) &  = & \frac{\alpha}{\omega}e^{-2\omega} \exp \left(- \frac{\alpha}{T - \frac{\alpha}{\omega}} \right)  \int_{-\eps_1}^{\eps_2} e^{-\omega u^2}   |(\phi^{-1})'(u)|du + \mathcal{R}(\omega)  +   O_T (e^{-(2 +c_\eps) \omega}) \\ 
&  = & \frac{\alpha}{\omega}e^{-2\omega} \exp \left(- \frac{\alpha}{T - \frac{\alpha}{\omega}} \right) \left\{ \int_{-\eps_1}^{\eps_2} e^{-\omega u^2}   |(\phi^{-1})'(0)|du
+   \mathcal{R}'(\omega) \right\}+ \mathcal{R}(\omega)  + O_T (e^{-(2 +c_\eps) \omega}) \\ 
&  = & \frac{\alpha}{\omega}e^{-2\omega} \exp \left(- \frac{\alpha}{T - \frac{\alpha}{\omega}} \right)  \left\{ \int_\R e^{-\omega u^2}   |(\phi^{-1})'(0)|du
+   \mathcal{R}'(\omega) \right\}+ \mathcal{R}(\omega)  + O_T (e^{-(2 +\tilde{c}_\eps) \omega}) ,  \quad   \tilde{c}_\eps >0 \\
&  = & \frac{\alpha}{\omega}e^{-2\omega} \exp \left(- \frac{\alpha}{T - \frac{\alpha}{\omega}} \right)  \left\{  \sqrt{\frac{\pi}{\omega}}
+   \mathcal{R}'(\omega) \right\}+ \mathcal{R}(\omega)  + O_T (e^{-(2 +\tilde{c}_\eps) \omega}) ,  \quad   \tilde{c}_\eps >0 \\
\ena
with, using~\eqref{e:Romegataylor}, 
\bna
|\mathcal{R}(\omega)| = \left| \frac{\alpha}{\omega}\int_{-\eps_1}^{\eps_2} e^{-\omega(2 +  u^2)}  R_{\omega} \circ \phi^{-1}(u) |(\phi^{-1})'(u)|du  \right|
\leq C   \frac{\alpha}{\omega^2}  e^{-2\omega }\int_{-\eps_1}^{\eps_2} |u| e^{-\omega u^2} du \leq  \frac{\alpha}{\omega^2} e^{-2\omega } \frac{C}{\omega} ,
\ena
and 
\bna
|\mathcal{R}'(\omega) |= \left| \int_{-\eps_1}^{\eps_2} e^{-\omega u^2}  \left((\phi^{-1})'(u) - (\phi^{-1})'(0)\right)du \right| \leq C  \int_{-\eps_1}^{\eps_2} |u| e^{-\omega u^2} du \leq \frac{C}{\omega} .
\ena
Using finally that $\exp \left(- \frac{\alpha}{T - \frac{\alpha}{\omega}} \right) = e^{-\frac{\alpha}{T}}\left(1 + O_T\left( \frac{\alpha}{\omega}\right)\right)$, we finally obtain 
\bna
\mathcal{I}(T,\lambda) 
&  = & \frac{\alpha}{\omega}e^{-2\omega} \exp \left(- \frac{\alpha}{T - \frac{\alpha}{\omega}} \right) \sqrt{\frac{\pi}{\omega}}
\left( 1 + O_T\left(\frac{1}{\sqrt{\omega}}\right) \right)   ,
\ena
which, recalling that  $\omega = \sqrt{\alpha \lambda}$, concludes the proof of the lemma.
\enp

\begin{lemma}
\label{l:jensen-norm}
Let $F : \R^+ \to \R^+$ be any function and let $G : F(\R^+) \to \R^+$ be a function such that $G^2$ is convex (it is for instance the case if $G$ is). 
Then, for all $u \in \H^F_\L \cap D(G\circ F(\L))$, we have
\bnan
\label{e:jensen-norm}
G\left(\frac{\nor{u}{\H^F_\L}^2}{\nor{u}{L^2}^2} \right) \nor{u}{L^2}  \leq \nor{G \circ F(\L)u}{L^2} ,
\enan
where the seminorm $\H^F_\L$ is defined by $\nor{u}{\H^F_\L}^2 = \sum_{j\in \N} F(\lambda_j)  |a_j|^2$ if $u = \sum_{j\in \N} a_j \varphi_j$.
\end{lemma}
\begin{remark}
Using the previous lemma with $F(s)=s+1$ and $G(s)=\frac{1}{\sqrt{s}}$, we obtain the interpolation inequality
\bna
\frac{\nor{u}{L^2}}{\nor{u}{\H^{-1}_\L}}\leq \frac{\nor{u}{\H^1_\L}}{\nor{u}{L^2}} ,
\ena
comparing two types of ``frequency functions'', the first of which being used e.g. in~\cite{Phung:04} for the classical heat equation.
\end{remark}
\bnp
Dividing by $\nor{u}{L^2}$ (if non zero, otherwise the inequality is clear), it is enough to prove \eqref{e:jensen-norm} assuming $\nor{u}{L^2}=1$. If so, we write $u = \sum_{j\in \N} a_j \varphi_j$ with $\sum_{j\in \N} |a_j|^2 = 1$. Using the Jensen inequality with the convex function $G^2$, we have 
\bna
\nor{G \circ F(\L)u}{L^2}^2 = \sum_{j\in \N} G\big(F(\lambda_j) \big)^2 |a_j|^2  \geq G^2 \left( \sum_{j\in \N} F(\lambda_j ) |a_j|^2\right) = G^2(\nor{u}{\H^F_\L}^2),
\ena
which concludes the proof of the lemma.
\enp

\begin{lemma}
\label{lmcroissancefrequency}
Let $F, G : \R^+ \to \R^+$ be two nondecreasing continuous functions such that $F(s)G(s) \to +\infty$ as $s \to + \infty$. Then, for all $u \in D(F(\L)G(\L))$, we have
\bna
\nor{F(\L)u}{L^2}\nor{G(\L)u}{L^2}\leq 2 \nor{F(\L)G(\L)u}{L^2}\nor{u}{L^2}.
\ena
\end{lemma}
Note that, replacing $F$ and $G$ by $1/F$ and $1/G$, the same statement is true as well if $F,G$ are nonvanishing, nonincreasing and such that $F(s)G(s) \to 0$ as $s \to + \infty$.
\begin{corollary}
\label{cor:croissancefrequency}
Denoting, for $\sigma>0$ by $\Lambda_\sigma(u)=\frac{\nor{u}{\H^\sigma_\L}}{\nor{u}{L^2}}$ the modified frequency functions, we have that for any $H : \R^+ \to \R^+$ nonvanishing, nonincreasing such that $H(s) \to 0$ as $s \to + \infty$, 
\bnan
\label{inegfrequency}
\Lambda_\sigma(H(\L)u) \leq 2 \Lambda_\sigma(u) .
\enan
\end{corollary}
The corollary is obtained by taking $F(s)=(s+1)^{\frac{\sigma}{2}}$ in Lemma~\ref{lmcroissancefrequency}, $G= 1/H$ and $u=H(\L)v$.
Remark that the frequency function $\Lambda$ used in the main part of the article is $\Lambda = \Lambda_1$.
\begin{remark}
\label{rkinverse}
The interpretation of the corollary is clearer. Indeed, in this context, $H(\L)$ is a compact operator of $L^2$ and \eqref{inegfrequency} only translates that the ``average frequency'' of $H(\L)u$ is smaller than the ``average frequency'' of $u$.
\end{remark}
\begin{remark}

It is very likely that the previous Lemma (or at least Corollary~\ref{cor:croissancefrequency}) is still true with the constant $2$ replaced by $1$. Indeed, when taking $H(s)= (s+1)^{-\beta}$, or $H(s)=e^{-t  s^\beta}$ with $\beta, t >0$, Corollary~\ref{cor:croissancefrequency} is true with a constant $1$; in the first case, it is proved using Sobolev interpolation and in the second one using the monotonicity of the frequency function for solutions of the heat equation, see Phung \cite{Phung:04}.
\end{remark}

\bnp[Proof of Lemma~\ref{lmcroissancefrequency}]
First, since $F$ and $G$ are nondecreasing, by decomposing $u= \sum a_j \varphi_j$ with frequency less than $\lambda$ and greater than $\lambda$, we notice that for any $\lambda>0$, we have
\bna
\nor{F(\L)u}{L^2}^2 & =&   \sum_j F(\lambda_j)^2 |a_j|^2
 = \sum_{\lambda_j \leq \lambda} F(\lambda_j)^2 |a_j|^2  + \sum_{\lambda_j > \lambda} F(\lambda_j)^2 |a_j|^2 \\
 & \leq & F(\lambda)^2 \sum_{\lambda_j \leq \lambda} |a_j|^2  + \sum_{\lambda_j > \lambda} F(\lambda_j)^2 \frac{G(\lambda_j)^2}{G(\lambda)^2}|a_j|^2  \\
 & \leq & F(\lambda)^2 \nor{u}{L^2}^2 +\frac{1}{G(\lambda)^2} \nor{F(\L)G(\L)u}{L^2}^2 .
 \ena
Similarly, we have
 \bna
\nor{G(\L)u}{L^2}^2 \leq G(\lambda)^2 \nor{u}{L^2}^2 +\frac{1}{F(\lambda)^2} \nor{F(\L)G(\L)u}{L^2}^2.
\ena
Multiplying these two estimates, we obtain, for all $\lambda>0$
\bna
\nor{F(\L)u}{L^2}^2 \nor{G(\L)u}{L^2}^2 \leq  F(\lambda)^2 G(\lambda)^2\nor{u}{L^2}^4 + 2\nor{u}{L^2}^2 \nor{F(\L)G(\L)u}{\L^2}^2 +\frac{1}{F(\lambda)^2 G(\lambda)^2} \nor{F(\L)G(\L)u}{L^2}^4 ,
\ena
which is 
\bna
\nor{F(\L)u}{L^2} \nor{G(\L)u}{L^2} \leq  F(\lambda) G(\lambda) \nor{u}{L^2}^2 + \frac{1}{F(\lambda) G(\lambda)} \nor{F(\L)G(\L)u}{L^2}^2 ,  \quad \text{ for all } \lambda \geq 0 .
\ena
Now, using that $s \mapsto F(s) G(s)$ tends to $+ \infty$ as $s \to + \infty$ and is nondecreasing, it is onto $\R^+ \to [F(0)G(0),  + \infty[$. 
Since $F(0)G(0) = \min FG \leq \frac{\nor{F(\L)G(\L)u}{L^2}}{\nor{u}{L^2}}$, there is $\lambda \geq 0$ such that $F(\lambda) G(\lambda)=\frac{\nor{F(\L)G(\L)u}{L^2}}{\nor{u}{L^2}}$, which together with the last estimate yields the sought result.
\enp

\section{A partially analytic example: Grushin type operators}
\label{s:non-analytic}

In this section, we are concerned with the setting of Example~\ref{ex:Grushin+++} and give a proof of Theorem~\ref{t:partially-anal}. As explained in Section~\ref{s:sketch-plan}, it only suffices to prove the analogue of Theorem~\ref{thmwavehypo} (with estimate~\eqref{th-estimate-partial} instead of~\eqref{th-estimate-k}), that is for the hypoelliptic wave equation; all other results are then deduced as in Section \ref{s:hypo-heat}.

\bigskip
The setting of Example~\ref{ex:Grushin+++} differs from the general setting of the paper (compact manifolds, analytic context) with two respects: (i) we do not suppose analyticity in all variables; (ii) the manifold $\M = [-1,1] \times \T$ has a boundary. Hence, there are four main differences in the proofs, the first of which being of geometric nature, the next two being linked to the analysis of~\cite{LL:15}, and the last one to hypoelliptic estimates: 
\begin{enumerate}
\item \label{diff2} the presence of the boundary makes it complicated to apply directly Theorem \ref{th:riff-trel} coming from~\cite{RiffordTrelatMorse}.
\item \label{diff1} the partial analyticity assumption does not allow to make changes of variables to define the relation $\lhd$. For any couple of points $x^0$, $x^1$, we thus have to find some global set $\Omega$ containing one (short) path linking them.
\item \label{diff3} the application of the results in~\cite{LL:15} yields an observation term of the form $\nor{u}{H^{1}_b(\tilde{\omega})}$ and we would expect it to be in $L^2$. 
\item \label{diff4} The available hypoelliptic estimates, similar to those of Theorem \ref{thmhypoestim}, do not apply directly in the presence of boundary.
\end{enumerate}

The problem imposed by Item \ref{diff2} is that because of the boundary, the shortest path between two points in $\Int(\M)$ does not necessarily exists inside of $\Int(\M)$. To understand this issue, it may help to think about the flat metric in $\R^n \setminus O$ where $O$ a convex obstacle. The boundary $\d\M =\{x_1 =\pm 1\}$ can look like a convex obstacle for the Grushin case for instance (see Figure 3.1 of \cite{BoscLaurent:13} for some drawing of geodesics in Grushin). The solution we propose is to apply the result of Rifford and Tr\'elat~\cite{RiffordTrelatMorse} only locally away from the boundary. The drawback is then that our path is only piecewise normal geodesic. But this will be sufficient thanks to the variant Proposition \ref{propvariantbypiece} of Proposition \ref{propdepgeodes}.

The solution to the issue of Item \ref{diff1} is the very simple geometry of $[-1,1]_{x_1}\times \T_{x_2}$, so we almost do not perform any change of variable. 

The solution to the issue of Item~\ref{diff3} is to use the fact that the operator $P$ is elliptic in $\{\zeta_a=0\}$ where $\zeta_a$ is the dual to the anaytic variable $\z_a = (t,x_2)$, see Section \ref{subsectL2partial}.

Concerning the issue of Item~\ref{diff4}, we prove the necessary estimates in Section \ref{subsectsubellbound}. Recall that the operator is elliptic close to the boundary. So, we are left to patch the usual elliptic estimates close to the boundary with internal hypoelliptic estimates.

\bigskip
All in all, the proof of Estimate~\eqref{th-estimate-partial} is as above in two steps: first, proving~\eqref{e:intro-estim-partial-wave} (hard part), and then performing energy estimates (soft part). The latter are done the same way as in Section~\ref{s:simple-case}, except that the hypoelliptic estimates of Corollary \ref{corhypoestim} have to be replaced by those of Theorem \ref{cor:HsHsLbord} (with boundary).

We now focus on the first part of the proof, that is, proving~\eqref{e:intro-estim-partial-wave} in the context of Example~\ref{ex:Grushin+++}. This corresponds to the above Step~1 (Section~\ref{subsectstepgeometry}) and Step~2 (Section \ref{sectstepsmallness}).

\subsection{The geometric context}

Denote $\pi:[-1,1]_{x_1}\times \R_{x_2} \rightarrow [-1,1]_{x_1}\times \T_{x_2}$ the natural covering map, $\pi(x_1,x_2)=(x_1,x_2 +\Z)$. The vector fields $X_1$ and $X_2$ can be lifted to $[-1,1]_{x_1}\times \R_{x_2}$, which allows to define the natural sub-Riemannian distance on $]-1,1[_{x_1}\times \R_{x_2}$. As for the case of $[-1,1]_{x_1}\times \T_{x_2}$, the latter can be extended up to the boundary as well as all the other notions naturally inherited. We keep the same notations without leading to any confusion.

We will need the following Geometric Lemma, the proof of which relies on an iterative use of a slight variant of the result of Rifford-Tr\'elat~\cite{RiffordTrelatMorse}, see Theorem~\ref{th:riff-trel}.
\begin{lemma}
\label{lmconstructgeodespiece}
Let $x^0=(x_1^0,x_2^0)$ and $x^1=(x_1^1,x_2^1)$ in $[-1,1]_{x_1}\times \T_{x_2}$. Then, for any $\e>0$, there exists a continuous path  $\gamma: [0,1]\mapsto [-1,1]_{x_1}\times \R_{x_2}$ so that  with $\gamma(s)=(x_1(s),x_2(s))$ we have
\begin{enumerate}
\item $\pi(\gamma(0))=x^0$ and $\dist(\pi(\gamma(1)),x^1)<\e$;
\item $x_1(s)\notin \{-1,1\}$ for $s\in ]0,1[$;
\item $\gamma$ is piecewise normal geodesic in $]-1,1[_{x_1}\times \R_{x_2}$;
\item \label{itemborddroit}if $x_1^1=-1$ (resp. $1$) then there is $\delta>0$ so that $\gamma(s)=(-s,x_2^1)$ (resp. $\gamma(s)=(s,x_2^1)$) for $s\in [1-\delta,1]$. Similarly, if $x_1^0=-1$ (resp. $1$) then there is $\delta>0$ so that $\gamma(s)=(-1+s,x_2^0)$ (resp. $\gamma(s)=(1-s,x_2^0)$) for $s\in [0,\delta]$;
\item $\length(\gamma)\leq \dist(x^0,x^1)+\e$.
\end{enumerate}
\end{lemma}
\bnp
Note first that the paths defined in Item \ref{itemborddroit} are normal geodesic paths (i.e. geodesics, since the metric is Riemannian near the boundary) corresponding to $(\xi_1,\xi_2)=(\pm 1/2, 0)$ since $f$ does not depend on $x_2$ near the boundary. Therefore, by defining $\gamma$ like this for $s\in [1-\delta,1]$, we have $\gamma(1-\delta)\in ]-1,1[_{x_1}\times \T_{x_2}$ and $\length(\gamma(s), s \in [1-\delta,1])= \delta$; hence up to changing the length of $\gamma$ by $\delta$, we are left to the case where $x^1$ does not belong to the boundary. The argument shows that we may assume as well that $x^0$ does not belong to the boundary.

Let now $\widetilde{\gamma}$ be a smooth path on $[-1,1]_{x_1}\times \T_{x_2}$ so that $\widetilde{\gamma}(0)=x^0$, $\widetilde{\gamma}(1)=x^1$ and $\length(\widetilde{\gamma})\leq \dist(x^0,x^1)+\e$. We select one continuous lifting of $\widetilde{\gamma}$ on $[-1,1]_{x_1}\times \R_{x_2}$, denoted by $\gamma_1$, so that $\pi(\gamma_1(s))=\widetilde{\gamma}(s)$ for $s\in [0,1]$. Moreover, we have $\length(\gamma_1)=\length(\widetilde{\gamma})$ that we denote by $L$. Since $x^0, x^1 \notin \{\pm 1\} \times \T$, then up to deforming a bit $\gamma_1$ without changing $\gamma_1(0)$ (and still denoting it with the same name), we can assume that $\dist(\gamma_1(s),\{\pm1\}\times \R)>\eta>0$ for all $s \in [0,1]$, up to having only the estimate $\length(\gamma_1)\leq \dist(x^0,x^1)+2\e$. Now, we choose $N\in\N$ large enough so that $(L+\e)/N< \eta$. Up to reparametrization, we can also assume that $\gamma_1 :[0,L] \to ]-1,1[\times \R$ has unit speed. 

Denote $t_i=\gamma_1(iL/N)$ for $i=0,\cdots,N$. 
In particular, we have $\dist(t_i,t_{i+1})\leq L/N$ since $\gamma_1$ has unit speed.

\medskip

We now define $m_i$ for $i=0,\cdots,N$ by induction, so that the following properties are satisfied:
\begin{itemize}
\item[(P1)] $\dist(m_i,t_i)<\e/2N$;
\item[(P2)] there is a minimizing normal geodesic between $m_i$ and $m_{i+1}$.
\end{itemize}
 Note that these properties imply in particular 
 \bnan
\label{distOK}
\dist(m_i, \{\pm1\}\times \R)\geq\dist(t_i,\{\pm1\}\times \R)-\dist(m_i,t_i)>\eta -\e/2N>(L+\e/2)/N.
\enan
Let us now construct the points $m_i$ by induction as follows:
\begin{itemize}
\item $m_0=\widetilde{\gamma}(0)$;
\item $m_i\rightarrow m_{i+1}$: by iteration hypothesis, $\dist(m_i,t_i)\leq \e/N$.
Note that the ball $B(m_i,(L+\e/2)/N)$ does not intersect the boundary thanks to \eqref{distOK} and  that $\dist(m_i,t_{i+1})\leq \dist(m_i,t_{i})+\dist(t_i,t_{i+1})< (L+\e/2)/N$. In particular, $t_{i+1}\in B(m_i,(L+\e/2)/N)$, and there exists one ball of radius $r<\e/2N$ so that $B(t_{i+1},r)\subset B(m_i,(L+\e/2)/N)$.

A slight variant of Theorem~\ref{th:riff-trel} of Rifford-Tr\'elat~\cite{RiffordTrelatMorse} implies that the image of the exponential map (given by $T^*(\R^2)\to \R^2$, $(m_0, \xi_0) \mapsto m(1)$ where $(m(t), \xi(t))$ is the Hamiltonian curve issued from $(m_0, \xi_0)$, see Definition~\ref{def:Hamilton}) from the point $m_i$ is dense in $B(m_i,(L+\e/2)/N)$. In particular, there exists one point, which we choose as $m_{i+1}\in B(t_{i+1},r)\subset B(m_i,(L+\e/2)/N)$ so that there exists a minimizing normal geodesic between $m_i$ and $m_{i+1}$, and (P2) is satisfied. Then, we have by construction $\dist(m_{i+1},t_{i+1})<r<\e/2N$, so the first induction assumption (P1) is also fulfilled.
\end{itemize}

Once the process is finished, we have by construction, $\dist(m_i,m_{i+1})\leq \dist(m_i,t_i)+\dist(t_i,t_{i+1})+\dist(m_{i+1},t_{i+1})\leq (L+\e)/N$.
Now, denote by $\gamma$, the path defined by concatenation of the above defined normal geodesic path linking $m_i$ and $m_{i+1}$. 
We have $\gamma(0)=m_0=\widetilde{\gamma}(0)=\gamma_1(0)$ by construction. Also, we have $\dist(\gamma(L),\gamma_1(L))=\dist(m_{N-1},t_{N-1})<\e/2N<\e$.
Moreover, since the geodesic linking $m_i$ and $m_{i+1}$ are minimizing, we have $\length(\gamma)=\sum_{i=0}^{N-1} \dist(m_i,m_{i+1})\leq L+\e$. This concludes the proof of the lemma.
\enp

\subsection{A proof of Estimate~\eqref{th-estimate-partial}}

Let us now sketch the proof of Estimate~\eqref{th-estimate-partial}. It is very similar to that of Theorem~\ref{thmwavehypo}, so that we only stress the main differences.

For this, we define $\Omega = I \times U$ with $I$ a bounded neighborhood of $(-T,T)$ (where $T$ is that of the statement of Theorem~\ref{t:partially-anal}) and $U$ a bounded neighborhood of $\gamma$ in $[-1,1]_{x_1} \times \R_{x_2}$. We consider the operator $P=\d_t^2+\L$ in this set, and use the splitting of variables in $\R^n = \R^{n_a} \times \R^{n_b}= \R^{d+1}$ with $n=3,n_a=2,n_b=1,d=2,$ as
$$
\z = (\z_a, \z_b) , \quad \z_a = (t, x_2), \quad \z_b = x_1 ,
$$ 
with $t$ being the time variable, and $x=(x_1,x_2)$ the space variable.

Now, we follow the general proof.
The geometrical context being made precise in Lemma~\ref{lmconstructgeodespiece}, it only remains to check that we can apply the equivalent of Proposition~\ref{propvariantbypiece}, with the scheme of proof described in Remark \ref{rem:propvariantbypiece}. Since the appropriate piecewise geodesic path is constructed in Lemma \ref{lmconstructgeodespiece}, we only need to check that the local results (the equivalent of Lemma~\ref{lemhypolocal}) can be applied in this setting. 

There are two differences:
\begin{itemize}
\item We are in a situation where the only analytic variables are $z_a = (t, x_2)$. So, all Fourier multipliers defined in Section~\ref{s:def-mult} (and therefore the associated relation $\lhd$) are taken with respect to these variables, see Remark~\ref{rkinvariance}. The symbol of the wave operator $P$ is 
$$p(t,x_1,x_2,\tau,\xi_1,\xi_2)=-\tau^2+\xi_1^2+f(x_1,x_2)^2\xi_2^2 .$$
But we check that we are still in the situation of Remark 1.10 of~\cite{LL:15} with $\z_a=(t,x_2)$ and $\z_b=x_1$. Indeed, $p(t,x_1,x_2,0,\xi_1,0)=\xi_1^2$ that is positive definite on $\R_{\xi_1}$. 
\item The equivalent of Lemma \ref{lemhypolocal} should be obtained in the presence of boundary. We have to check that we can apply \cite[Theorem~5.12]{LL:15}, namely ``propagation up to the boundary'' $\left\{x_1= \pm1\right\}$. Let us only explain the construction near the boundary $\left\{x_1= 1\right\}$, the other case being similar.
One important thing is that we are in the geometric situation described in Lemma~\ref{lmconstructgeodespiece}: $P$ is already under the form of \eqref{Pnormal} and the choice of the geodesic close to the boundary of Item \ref{itemborddroit} of Lemma \ref{lmconstructgeodespiece} is already the same straight line. Indeed, we almost do not need to perform any change of coordinates, but only a translation. We can directly construct the noncharacteristic hypersurfaces of Lemma \ref{lmconstruction} with $l_0=\e_0$, $(t,\check{x})=(t,x_2)$ (tangential) and $x_d=x_1+1-\e_0$ (normal). Everything works then as in the interior case precised before, except for the last hypersurface $S_1= \{\phi_1=0\}$, which touches the the boundary $\left\{x_1=1\right\}$ tangentially. For this last step, we apply the local propagation result up to the boundary~\cite[Theorem~5.12]{LL:15}. We only need to check that the additional assumptions of this result are fulfilled:
\begin{itemize}
\item The analytic variable $\z_a = (t,x_2)$ are tangential with respect to the boundary $\left\{x_1=1\right\}$;
\item Assumption 5.1 in \cite{LL:15} is fulfilled for $P$ because close to the boundary, $p=q_x((\tau,\xi_2))+\widetilde{q}_x(\xi_1)$ where $q_x((\tau,\xi_2))=-\tau^2+f(x)^2\xi_2^2$ and $\widetilde{q}_x(\xi_1)=\xi_1^2$ are both quadratic forms independent on $t$ and $x_2$;
\item The boundary $\left\{x_1=1\right\}$ is non characteristic for $P$;
\item To apply Theorem 5.12 of \cite{LL:15} in this context,
calling $(x',x_n)$ the variables in that reference (the domain is locally $\{x_n>0\}$ in \cite{LL:15}, it is $\{x_1<1\}$ here), 
 one needs to set $x_n= 1-x_1$ and $x'=(t,x_2)$ so that $\left\{x_1 \leq 1\right\}$ is transformed into $\R^{n}_+$.
 The defining function of the last hypersurface $\phi_1(t,\check{x},x_d) = G((t,\check{x}),1) -x_d$ is changed, for the application of~\cite[Theorem~5.12]{LL:15} into $\tilde \phi_1(x',x_n) =G(x',1) -(1 - x_n)$. The assumption $\d_{x_n} \tilde \phi_1 =- \d_{x_d}\phi_1 =1>0$ of~\cite[Theorem~5.12]{LL:15} is hence satisfied.
\end{itemize}
\end{itemize}
This variant of Proposition~\ref{propvariantbypiece}, with an application of the boundary Theorem 5.12 of \cite{LL:15} as a last step, leads to the relation
\bna
\nor{M^{\beta\mu}_{\mu} \sigma_{r,\mu}  u}{1}\leq C e^{\kappa \mu}\left(\nor{M^{\alpha\mu}_{\mu} \vartheta_{\mu} u}{1} + \nor{Pu}{L^2(\Omega)}\right)+Ce^{-\kappa' \mu}\nor{u}{1}
\ena
where $M^{\beta\mu}_{\mu}$ is defined with the analytic variables $z_a = (t, x_2)$ and naturally extended to the boundary case since $z_a$ is tangential. The function $\vartheta \in C^{\infty}_0(\M)$ is chosen supported close to $x^0=(x_1^0,x_2^0)$ and can therefore be taken in $C^{\infty}_0(]-T,T[\times\omega)$. The function $\sigma_{r}\in C^{\infty}_0(\M)$ is equal to $1$ in a small ball centered at $x^1=(x_1^1,x_2^1)$ of size $r>0$.

Lemma \ref{lmobsL2} below allows to obtain with different constants
\bna
\nor{M^{\beta\mu}_{\mu} \sigma_{r,\mu}  u}{1}\leq C e^{\kappa \mu}\left(\nor{u}{L^2(]-T,T[\times \omega)} + \nor{Pu}{L^2(\Omega)}\right)+Ce^{-\kappa' \mu}\nor{u}{1} .
\ena
Again, as in~\cite[Section~4.2]{LL:15}, this leads, after a rough estimate of the high frequency, to
\bna
\nor{u}{L^2(]-\e,\e[\times B(x^1,r))}\leq C e^{\kappa \mu}\left(\nor{u}{L^2(]-T,T[\times \omega)} + \nor{Pu}{L^2(\Omega)}\right)+\frac{1}{\mu}\nor{u}{H^{1}(]-T,T[\times \M)} .
\ena
This is the equivalent to Proposition \ref{prop:local} which leads to a result similar to Corollary \ref{cor:H1th} by a compactness argument.

As explained above, the last step to get estimates \eqref{th-estimate-partial} corresponds to the energy estimates of Step 3 Section~\ref{s:simple-case}. There, they relied on the hypoelliptic estimates of Corollary \ref{cor:HsHsL}. The equivalent in the present situation with boundary is provided by Theorem \ref{cor:HsHsLbord}.

\subsection{An observation term in $L^2$ in quantitative unique continuation estimates}
\label{subsectL2partial}
In this section, we explain how the observation term $\nor{u}{H^{1}_b(\tilde{\omega})}=\sum_{|\beta| \leq 1}\nor{ D_b^\beta u}{L^2(\tilde\omega)}$ in unique continuation estimates as~\eqref{e:estimate-hyp} can actually be replaced by the weaker norm $\nor{u}{L^2(\tilde{\omega})}$ under suitable assumptions.
\begin{lemma}
\label{lmobsL2}
Let $\Omega$ be a bounded open set of $\R^{n}$ with $n = n_a+n_b$. Let $P$ be a differential operator of order $2$, defined in a neighborhood of $\Omega$, with real principal symbol and coefficients independent on the variable $\z_a$, and being elliptic in $\{\zeta_a=0\}$. Let $\omega\Subset \Omega$ and $\vartheta\in C^{\infty}_0(\omega)$. Then, for all $\alpha >0$, there exists $C>0$ such that for every $u\in C^{\infty}_0(\R^{n})$ and $\mu \geq 1$, we have 
\bna
\nor{M_{\alpha \mu}^{\mu}\vartheta_{\mu}u}{H^1}\leq C\left\langle \mu\right\rangle\nor{u}{L^2(\omega)}+C\nor{Pu}{L^2(\Omega)}+Ce^{-c\mu}\nor{u}{H^1}.
\ena
\end{lemma}
Recall that the regularization process $\vartheta \to \vartheta_{\mu}$ and the Fourier multiplier $M_{\alpha \mu}^{\mu}$ are defined at the beginning of Section~\ref{s:def-mult}.
\bnp
Since $P$ is elliptic (say positive to fix the ideas) in $\zeta_a=0$ and $\overline{\Omega}$ is compact, we can find $A>0$ (fixed for the rest of the proof) so that $A|\zeta_a|^2+p(\z_b , \zeta_a, \zeta_b)$ is elliptic on $\overline{\Omega} \times \R^{n_a+n_b}$ (where $p$ is the principal symbol of $P$), see for instance \cite[Lemma~A.1]{LL:15}. Using then the G\aa rding inequality, there exists $C>0$ so that
\bna
\nor{v}{H^1}^2&\leq& C \Re \left((A|D_a|^2+P)v,v\right)_{L^2}+C\nor{v}{L^2}^2\\
&\leq & C\nor{|D_a|v}{L^2}^2+C \Re\left( Pv,v\right)_{L^2}+C\nor{v}{L^2}^2
\ena
for every $v\in C^{\infty}_0(\Omega)$.
Let $\varphi, \chi\in C^{\infty}_0(\Omega)$ being real valued and such that $\varphi=1$ on a neighborhood of $\supp(\vartheta)$ and $\chi=1$  on a neighborhood of $\supp\varphi$. Applying this estimate to $v= \varphi w$ for $w\in C^{\infty}_0(\R^{n})$, we obtain that 
\bna
\nor{\varphi w}{H^1}^2&\leq & C \nor{\nabla_a (\varphi w)}{L^2}^2+C \Re\left( P(\varphi w),\varphi w\right)_{L^2}+C\nor{\varphi w}{L^2}^2\\
&\leq & C \nor{\nabla_a w}{L^2}^2+C \Re\left( P(\varphi w),(\varphi w)\right)_{L^2}+C\nor{w}{L^2}^2
\leq C \nor{\nabla_a w}{L^2}^2+C \nor{P(\varphi  w)}{H^{-1}}\nor{\varphi w}{H^1}+C\nor{w}{L^2}^2
\ena
Writing $P(\varphi w)=\varphi Pw+[P,\varphi] w$, where $[P,\varphi]$ is of order $1$, we have  
\bna
\nor{P(\varphi  w)}{H^{-1}} \leq \nor{\varphi P w}{H^{-1}}+ C \nor{w}{L^2} ,
\ena
so that, after absorption, we have proved the existence of $C>0$ such that for every $w\in C^{\infty}_0(\R^{n})$ (and so for $w\in\mathcal{S}(\R^{n})$) , we have
\bnan
\label{e:estim-w-L2}
\nor{\varphi w}{H^1}^2 \leq  C \nor{\nabla_a w}{L^2}^2+C \nor{\varphi P w}{H^{-1}}^2+C\nor{w}{L^2}^2.
\enan
We apply the previous estimate to $w=M_{\alpha\mu}^{\mu}\vartheta_{\mu}u$. For the first term in the right handside of~\eqref{e:estim-w-L2}, we have $\nor{\zeta_a m_\mu(\zeta_a/\mu)}{L^\infty(\R^{n_a})} \leq \mu \nor{\zeta_a m_\mu(\zeta_a)}{L^\infty(\R^{n_a})}\leq C\mu$, so that 
\bna
\nor{\nabla_a (M_{\alpha\mu}^{\mu}\vartheta_{\mu}u)}{L^2}\leq  C\left\langle \mu \right\rangle \nor{\vartheta_{\mu}u}{L^2}.
\ena
For this term, we further use that $\vartheta\in C^{\infty}_0(\omega)$, which according to Lemma \ref{Lemma23} Item~2, gives
\bna
\nor{\vartheta_{\mu}u}{L^2}\leq \nor{u}{L^2(\omega)}+ Ce^{-c\mu} \nor{u}{L^2}. 
\ena
Note that this previous inequality also rules the term $\nor{w}{L^2}$ in the right handside of~\eqref{e:estim-w-L2}. 

It only remains to estimate the term $\nor{\varphi P w}{H^{-1}}$ in the right handside of~\eqref{e:estim-w-L2}. Using that $P$ is invariant on $\z_a$ (and hence commutes with $M_{\alpha\mu}^{\mu}$), it is
\bna
\nor{\varphi P M_{\alpha\mu}^{\mu}\vartheta_{\mu}u}{H^{-1}}&\leq& \nor{ [P,\vartheta_{\mu}]u}{H^{-1}}+\nor{\varphi  M_{\alpha\mu}^{\mu}\vartheta_{\mu}Pu}{H^{-1}}\\
&\leq &\nor{ [P,\vartheta_{\mu}]u}{H^{-1}}+\nor{\varphi  M_{\alpha\mu}^{\mu}\vartheta_{\mu}\chi Pu}{H^{-1}}+\nor{\varphi  M_{\alpha\mu}^{\mu}(1-\chi) \vartheta_{\mu}Pu}{H^{-1}}.
\ena
Note that, {\em a priori}, since $P$ is not defined on the whole $\R^{n}$ but only in a neighborhood of $\Omega$, the term $Pu$ does not have any meaning. Yet, since $P$ is invariant in $\z_a$, the differential operator $\vartheta_{\mu} P$ is a well defined operator on $\R^{n}$, so as $\chi P$ and $[P,\vartheta_{\mu}]$. In the end, all terms involved are well defined for all $u \in C^\infty_0(\R^n)$, even if not supported inside of $\Omega$.

Now, using $\vartheta\in C^{\infty}_0(\omega)$, we have by (a dual version of) Lemma \ref{Lemma23} that
\bna
\nor{[P,\vartheta_{\mu}]u}{H^{-1}}\leq \nor{u}{L^2(\omega)}+ Ce^{-c\mu} \nor{u}{L^2} .
\ena
Finally, since $\supp (\varphi)\cap \supp (1-\chi)=\emptyset$, \cite[Lemma~2.10]{LL:15} also yields
\bna
\nor{\varphi  M_{\alpha\mu}^{\mu}(1-\chi) }{H^{-1}\rightarrow H^{-1}}\leq Ce^{-c\mu}.
\ena
Combining the last five inequalities together with~\eqref{e:estim-w-L2}, we are led, after absorption, to the estimate
\bna
\nor{\varphi M_{\alpha\mu}^{\mu}\vartheta_{\mu}u}{H^1}\leq C\nor{Pu}{L^2(\Omega)}+C\left\langle \mu \right\rangle \nor{u}{L^2(\omega)}+ Ce^{-c\mu} \nor{u}{H^1}.
\ena
The property $\supp (1-\varphi)\cap \supp (\vartheta)=\emptyset$ with \cite[Lemma~2.10]{LL:15} gives the similar estimate
\bna
\nor{(1-\varphi) M_{\alpha\mu}^{\mu}\vartheta_{\mu}u}{H^1}\leq Ce^{-c\mu} \nor{u}{H^1},
\ena
which allows to conclude the proof of the lemma.
\enp

\appendix

\section{On the optimality: Proof of Proposition~\ref{Prop:BCG}}
\label{s:proof-Prop:BCG}
In this appendix, we discuss the optimality of the results presented in the main part of the paper, in the situation of Example~\ref{ex:Grushin++}, i.e. we give a proof of Proposition~\ref{Prop:BCG}. The estimates we use are mainly extracted from the article \cite{BeauchardCanGugl} by Beauchard, Cannarsa and Guglielmi. They are slightly spread out in this reference so that the proof below mainly explains where in \cite{BeauchardCanGugl} to pick the results. This is mainly the proof of Theorem 5, Section 3.2 and 3.3 in this reference. 

\bnp[Proof of Proposition~\ref{Prop:BCG}]
First, Fourier transforming the operator $\L_\gamma = -(\d_{x_1}^2-x_1^{2\gamma}\d_{x_2}^2)$ in the $x_2$ variable, we obtain a family $A_{n,\gamma}$ of $1$-dimensional operators defined for $n \in \Z$ by $(A_{n,\gamma}f)(x_1)= -f''(x_1) + (n\pi)^{2} x_1^{2\gamma}f(x_1)$ on $]-1,1[$, with Dirichlet boundary conditions.
 
The sequence of eigenfunctions $\varphi_n$ is then taken of the form $\varphi_n(x_1,x_2)=\sqrt{2}v_n(x_1)\sin(n\pi x_2)$ where $v_n$ is the first normalized eigenvector of $A_{n,\gamma}$ (see Lemma 2 of \cite{BeauchardCanGugl}). We have $A_{n,\gamma}v_n= \lambda_{n,\gamma}v_n$, with $\lambda_{n,\gamma}$ the lowest eigenfunction of $A_{n,\gamma}$, and hence $\L_{\gamma}\varphi_n= \lambda_{n,\gamma}\varphi_n$.

The following estimates hold:
\begin{enumerate}
\item \label{item1BCG} $\frac{1}{C} n^{\frac{2}{1+\gamma}}\leq \lambda_{n,\gamma}\leq C n^{\frac{2}{1+\gamma}}$, see Proposition 4 of \cite{BeauchardCanGugl}.
\item for $0<a<b<1$, we have $\int_a^b v_n(x)^2dx\leq Ce^{-C_1(\gamma) n a^{\gamma+1} }c_n$ with $c_n\leq n^{\beta}$ for some appropriate $\beta$: this is inequality (35) of \cite{BeauchardCanGugl} once we have checked that for $n$ large enough $\mu_n=C(\gamma)n$ (written in (33)) and the definition of $x_n$ in (26). For $\gamma=1$, a more precise result is stated \cite[Lemma~4]{BeauchardCanGugl}, where the constant is computed, namely $\lambda_{n,\gamma}\approx n\pi$ and $\int_a^b v_n(x)^2dx\approx \frac{e^{-a^2 n\pi}}{2a\pi\sqrt{n}}$ for $a>0$.
\end{enumerate}
So, in any cases, if $a>0$, there are $C,c>0$, so that $\int_a^b v_n(x)^2dx\leq C e^{-c n}\leq C e^{-c \lambda_{n,\gamma}^{\frac{1+\gamma}{2}}}$ where we have used Item \ref{item1BCG}. Then, since $\overline{\omega} \cap \left\{x_1=0\right\} = \emptyset$, there exists $a>0$ so that $\nor{\varphi_n}{L^2(\omega)}^2\leq C\int_a^b v_n(x)^2dx$. To finish the proof, it is enough to notice that $v_n$ was chosen normalized in $L^2(]-1,1[)$ so that $\varphi_n$ is normalized in $L^2(\M)$.
\enp
\begin{remark}
It is very interesting to compare the estimates obtained by \cite{BeauchardCanGugl} with respect to those obtained in the present paper, even if the techniques are quite different.
The scheme of proof we followed may be summarized in two steps:
\begin{enumerate}
\item We prove an observability estimate where the cost is, more or less, exponential of the usual Sobolev frequency. This step is performed in \cite{BeauchardCanGugl} by using only the analyticity of the coefficients in the $x_2$ variable. In their Proposition 5, they prove a $1D$ Carleman estimate and the cost is of the order of $e^{cn}$ where $n$ is the frequency in the $x_2$ (the analytic frequency). 
\item Then, we use the decay of the heat flow using hypoelliptic estimates. For this, $\frac{1}{C} n^{\frac{2}{1+\gamma}}\leq \lambda_{n,\gamma}\leq C n^{\frac{2}{1+\gamma}}$ may be seen as a counterpart of the hypoelliptic estimates of Theorem \ref{thmhypoestim}. Indeed, these estimates roughly say that in the worst case, the operator $\L$ counts (when we want estimates from below) as $\frac{2}{1+\gamma}=\frac{2}{k}$ when compared to the usual derivatives (that is to the usual Sobolev norms).
\end{enumerate}

\end{remark}

\section{Subelliptic estimates}
\label{app:sectsub}
\subsection{$H^s$ subelliptic estimates on compact manifolds}
\label{app:Hs-comm}

In this appendix, we draw classical consequences of the subelliptic estimate~\eqref{estimhypo} of Rothschild-Stein \cite{RS:76} and Fefferman-Phong \cite{FP:83}, that are used in the main part of the paper.
The following corollary of the subelliptic estimate~\eqref{estimhypo} might be written elsewhere, but we did not find any reference. The short proof below stresses that the sole subelliptic estimate we rely on in the paper is~\eqref{estimhypo}.

\begin{corollary}
\label{corhypoestim}
Under Assumption \ref{assumLiek}, for any $s\geq 0$ there is $C>0$ such that we have
\bnan
\label{estimhypocor}
\nor{u}{H^{s+\frac{1}{k}}(\M)}^2\leq C\sum_{i=1}^m\nor{X_i u}{H^s(\M)}^2+C\nor{u}{L^2(\M)}^2 , \\
\label{estimhypocor-2}
\nor{u}{H^{s+\frac{2}{k}}(\M)}^2\leq C\nor{\L u}{H^s(\M)}^2+C\nor{u}{L^2(\M)}^2 ,
\enan
for any $u\in C^{\infty}(\M)$. 
\end{corollary}
The proof we give is inspired by~\cite{FP:83} (see the beginning of the proof of Theorem~1). For this, we let $\Lambda$ be an elliptic inversible pseudodifferential operator of order one in $\M$, being selfadjoint in $L^2(\M)$ (see e.g.~\cite[Remark~2.11]{LL:16} after having endowed $\M$ with a Riemannian metric). Recall that the power operator $\Lambda^s$ is an elliptic inversible pseudodifferential operator of order $s$ in $\M$, being also selfadjoint in $L^2(\M)$. All $H^s$ norms are equivalent to $\nor{\cdot }{H^s(\M)} = \nor{ \Lambda^s \cdot }{L^2(\M)}$.

\bnp
We start proving Estimate~\eqref{estimhypocor}, which is simpler due to the fact that $X_i$ is only of order $1$ and therefore $[X_i , \Lambda^s]$ is of order $s$ and hence an admissible remainder term (compared to the estimated norm $H^{s+\frac{1}{k}}$). 
Using the $L^2$ estimate \eqref{estimhypo}, we have
\bna
\nor{u}{H^{s+\frac{1}{k}}(\M)}^2&\leq& C \nor{\Lambda^s u}{H^{\frac{1}{k}}(\M)}^2\leq C\sum_{i=1}^m\nor{X_i \Lambda^s u}{L^2(\M)}^2+C\nor{\Lambda^s u}{L^2(\M)}^2\\
&\leq & C\sum_{i=1}^m\nor{\Lambda^s X_i u}{L^2(\M)}^2+C\sum_{i=1}^m\nor{[\Lambda^s, X_i] u}{L^2(\M)}^2+C\nor{\Lambda^s u}{L^2(\M)}^2\\
&\leq&  C\sum_{i=1}^m\nor{ X_i u}{H^s(\M)}^2+C\nor{u}{H^s(\M)}^2.
\ena
An interpolation estimate gives $\nor{u}{H^s(\M)}\leq \e \nor{u}{H^{s+\frac{1}{k}}(\M)}+C_{\e}\nor{u}{L^2(\M)}$ for any $\e>0$, which yields~\eqref{estimhypocor} after having taken $\eps$ sufficiently small.

\medskip
Concerning Estimate~\eqref{estimhypocor-2}, we have to be more careful since the commutator $[\L , \Lambda^s]$ is of order $s+1$ and hence not an admissible remainder term (compared to the estimated norm $H^{s+\frac{2}{k}}$). 
Following~\cite{FP:83}, we apply the $L^2$-based estimate~\eqref{estimhypo2} to $\Lambda^{s+\frac{1}{k}}u$, yielding
\bnan
\label{computcommut}
\nor{u}{H^{s+\frac{2}{k}}(\M)}^2&\leq& C \Re(\L\Lambda^{s+\frac{1}{k}}u,\Lambda^{s+\frac{1}{k}}u)_{L^2(\M)}+C\nor{\Lambda^{s+\frac{1}{k}}u}{L^2(\M)}^2\nonumber\\
&\leq &C \Re(\Lambda^{s+\frac{1}{k}}\L u,\Lambda^{s+\frac{1}{k}}u)_{L^2(\M)}+ C \Re([\L,\Lambda^{s+\frac{1}{k}}]u,\Lambda^{s+\frac{1}{k}}u)_{L^2(\M)} +C\nor{u}{H^{s+\frac{1}{k}}(\M)}^2\nonumber\\
&\leq &C \Re(\Lambda^{s}\L u,\Lambda^{s+\frac{2}{k}}u)_{L^2(\M)}+C\Re ([\L,\Lambda^{s+\frac{1}{k}}]u,\Lambda^{s+\frac{1}{k}}u)_{L^2(\M)} +C\nor{u}{H^{s+\frac{1}{k}}(\M)}^2\nonumber\\
&\leq &\frac{1}{2}\nor{u}{H^{s+\frac{2}{k}}(\M)}^2+C \nor{\L u}{H^{s}(\M)}^2+ C \Re(\Lambda^{s+\frac{1}{k}}[\L,\Lambda^{s+\frac{1}{k}}]u,u)_{L^2(\M)} +C\nor{u}{H^{s+\frac{1}{k}}(\M)}^2 ,
\enan
where we have used Cauchy-Schwarz inequality in the last step.
The term with the commutator has to be taken carefully since it is a priori of order $2s+\frac{2}{k}+1$. But the following simple remark is in order: this pseudodifferential operator has purely imaginary principal symbol. Hence, according to pseudodifferential calculus, it can be written as $\Lambda^{s+\frac{1}{k}}[\L,\Lambda^{s+\frac{1}{k}}]=T_1+T_2$ where $T_1$ is a {\em skew-adjoint} pseudodifferential operator of order ${2s+\frac{2}{k}+1}$, and $T_2$ is a pseudodifferential operator of order ${2s+\frac{2}{k}}$. In particular, we have
\bna
\left|\Re(\Lambda^{s+\frac{1}{k}}[\L,\Lambda^{s+\frac{1}{k}}]u,u)_{L^2(\M)}\right|=\left|\Re(T_2 u,u)_{L^2(\M)}\right|\leq \nor{u}{H^{s+\frac{1}{k}}(\M)}^2.
\ena
So, at this stage, we have proved 
\bna
\nor{u}{H^{s+\frac{2}{k}}(\M)}^2&\leq &C \nor{\L u}{H^{s}(\M)}^2 +C\nor{u}{H^{s+\frac{1}{k}}(\M)}^2.
\ena
This concludes the proof of~\eqref{estimhypocor-2}, after an interpolation argument as above.
\enp

To conclude this section, we prove the continuous injection $\H_\L^{s} \subset H^{s/k} (\M)$ for all $s\geq 0$.
We shall use the following classical operator theoretic (interpolation) result for which we refer e.g. to~\cite[Corollary~12.15]{SRV:10}.
Given two selfadjoint nonnegative operators $(A, D(A))$ and $(B, D(B))$ on a Hilbert space, we have
\bnan
\label{e:operator-momotone-functions}
\|Au\| \leq \|Bu\| \text{ for all }u \in D(B) \Longrightarrow \|A^\alpha u\| \leq \|B^\alpha u\| \text{ for all } \alpha \in [0,1]  \text{ and }  u \in D(B^\alpha) .
\enan
Note that this result already yields the simple inequality: for $s\geq 0$, there is $C>0$ such that for all $u \in H^s(\M)$, 
\bna
\label{e:sobolev-simple}
\| u \|_{\H_\L^{s}}\leq C\nor{u}{H^{s}(\M)},
\ena
consequence of that obtained for $s \in 2\N$. It also yields by duality for all $s\geq 0$ the existence of $C>0$ such that for all $u \in \H_\L^{-s}$, 
\bnan
\label{e:sobolev-simple-dual}
\| u \|_{H^{-s}(\M)}\leq C\nor{u}{\H_\L^{-s}} .
\enan

\begin{corollary}
\label{cor:HsHsL}
For all $s \geq 0$, there exists $C>0$ such that for all $u \in \H_\L^{s}$, we have
\bna
\| u \|_{H^{\frac{s}{k}} (\M)}\leq C\nor{u}{\H_\L^{s}} .
\ena
\end{corollary}
Note that it also yields by duality for all $s\geq 0$ the existence of $C>0$ such that for all $u \in H^{-\frac{s}{k}}(\M)$, 
\bna
\| u \|_{\H_\L^{-s}}\leq C\nor{u}{H^{-\frac{s}{k}}(\M)}.
\ena

\bnp

We first prove the result for $s=2 p$, $p \in \N$, and then conclude by interpolation. We prove by induction that for all $p \in \N$, we have
\bnan
\label{e:induc-Hs}
\| u \|_{H^{\frac{2p}{k}} (\M)}\leq C\nor{u}{\H_\L^{2p}} = C\nor{(\L + 1)^p u}{L^2(\M)}, \quad\text{ for all } u \in C^\infty(\M) .
\enan
The case $p=0$ is clear. Assume now that this is satisfied for $p$, and estimate $\| u \|_{H^{\frac{2(p+1)}{k}} (\M)} = \| u \|_{H^{\frac{2p}{k}+\frac{2}{k}} (\M)}$.
After having used~\eqref{estimhypocor-2}, we have
\bna
\| u \|_{H^{\frac{2(p+1)}{k}} (\M)}^2  
\leq C\nor{ \L u}{H^{\frac{2p}{k}}(\M)}^2  + C\nor{u}{L^2(\M)}^2 , 
\ena
which, using the induction assumption~\eqref{e:induc-Hs} to $\L u$, yields
\bna
\| u \|_{H^{\frac{2(p+1)}{k}} (\M)}^2  
\leq C\nor{(\L + 1)^p \L u}{L^2(\M)}^2  + C\nor{u}{L^2(\M)}^2 .
\ena
Using the functional calculus~\eqref{e:fct-calcul}, this implies
\bna
\| u \|_{H^{\frac{2(p+1)}{k}} (\M)}^2  
\leq C\nor{(\L + 1)^{p+1} u}{L^2(\M)}^2 = C\nor{u}{\H_\L^{2(p+1)}}^2 , \quad\text{ for all } u \in C^\infty(\M) ,
\ena
which is the sought estimate.

Now for $s\geq 0$, $s\notin \N$, pick $p\in \N$ such that $s \in [0,p]$, write~\eqref{e:induc-Hs} as $\| \Lambda^{\frac{2p}{k}} u \|_{L^2(\M)}\leq   C\nor{(\L + 1)^p u}{L^2(\M)}$ and apply~\eqref{e:operator-momotone-functions} to $A= \Lambda^{\frac{2p}{k}}$, $B=(\L + 1)^p$, and $\alpha = \frac{s}{p}\in [0,1]$ to obtain the result.
\enp

\subsection{Subelliptic estimates for manifolds with boundaries}
\label{subsectsubellbound}
In this section, we assume that $\M$ is a compact manifold with a nonempty boundary $\d \M$, and write $\M = \Int(\M) \cup \d \M$, with a disjoint union.
We assume that the coefficients of $X_j$'s are smooth up to the boundary, and that $\vect(X_1 , \cdots , X_m)(x) = T_x \M$ for $x \in \M \setminus K$, where $K$ is a compact subset of $\Int(\M)$ (i.e. the operator $\L$ is elliptic in the neighborhood of the boundary) and that Assumption~\ref{assumLiek} is satisfied on $K$.
\begin{theorem}
\label{cor:HsHsLbord}
Denote by $\Delta_D$ the Laplace-Dirichlet operator associated to some/any Riemannian metric equal to that issued from the vector fields $(X_1 , \cdots , X_m)$ in $ \M \setminus K$.
Then, for all $s \geq 0$, we have $ \H_\L^{s} := D(\L^{\frac{s}{2}}) \subset D\left((-\Delta_D)^{\frac{s}{2k}}\right)$. Moreover, there exists $C>0$ such that for all $u \in \H_\L^{s}$, we have
\bna
\nor{\left(-\Delta_D\right)^{\frac{s}{2k}} u}{L^2(\M)}\leq C\nor{u}{\H_\L^{s}} .
\ena
\end{theorem}
Note that the space $D\left((-\Delta_D)^{s}\right)$ does not depend on the metric chosen to define $\Delta$ inside $\Int(\M)$ but only on its values in a neighborhood of $\d\M$.
Remark that we also have the converse simple inclusion $D\left((-\Delta_D)^{s}\right) \subset D(\L^{s})$.

\bigskip

We now explain how the estimates in the previous section have to be modified to yield the statement of Theorem~\ref{cor:HsHsLbord}. 

We first let $\Lambda^{\frac{1}{k}}$ (with a slight abuse of notation: $\Lambda^{\frac{1}{k}}$ is not the $\frac{1}{k}$ power of an operator) be 
\begin{itemize}
\item a pseudodifferential operator of order $\frac{1}{k}$ on $\M$, with kernel compactly supported in $\M \times \M$,
\item formally selfadjoint on $L^2(\M)$,
\item locally elliptic on a neighborhood $\calN$ of $K$ (i.e. $K \Subset \Int(\calN) \Subset \Int(\M)$),
\item with kernel compactly supported in an $\eps$-neighborhood of $\diag(\M \times \M)$, with $\eps \ll \dist(K , \calN^c)$. This implies that for $n\in \N$, $\left(\Lambda^{\frac{1}{k}}\right)^n$ has a kernel compactly supported in an $ n\eps$-neighborhood of $\diag(\M \times \M)$. Here, we will use the abuse of notation $\Lambda^{\frac{n}{k}}$ instead of $(\Lambda^{\frac{1}{k}})^n$. This will not lead to any confusion since we will only use  $\Lambda^{s}$ for $s$ of the form $\frac{n}{k}$ with $n\in \N$.
\end{itemize}

When a maximal Sobolev exponent $s_0$ is fixed, we will need to use such operators for $n\leq n_0$ with $n_0=s_0k$ and make proofs by induction using several cutoff functions. At each step, we shall need to make some estimates on some $C\e$ neighborhood of the zone where we get the information, with $C$ depending on $s_0$ and some geometric properties of the cutoff functions. At the end, once the number of steps is fixed, we can select $\e$ small enough (and the associated $\Lambda^{\frac{1}{k}}$) so that all the reasoning is valid. To make the presentation more readable, we have chosen not to keep track of all the constants and the geometrical conditions involved. Yet, the proof will make it clear that there is $\e_0>0$ depending on $s_0$, $\calN$, $K$ and $\M$ so that all the support conditions of the following proof are fulfilled if $0<\e<\e_0$.

That Assumption~\ref{assumLiek} is satisfied on $K$ (and hence on $\calN$, since $\L$ is elliptic on $\M \setminus K$) yields, for all $\chi \in C^\infty_0(\M)$ such that $\chi=1$ in a neighborhood of $K$, the existence of $C>0$ such that for all $u \in C^\infty_0(\M)$, we have (see again \cite{RS:76} Theorem~17 and estimate (17.20) p311)
\bna
\nor{\Lambda^\frac{1}{k} \chi u}{L^2(\M)}^2\leq C\sum_{i=1}^m\nor{X_i u}{L^2(\M)}^2+C\nor{u}{L^2(\M)}^2 , 
\ena
and hence, still for  $u \in C^\infty_0(\M)$,
\bnan
\label{estimhypo2-bis}
\nor{\Lambda^\frac{1}{k} \chi u}{L^2(\M)}^2\leq C(\L u , u)_{L^2(\M)} +C\nor{u}{L^2(\M)}^2 .
\enan
We now decompose the proof in several lemmata. 

Several times in the proof, we shall use the following fact of pseudodifferential calculus.
Given $n\leq n_0$ and a function $\varphi\in C^{\infty}_0(\calN)$, we remark that $\Lambda^{\frac{n}{k}}$ is elliptic of order $\frac{n}{k}$ in a neighborhood of $\supp(\varphi)$. As a consequence, the classical parametrix construction (see for instance~\cite[Proof of  Theorem~18.1.9]{Hoermander:V3}) allows, for any $N\in \N$, to construct a pseudodifferential operator $\tilde{\Lambda}^{-\frac{n}{k}}_N$ of order $-\frac{n}{k}$, elliptic on a neighborhood of $\supp(\varphi)$, such that $\tilde{\Lambda}^{-\frac{n}{k}}_N \Lambda^{\frac{n}{k}}=\varphi(x)+R^1_N$ and $\Lambda^{\frac{n}{k}}\tilde{\Lambda}^{-\frac{n}{k}}_N=\varphi(x)+R^2_N$ with $R^i_N$, $i=1,2$, pseudodifferential operators of order $-N$ with kernel compactly supported in $\M\times \M$. In the applications, $N$ will always be fixed, sufficiently large.

\begin{lemma}
\label{:l:hypo-bord-1}
For all $\chi_0 \in C^\infty_0(\calN)$ such that $\chi_0=1$ in a neighborhood of $K$, and all $\chi_1 \in C^\infty_0(\calN)$ such that $\chi_1=1$ in a neighborhood of $\supp(\chi_0)$, for all $s\in \N/k$ with $s+2/k\leq s_0$, for $\e$ small enough, there is $C>0$ such that for all $u \in C^\infty_0(\M)$, we have
\bna
\nor{\Lambda^{s+\frac{2}{k}} \chi_0 u}{L^2(\M)}^2\leq C\| \Lambda^s \chi_0 \L u\|_{L^2(\M)}^2 +C\nor{\Lambda^{s+\frac{1}{k}}\chi_1 u}{L^2(\M)}^2 +\nor{\chi_1 u}{L^2(\M)}^2
\ena
\end{lemma}
\bnp
The proof is almost the same as that of Estimate~\eqref{estimhypocor-2} in Corollary~\ref{corhypoestim}, and only relies on the application of~\eqref{estimhypo2-bis} (instead of~\eqref{estimhypo2}). 
First, estimate~\eqref{estimhypo2-bis} applies for $\chi u$ replaced by $\Lambda^{s+\frac{1}{k}} \chi_0 u$ since $\Lambda^{s+\frac{1}{k}} \chi_0 u = \tilde{\chi}_0\Lambda^{s+\frac{1}{k}} \chi_0 u$ for $\tilde{\chi}_0 = 1$ in an $n\eps$-neighborhood of $\supp (\chi_0)$, with $n=sk+1$. This yields 
\bna
\nor{\Lambda^{s+\frac{2}{k}} \chi_0 u}{L^2(\M)}^2\leq C\left(\L \Lambda^{s+\frac{1}{k}} \chi_0 u , \Lambda^{s+\frac{1}{k}} \chi_0 u\right)_{L^2(\M)} +C\nor{\Lambda^{s+\frac{1}{k}} \chi_0 u}{L^2(\M)}^2 .
\ena
Then, a computation similar to \eqref{computcommut}, and the only difference comes from the estimate of the remainder term 
$$
\Re(\chi_0 \Lambda^{s+\frac{1}{k}}[\L,\Lambda^{s+\frac{1}{k}}\chi_0 ]u,u)_{L^2(\M)} = \Re((T_1 +T_2)u,u)_{L^2(\M)} =  \Re(T_2 u,u)_{L^2(\M)} ,
$$ 
where $T_2$ is a pseudodifferential operator of order ${2s+\frac{2}{k}}$, with kernel supported in $\supp(\chi_0) \times \supp(\chi_0)$. Given $\varphi\in C^{\infty}_0(\calN)$ such that $\varphi=1$ on $\supp(\chi_1)$, we may define the associated parametrix $\tilde{\Lambda}^{-(s+\frac{1}{k})}$ of $\Lambda^{s+\frac{1}{k}}$ as above. Writing $T_2=\chi_1 \varphi T_2  \varphi\chi_1 $, we now have $T_2 = \chi_1 \Lambda^{s+\frac{1}{k}} \tilde{\Lambda}^{-(s+\frac{1}{k})*} T_2 \tilde{\Lambda}^{-(s+\frac{1}{k})} \Lambda^{s+\frac{1}{k}} \chi_1+ R$, where $R$ is a smoothing operator with kernel compactly supported in $\M \times \M$. The boundedness of $\tilde{\Lambda}^{-(s+\frac{1}{k})*} T_2 \tilde{\Lambda}^{-(s+\frac{1}{k})}$ (as a pseudodifferential operator of order zero) and $R$ on $L^2(\M)$ then implies 
$$
| \Re(T_2 u,u)_{L^2(\M)}|
= | \Re(T_2 \chi_1u , \chi_1 u)_{L^2(\M)}| \leq C\nor{\Lambda^{s+\frac{1}{k}} \chi_1  u}{L^2(\M)}^2+C\nor{\chi_1 u}{L^2(\M)}^2,
$$
 which concludes the proof. Note also that the term $\nor{\Lambda^{s+\frac{1}{k}} \chi_0 u}{L^2(\M)}^2$ has been bounded by $\nor{\Lambda^{s+\frac{1}{k}} \chi_1 u}{L^2(\M)}^2$ the same way using that $ \Lambda^{s+\frac{1}{k}} \chi_0 =\Lambda^{s+\frac{1}{k}} \chi_0\varphi \chi_1 =\Lambda^{s+\frac{1}{k}} \chi_0 \tilde{\Lambda}_1^{-(s+\frac{1}{k})}\Lambda^{s+\frac{1}{k}}\chi_1 +R\chi_1$ with $R$ smoothing.
\enp

Before going further, recall that we can localize~\eqref{estimhypo2-bis} under the following form.
\begin{lemma}
\label{:l:hypo-bord-1bis}
For all $\chi_0 \in C^\infty_0(\calN)$ such that $\chi_0=1$ in a neighborhood of $K$, and all $\chi_1 \in C^\infty_0(\calN)$ such that $\chi_1=1$ in a neighborhood of $\supp(\chi_0)$, there is $C>0$ such that for all $u \in C^\infty_0(\M)$, we have
\bna
\nor{\Lambda^{\frac{1}{k}}\chi_0 u}{L^2(\M)}^2 \leq C\left(\chi_0 \L u , \chi_0 u \right)_{L^2(\M)}  + C\nor{\chi_1u}{L^2(\M)}^2 
\ena
\end{lemma}
\bnp
We apply Estimate~\eqref{estimhypo2-bis} with $\chi$ such that $\chi=1$ on a neighborhood of $\supp(\chi_0)$:
\bna
\nor{\Lambda^\frac{1}{k} \chi_0 u}{L^2(\M)}^2 & = & \nor{\Lambda^\frac{1}{k}\chi  \chi_0 u}{L^2(\M)}^2 
\leq  C(\L \chi_0 u , \chi_0 u)_{L^2(\M)} +C\nor{\chi_0 u}{L^2(\M)}^2 \\
& \leq & C( \chi_0 \L u , \chi_0 u)_{L^2(\M)} + \Re ([\L, \chi_0] u , \chi_0 u)_{L^2(\M)} +C\nor{\chi_0 u}{L^2(\M)}^2 ,
\ena
where $[\L, \chi_0]$ is a skew-adjoint first order differential operator: the principal part of $ \chi_0 [\L, \chi_0]$ is thus skew-adjoint, and hence 
$$
\Re ([\L, \chi_0] u , \chi_0 u)_{L^2(\M)} =\Re (\chi_0[\L, \chi_0]\chi_1 u , \chi_1 u)_{L^2(\M)}\leq  C \nor{\chi_1 u}{L^2(\M)}^2 ,
$$
which, together with the preceding estimate, proves the lemma.
\enp

\begin{lemma}
\label{:l:hypo-bord-2}
For all $\chi_0, \chi_1, \chi_2 \in C^\infty_0(\calN)$ such that $\chi_0=1$ in a neighborhood of $K$, $\chi_1=1$ in a neighborhood of $\supp(\chi_0)$, and $\chi_2=1$ in a neighborhood of $\supp(\chi_1)$, for all $p \in \N$ with $\frac{p}{k}+\frac{2}{k}\leq s_0$ and $\e$ small enough, there is $C>0$ such that for all $u \in C^\infty_0(\M)$, we have
\bna
\nor{\Lambda^{\frac{p}{k}+\frac{2}{k}} \chi_0 u}{L^2(\M)}^2\leq C\nor{ \Lambda^{\frac{p}{k}} \chi_1\L u}{L^2(\M)}^2  + C\nor{\chi_2 u}{L^2(\M)}^2 .
\ena
\end{lemma}
\bnp
We prove the statement by induction. The case $p=0$ follows directly from the estimate of Lemma~\ref{:l:hypo-bord-1} with $s=0$, combined with Lemma~\ref{:l:hypo-bord-1bis} (using an additional cutoff function as done below).
Assume now this is true for $p$, then, the estimate of Lemma~\ref{:l:hypo-bord-1}  with $s=\frac{p+1}{k}$ yields, for some $\tilde{\chi}_1$ such that $\tilde{\chi}_1 = 1$ on a neighborhood of $\supp(\chi_0)$ and $\chi_1 = 1$ in a neighborhood of $\supp(\tilde{\chi}_1)$,
\bna
\nor{\Lambda^{\frac{p+1}{k}+\frac{2}{k}} \chi_0 u}{L^2(\M)}^2
& \leq & C\| \Lambda^{\frac{p+1}{k}}\chi_0 \L u\|_{L^2(\M)}^2 +C\nor{\Lambda^{\frac{p+1}{k}+\frac{1}{k}}\tilde{\chi}_1u}{L^2(\M)}^2+C\nor{\tilde{\chi}_1u}{L^2(\M)}^2 .
\ena
Using then the induction assumption for $p$ for the term $\Lambda^{\frac{p}{k}+\frac{2}{k}} \tilde{\chi}_1u$ gives, since $\chi_1 = 1$ in a neighborhood of $\supp(\tilde{\chi}_1)$,
\bnan
\label{e:technical-bdr}
\nor{\Lambda^{\frac{p+1}{k}+\frac{2}{k}} \chi_0 u}{L^2(\M)}^2 
& \leq & C\| \Lambda^{\frac{p+1}{k}}\chi_0 \L u\|_{L^2(\M)}^2 +C\| \Lambda^{\frac{p}{k}} \chi_1 \L u\|_{L^2(\M)}^2 +C\nor{\chi_1 u}{L^2(\M)}^2 .
\enan
We now use pseudodifferential calculus and the parametrices of $\Lambda^{\frac{p+1}{k}}$ and $\Lambda^{\frac{p}{k}}$to write, for $\varphi = 1$ on $\supp(\chi_1)$,
\bna
&&\Lambda^{\frac{p+1}{k}}\chi_0 =\Lambda^{\frac{p+1}{k}} \chi_0 \varphi \chi_1=\Lambda^{\frac{p+1}{k}} \chi_0 \tilde{\Lambda}^{-\frac{p+1}{k}}\Lambda^{\frac{p+1}{k}} \chi_1+R \chi_1   , \\
&& \Lambda^{\frac{p}{k}} \chi_1=\Lambda^{\frac{p}{k}} \varphi \chi_1 =\Lambda^{\frac{p}{k}}  \tilde{\Lambda}^{-\frac{p+1}{k}}\Lambda^{\frac{p+1}{k}}\chi_1 +R \chi_1 ,
\ena
and hence
\bna
&& \| \Lambda^{\frac{p+1}{k}}\chi_0 \L u\|_{L^2(\M)} \leq C \| \Lambda^{\frac{p+1}{k}} \chi_1\L u\|_{L^2(\M)} +C \| R \chi_1 \L u\|_{L^2(\M)} 
 \leq C \| \Lambda^{\frac{p+1}{k}} \chi_1\L u\|_{L^2(\M)} +C \| \chi_2 u\|_{L^2(\M)} , \\
&& \| \Lambda^{\frac{p}{k}} \chi_1 \L u\|_{L^2(\M)} \leq C \| \Lambda^{\frac{p+1}{k}} \chi_1\L u\|_{L^2(\M)} +C \| \chi_2 u\|_{L^2(\M)} ,
\ena
which, combined with~\eqref{e:technical-bdr} concludes the proof of the statement for $p+1$, and hence of the lemma.
\enp

\begin{lemma}
\label{:l:hypo-bord-3} 
For all $\chi_0 \in C^\infty_0(\calN)$ such that $\chi_0=1$ in a neighborhood of $K$, and all $\chi_1 \in C^\infty_0(\M)$ such that $\chi_1=1$ in a neighborhood of $\supp(\chi_0)$, for all $p \in \N$ and $\e$ small enough, there is $C>0$ such that for all $u \in C^\infty_0(\M)$, we have
\bna
\nor{\Lambda^{\frac{2p}{k}}\chi_0 u}{L^2(\M)}^2
\leq C \sum_{j=0}^p \nor{ \chi_1 \L^j u}{L^2(\M)}^2 .
\ena
\end{lemma}
\bnp
Again, we prove this by an induction argument. For $p=1$, this is the estimate of Lemma~\ref{:l:hypo-bord-2} with $p=0$.

Assume the result for $p$. Using Lemma~\ref{:l:hypo-bord-2}, we obtain, for some $\tilde{\chi}_1, \tilde{\chi}_2$ such that $\tilde{\chi}_1 = 1$ on a neighborhood of $\supp(\chi_0)$, $\tilde{\chi}_2 = 1$ on a neighborhood of $\supp(\tilde{\chi}_1)$ and $\chi_1 = 1$ in a neighborhood of $\supp(\tilde{\chi}_2)$,
\bna
\nor{\Lambda^{\frac{2p+2}{k}} \chi_0 u}{L^2(\M)}^2\leq
 C\nor{\Lambda^{\frac{2p}{k}} \tilde{\chi}_1\L u}{L^2(\M)}^2  + C\nor{\tilde{\chi}_2u}{L^2(\M)}^2 .
\ena

Applying the induction assumption with $u$ replaced by $\L u$ (and $\chi_0$ replaced by $\tilde{\chi}_1$) yields (since $\chi_1 = 1$ in a neighborhood of $\supp(\tilde{\chi}_1)$),
\bna
\nor{\Lambda^{\frac{2p+2}{k}}\chi_0 u}{L^2(\M)}^2 
\leq C \sum_{j=0}^p \nor{ \chi_1 \L^j (\L u)}{L^2(\M)}^2 
 + C\nor{\tilde{\chi}_2 u}{L^2(\M)}^2 ,
\ena
which concludes the proof of the lemma since $\chi_1 = 1$ in a neighborhood of $\supp(\tilde{\chi}_2)$.
\enp

Combining Lemma~\ref{:l:hypo-bord-3} together with classical ellipticity at the boundary, we are now ready for proving the following results.
\begin{proposition}
\label{prop:HkHkLbord}
For all $m \in \N$, there is $C>0$ such that for all $u \in D(\L^{mk})$, we have
\bna
\nor{u}{H^{2m}(\M)}^2 \leq C \sum_{j=0}^{mk} \| \L^j u\|_{L^2(\M)}^2 \leq C \| (\L+1)^{mk} u\|_{L^2(\M)}  .
\ena
\end{proposition}
From this proposition, we directly obtain by interpolation the statement of Theorem~\ref{cor:HsHsLbord} (see e.g. the proof of Corollary~\ref{cor:HsHsL}) using that $\L^j u|_{\d \M} = \Delta u|_{\d \M}$ for all $u \in C^\infty(\M)$.

\bnp[Proof of Proposition~\ref{prop:HkHkLbord}]
First write Lemma~\ref{:l:hypo-bord-3} with $p=mk$, $\chi_0\in C^{\infty}_0(\calN)$ with $\chi_0 =1$ in a neighborhood $\widetilde{\calN}\subset \calN$ of $K$, yielding
\bnan
\label{e:internal-subelliptic}
 \nor{u}{H^{2m}(\widetilde{\calN})}^2
  \leq \nor{\Lambda^{2m}\chi_0 u}{L^2(\M)}^2
\leq C \sum_{j=0}^{mk} \nor{ \chi_1 \L^j u}{L^2(\M)}^2  .
\enan
Then, concerning estimates near the boundary, we first have the following statement. For all $\theta_0 \in C^\infty(\M)$ such that $\theta_0=1$ in a neighborhood of $\d\M$, $\supp(\theta_0) \cap K = \emptyset$, and all $\theta_1 \in C^\infty(\M)$ such that $\theta_1=1$ in a neighborhood of $\supp(\theta_0)$, since $\L$ elliptic in $\supp(\theta_1)$, there is $C>0$ such that for all $u \in C^\infty(\M)$ with $u|_{\d \M}=0$, we have
\bnan
\label{estimelliptic}
\nor{\theta_0 u}{H^{m+2}(\M)}^2 \leq C \|\theta_1 \L u\|_{H^m(\M)}^2  + C\nor{\theta_1u}{L^2(\M)}^2  .
\enan
This is actually a corollary of the usual proof of elliptic regularity up to the boundary. Yet, to check it directly, we can apply the global elliptic regularity result (see \cite[Theorem 5 p323]{Evans:98}) to $\theta_0 u$ with a global elliptic operator $\widetilde{\L}$ equal to $\L$ on $\supp(\theta_1)$. Let $V$ be an open subset with $\supp (\nabla \theta_1)\Subset V\Subset U\Subset \M$, where $U$ in an open set so that $\theta_1=1$ on $U$. This yields
\bna
\nor{\theta_0 u}{H^{m+2}(\M)}&\leq &C \|\theta_0 \widetilde{\L}u\|_{H^m(\M)} +C \|u\|_{H^{m+1}(V)} + C\nor{\theta_0 u}{L^2(\M)}\\
 &\leq &C \|\theta_1 \L u\|_{H^m(\M)} + C\nor{\theta_1 u}{L^2(\M)}.
\ena
where we have used interior regularity (see \cite[Theorem 2 p314]{Evans:98}) that gives $\|u\|_{H^{m+1}(V)}\leq C \|\widetilde{\L} u\|_{H^m(U)}+C \nor{u}{L^2(U)}\leq C \|\theta_1 \L u\|_{H^m(\M)}+C \nor{\theta_1u}{L^2(\M)}$. This proves \eqref{estimelliptic}.

Then, from~\eqref{estimelliptic}, another induction argument as in the proof of Lemma~\ref{:l:hypo-bord-3} gives, for all $m\in \N$ and $\theta_0 , \theta_1$ as above, the existence of $C>0$ such that for all $u \in C^\infty(\M)$ with $u|_{\d \M}=  \L u |_{\d \M}=\cdots =  \L^{m-1} u |_{\d \M} = 0$,
\bna
\nor{\theta_0 u}{H^{2m}(\M)}^2
\leq C \sum_{j=0}^m \nor{ \theta_1 \L^j u}{L^2(\M)}^2 .
\ena
Combining this for a function $\theta_0$ equal to $1$ on a neighborhood of $\M \setminus \widetilde{\calN}$ together with~\eqref{e:internal-subelliptic} now implies for all $u \in C^\infty(\M)$ such that $\L^j u |_{\d \M} = 0$, for $0\leq j\leq m-1$, 
\bna
 \nor{u}{H^{2m}(\M)}^2
\leq C \sum_{j=0}^{mk} \nor{ \chi_1 \L^j u}{L^2(\M)}^2 +C \sum_{j=0}^m \nor{ \theta_1 \L^j u}{L^2(\M)}^2 .
 \ena
Since the set of such functions $u$ is dense in $D(\L^{mk})$, this yields the sought result.
\enp

\section{Sub-Riemannian norm of normal vectors}
\label{app:davide}
\begin{lemma}
\label{l:davide}
Let $X_i \in \R^d$ for $i=1, \cdots, m$ and, for $v \in \vect(X_i , i=1, \cdots, m)$, set
$$
g(v) = \inf\left\{\sum_{i=1}^m u_i^2 , \quad (u_1,\cdots,u_m)\in \R^m,\sum_{i=1}^mu_i X_i =v\right\} .
$$
Then, for any $\xi \in (\R^d)^*$, and for $v_0 = \sum_{i=1}^m 2 \langle \xi , X_i \rangle X_i$, we have $g(v_0) = 4 \ell(\xi)$ where $\ell(\xi) =\sum_{i=1}^m \langle \xi , X_i \rangle^2$.
\end{lemma}
Note that this is clear if the family $(X_i)_{i=1, \cdots, m}$ is linearly independant.
\bnp
We want to compute the minimum 
$$
g(v_0) = \inf\left\{\sum_{i=1}^m u_i^2 , \quad (u_1,\cdots,u_m)\in \R^m,\sum_{i=1}^mu_i X_i = \sum_{i=1}^m 2 \langle \xi , X_i \rangle X_i\right\} .
$$
First that taking $u_i = 2 \langle \xi , X_i \rangle$ in this definition direcly yields that 
\bnan
\label{ineq-gl}
g(v_0) \leq \sum_{i=1}^m (2 \langle \xi , X_i \rangle )^2 = 4 \ell(\xi) .
\enan
Then it only remains to prove the converse inequality.
To this aim, remark that $$ \ell(\xi) = \max \left\{ \sum_{i=1}^m u_i \langle \xi , X_i \rangle- \frac14 \sum_{i=1}^m u_i^2 , \quad (u_1,\cdots,u_m)\in \R^m \right\}.$$ As a consequence, we have
$$
\ell(\xi) \geq \Big< \xi , \sum_{i=1}^m u_i X_i \Big>- \frac14 \sum_{i=1}^m u_i^2 , \quad \text{ for all }(u_1,\cdots,u_m)\in \R^m .
$$
Hence, for all $(u_1,\cdots,u_m)\in \R^m$ such that $\sum_{i=1}^mu_i X_i = \sum_{i=1}^m 2 \langle \xi , X_i \rangle X_i$, we obtain 
$$
 \ell(\xi) \geq \Big< \xi ,  \sum_{i=1}^m 2 \langle \xi , X_i \rangle X_i \Big> - \frac14 \sum_{i=1}^m u_i^2
  = \sum_{i=1}^m 2 \langle \xi , X_i \rangle^2 - \frac14 \sum_{i=1}^m u_i^2 
= 2 \ell(\xi)- \frac14 \sum_{i=1}^m u_i^2,
$$
that is $4 \ell(\xi) \leq \sum_{i=1}^m u_i^2$, and hence $4\ell(\xi) \leq g(v_0)$. This, together with~\eqref{ineq-gl} concludes the proof of the lemma.
\enp

\small
\bibliographystyle{alpha}
\bibliography{bibli}
\end{document}